\documentclass[11pt,reqno]{amsart}
\usepackage{graphicx} 
\usepackage{amsmath,amssymb,amsthm}
\usepackage{url}
\usepackage{lmodern}
\usepackage{color}
\usepackage{tikz,tikz-cd}
\usepackage{dsfont}
\usepackage{comment}
\usepackage{cleveref}
\usepackage{enumitem}
\usepackage{thmtools}
\usepackage{mathtools}
\usepackage{thm-restate}
\usepackage{mathrsfs}
\usepackage{appendix}
\usepackage{xcolor}
\usepackage{array}

\newcommand{\AppMatrix}[1]{%
  \begingroup
  \setlength{\arraycolsep}{3pt}%
  \setbox0=\hbox{$\displaystyle #1$}%
  \ifdim\wd0>\linewidth
    \resizebox{\linewidth}{!}{\copy0}%
  \else
    \box0
  \fi
  \endgroup}
\calclayout

\newtheorem{theorem}{Theorem}[section]
\newtheorem{lemma}{Lemma}[section]
\newtheorem{proposition}{Proposition}[section]

 \theoremstyle{remark}
 \newtheorem{remark}[theorem]{Remark}

\newtheorem{prop}{Proposition}[section]
\newtheorem{lem}{Lemma}[section]
\newtheorem{defi}{Definition}[section]

\newcommand{\R}{\mathbb{R}}

\newcommand{\Lf}{\textup{\textbf{L}}}

\newcommand{\Af}{\textup{\textbf{A}}}

\newcommand{\Qf}{\textup{\textbf{Q}}}
\renewcommand{\Re}{\textup{Re}}
\renewcommand{\Im}{\textup{Im}}

\newcommand{\C}{\mathbb{C}}
\theoremstyle{definition}

\newcommand{\tr}{\mathrm{tr}}
\newcommand{\sov}{\mathfrak{so}(5)}

\newcommand{\ditto}{---\raisebox{-0.5ex}{\texttt{"}}---}

\title[Mode stability of a self-similar Yang-Mills blow-up]{On the mode stability of a self-similar blow-up of Yang-Mills heat flow in five dimensions}
\author{Josef E. Greilhuber}
\address{MIT, Department of Mathematics, 77 Massachusetts Avenue, Cambridge, MA 02142, USA}
\email{jgreil@mit.edu}
\author{David Wallauch}
\address{EPFL SB MATH PDE, Batiment MA, Station 8, CH-1015 Lausanne}
\email{david.wallauch@epfl.ch}
\date{April 2026}

\begin{document}
\begin{abstract}
 We study the spectral stability properties of an explicitly known self-similar blowup profile to the Yang-Mills heat flow on $\R^5\times SO(5)$ often referred to as the Weinkove solution. We prove that this solution is mode stable in a dynamically evolving DeTurck gauge. To arrive at that result, we first reduce the problem to countably many vector-valued Schrödinger operators on the half line. We then proceed by proving uniform estimates for all but finitely many of said operators. The spectra of the remaining operators are studied individually by means of GGMT bounds or by comparison arguments. Two of these operators are distinguished as they have negative eigenvalues originating from symmetries. Among these two there is a scalar operator, previously analyzed by Donninger and Sch\"orkhuber, and an operator acting on a three-dimensional vector bundle. The latter operator's matrix structure prohibits a straightforward supersymmetric removal. To circumvent this, we transform the equation to a Schrödinger operator with drift terms to simplify the structure of the eigenfunction. We then develop tools to study such operators and conclude that said operator has no second negative eigenvalue.
\end{abstract}
\maketitle

\section{Introduction}
This article is concerned with Yang-Mills theory, a class of partial differential equations used to describe the behavior of elementary particles. The Yang--Mills equations are posed on a principal $G$-bundle, where $G$ is in general a compact Lie group. In this article, we focus on the special case where the underlying Lie group is $SO(5)$. To introduce the corresponding Yang-Mills equations properly, we  start by considering Lie algebra-valued 1-forms $ A=A_j dx^j$, with $A_j: \R^5\to \mathfrak{so}(5)$. Throughout this work, we view $\mathfrak{so}(5)$ as the Lie algebra consisting of skew-symmetric $5\times 5$ matrices over $\R$, endowed with the standard commutator bracket $[\cdot,\cdot].$  Then, for 1-forms with values in $\mathfrak{so}(5)$, the 
covariant derivative with respect to $A$ in the adjoint representation
is given by
\[D^j_A \Phi = \partial^j \Phi + [A^j,\Phi].\]
Furthermore, the curvature tensor is of the form
\begin{align*}
    F_{ij}^A=\partial_i A_j-\partial_j A_i +[A_i,A_j].
\end{align*}
With this at hand, we introduce the Yang-Mills functional\footnote{Throughout this work we employ Einstein's summation convention.}
\begin{equation}\label{Eq: functional}
    \mathrm{YM}(A)=-\frac{1}{4}\int_{\R^5}  \text{tr}(F_{ij}^A F^{ij}_A) dx
\end{equation}
and note that the Euler-Lagrange equations associated to \eqref{Eq: functional} are given by
\begin{align}\label{Eq: static equation}
D^i_A F_{ij}^A=0.
\end{align}
One approach to constructing solutions to the static equations \eqref{Eq: static equation} is to utilize the gradient flow 
\begin{align}\label{Eq: heat flow}
    \partial_0 A_j = D^i_A F_{ij}^A,
\end{align}
called Yang-Mills heat flow. This system of equations has been studied intensively by the mathematical community in several general settings, see for instance the monograph \cite{Fee16} for an extensive overview. For general well-posedness results, we refer the reader to \cite{Rad92} for global well-posedness in 2 and 3 spatial dimensions. The four dimensional case is more involved and local well-posedness was established by Struwe in \cite{Str94}. A first global well-posedness result is due to  \cite{SchStrTah98} for equivariant connections. The full global problem was resolved much later by Waldron \cite{Wal19}. However, even though solutions exist globally in $4$ dimensions, blowup at infinity is possible \cite{SirWerZheZho26}. Starting from dimension $d\geq 5$, the scaling invariance 
$$
A \mapsto A_\lambda:= \lambda A(\lambda^2 t, \lambda x),\qquad \lambda>0
$$
identifies Eq.~\eqref{Eq: heat flow} as energy-supercritical. Thus, even though local well-posedness is classical, see Remark 2.29 in \cite{OhTat22}, one no longer expects global existence of general solutions. Instead, finite time blowup is expected and was established in various works cf. \cite{Nai94,Gro01,Gas02,Wan25}. One way this phenomenon can occur is through homothetically shrinking solitons. These are solutions of the form 
$$
U(t,x)=\frac{1}{\sqrt{T-t}}W_j\left(\frac{x-x_0}{\sqrt{T-t}}\right) dx^j.
$$
An explicit example of such a solution was found by Weinkove \cite{Wei04} and is given by
\begin{align}\label{Eq; Weinkove}
W^T(t,x)=\frac{1}{\sqrt{T-t}} W_i(\frac{x}{\sqrt{T-t}}) dx^i
\end{align}
with
$$
W_i(\xi)= \frac{\sigma_i(\xi)}{\frac{\sqrt{3}}{2\sqrt{2}}|\xi|^2+9-7\sqrt{3/2}},\qquad (\sigma_i)_{\nu}^\mu(\xi)=\delta_{i\nu}\xi^\mu-\delta_{i\mu}\xi^\nu.
$$
Note that $W^T$ belongs to the previously mentioned subclass of equivariant solutions, i.e. it is of the form
\[
A_j(t,x)=u(t,|x|)\,\sigma_j(x).
\]
Interestingly, due to \cite{BizWas15}, equivariant self-similar shrinkers are known not to exist in dimensions $d\geq 10$. 
\\
In addition to the scaling invariance and the translation symmetries of Eq.~\eqref{Eq: heat flow} another important property of the Yang-Mills heat flow is its gauge invariance.
Namely, if $A=A_i dx^i$ is a solution, then, for any $g:\R^5 \to SO(5)$, so is the form $A_\mu^g dx^\mu$ given by
$$
A_i^g = g^{-1}(\partial_i g) + g^{-1} A_i g.
$$
In contrast to the other symmetries of the Yang-Mills heat flow, applying a gauge transformation to a connection does not yield a distinct solution, instead it is just a different geometrical description of the same object. Consequently, when studying  Eq.~\eqref{Eq: heat flow} one commonly fixes a gauge. A standard example of this is to impose the Coulomb gauge
$$
\partial^j A_j =0.
$$
However, for our purposes a different choice of gauge, called the DeTurck gauge, is better suited. To introduce it, we first note the dependency on time motivates us to introduce an auxiliary variable $A_0(t,x) dt$ 
 with corresponding curvature term
$$
F_{0j}=\partial_0 A_j-\partial_j A_0+[A_0,A_j]= \partial_0 A_j -D_j^A A_0.
$$
This allows us to recast Eq.~\eqref{Eq: heat flow} in the gauge-covariant form
\begin{equation*}
    \partial_0 A_j =D^i F_{ij} +D_j A_0. 
\end{equation*}
Now, we have, to some extent, transferred the freedom of gauge to the freedom of choosing $A_0$. For instance, common choices are $A_0=0$, which corresponds to the temporal gauge and gives the gradient flow \eqref{Eq: heat flow} or 
$$
A_0=\partial^i A_i
$$
which yields the standard Yang-Mills--DeTurck flow
\begin{equation}\label{Eq: deT flow2}
    \partial_0 A_j =D^i_A F_{ij} +D_j^A(\partial^i A_i). 
\end{equation}
Since we are interested in properties of $W$ we choose a modified version of $A_0$. Namely, we let $\widetilde W_i(t,x)=\frac{1}{\sqrt{T-t}} W_i(\frac{x}{\sqrt{T-t}})$ and
set
\[
A_0=\partial^iA_i+[\widetilde W^i,A_i]
\]
which yields the equation 

\begin{align}\label{Eq: deT flow}
    \partial_0 A_j=D^i_A F_{ij} +D_j^A(\partial^i A_i +[\widetilde W^i,A_i])
\end{align}
which we call the Yang-Mills--DeTurck flow (with background $W$). To avoid confusion, we specify at this point that whenever we refer to the Yang-Mills--DeTurck flow, we will always mean \eqref{Eq: deT flow}. We also want to emphasize that this adapted choice of gauge does not impose any constraints on the dynamics of the Yang-Mills flow compared to another choice of gauge and in particular, on intervals of smooth existence, \eqref{Eq: deT flow} is equivalent, through a time-dependent gauge transformation, to the Yang–Mills heat flow in temporal gauge. Instead, it is just a useful choice as it renders the flow linearised around $W$ parabolic. 
We  also note that the Weinkove solution lies in the Coulomb gauge, i.e.~ satisfies $\partial^i \widetilde W_i=0$, hence it also solves Eq.~\eqref{Eq: deT flow} which expanded reads as
\begin{align*}
    \partial_t A_j&=D^i_A (\partial_i A_j-\partial_j A_i +[A_i,A_j])+D_j^A( \partial^i A_i+[\widetilde W^i,A_i])
    \\
    &=\Delta A_j+2[A^i, \partial_i A_j]-[A^i, \partial_j A_i]+[A^i,[A_i,A_j]]+\partial_j[\widetilde W^i,A_i]+[A_j,[\widetilde W^i,A_i]].
\end{align*}

\subsection{Mode stability}
One of the key components in the stability analysis of self-similar solutions is mode stability. To introduce mode stability properly, we begin by considering the similarity variables 
\[
\tau =-\log(T-t),\qquad \xi=\frac{x}{\sqrt{T-t}}
\]
and, in accordance with the scaling symmetry, define $B_j(\tau,\xi)=e^{-\frac{\tau}{2}} A_j(t,x)$.
Then, as
\[
\partial_t=e^{\tau} (\partial_\tau +\frac{1}{2}\xi^i\partial_i),\qquad \partial_{x_i}=e^{\frac{\tau}{2}}\partial_{\xi_i}
\]
 the equations  \eqref{Eq: deT flow} transform into
\begin{equation}\label{Eq:free flow in simvar}
\begin{split}
    \partial_0 B_j= \Delta B_j &-\frac12\Lambda B_j-\frac{1}{2}B_j+2[B^i, \partial_i B_j]
    \\
    &-[B^i, \partial_j B_i]+[B^i,[B_i,B_j]]+\partial_j[W^i,B_i]+[B_j,[W^i,B_i]]
    \end{split}
\end{equation}
with 
$\Lambda f(x)= x\cdot \nabla f(x) $. 
The idea behind mode stability is to characterise spectral instabilities of the generator of the semigroup of the linearised flow. Thus, we linearise \eqref{Eq:free flow in simvar} around 
$
W_i
$
to (formally) obtain 
\[
\partial_0 U_j =(\Lf U)_j
\]
with 
\begin{equation}\label{Eq: linearised op}
\begin{split}
        (\Lf U)_j &:=  \Delta U_j -\frac12(\Lambda+1) U_j+2[W^i, \partial_iU_j]+[W^i,[W_i,U_j]]+2[U^i,F_{ij}^W]
   \end{split}
\end{equation}
and insert the mode  ansatz $U(\tau,\xi)= e^{\lambda \tau} \Phi(\xi)$ into this equation.
With this ansatz, we arrive at the spectral equation
\begin{align*}
    (\lambda-\Lf)\Phi&=0.
\end{align*}
We aim to characterise all \(\lambda\in\mathbb C\) with \(\Re\lambda\geq0\) for which this equation admits a nonzero admissible profile \(\Phi\). Admissibility naturally depends on the underlying space.
To also account for possible oscillatory modes, we allow the spectral parameter \(\lambda\) to be complex and consider $\mathbf L$ on $1$-forms taking values in the complexified Lie algebra
\[
\mathfrak{so}(5)_{\mathbb C}
:=\mathfrak{so}(5)\otimes_{\mathbb R}\mathbb C
\cong\{A\in\mathbb C^{5\times5}:A^T=-A\}.
\]
To ease notation, we suppress the subscript \(\mathbb C\) and write \(\mathfrak{so}(5)\) for its complexification.
Next, we define a weight $\rho(x)=e^{-\frac{|x|^2}{4}}$ as well as an associated $L^2$ inner product
\[
(U,V)_{L^2_\rho}
:=\int_{\mathbb R^5}\operatorname{tr}\!\bigl(U^i(x)V_i(x)^*\bigr)\,\rho(x)\,dx.
\]
where \(A^*=\overline A^{\,T}\). We define \(H^2_\rho\) by taking the closure of compactly supported smooth \(\mathfrak{so}(5)_{\mathbb C}\)-valued one-forms with respect to the corresponding norm.
\begin{align*}
  \|U\|_{H^2_\rho(\R^5,\Omega(\mathfrak{so}(5)))}^2:=  (U,U)_{L^2_\rho(\R^5)}+(\partial^jU,\partial_j U)_{L^2_\rho(\R^5)}+(\partial^i \partial^jU,\partial_i \partial_jU)_{L^2_\rho(\R^5)}.
\end{align*}
Then, thanks to Lemma \ref{lem: self-adjoint}, we know
$$
\left(\Lf, H^2_\rho(\R^5,\Omega(\mathfrak{so}(5)))\right)
$$
is a self-adjoint operator on $L^2_\rho(\R^5)$, consequently, 
$$
D(\Lf):= H^2_\rho(\R^5,\Omega(\mathfrak{so}(5)))
$$
is the natural underlying space for our spectral analysis. 
Note that the translation symmetries of the Weinkove solution yield unstable modes, 
namely, if we define 
\[
(h_{1})_j=(\Lambda+1)W_j,
\qquad \text{ and } \qquad 
\big(h_{\frac{1}{2}}^{(i)}\big)_j=F_{ij}^W
=\partial_iW_j-D_j^W W_i,
\qquad i=1,\dots,5,
\]
then
\[
\Lf h_1=h_1,
\qquad
\Lf h_{\frac12}^{(i)}=\frac12 h_{\frac12}^{(i)}.
\]
Here $h_1$ arises from variation of the blowup time, whereas $h_{\frac12}^{(i)}$
is the infinitesimal gauge-corrected translation in the $i$-th direction.
This allows us to finally state our main theorem.

\begin{theorem}
\label{thm:main}
Consider the operator \[-\Lf: H^2_\rho(\R^5,\Omega(\mathfrak{so}(5)))\subset L^2_\rho(\R^5,\Omega(\mathfrak{so}(5)))\to L^2_\rho(\R^5,\Omega(\mathfrak{so}(5)))\]
where
\begin{align*}
      (\Lf U)_j &:=  \Delta U_j -\frac12(\Lambda+1) U_j+2[W^i, \partial_iU_j]+[W^i,[W_i,U_j]]+2[U^i,F_{ij}^W].
\end{align*}
Then the spectrum of the operator $(-\Lf, H^2_\rho(\R^5,\Omega(\mathfrak{so}(5))))$ is discrete and
    $$
\sigma (-\Lf)\cap \{z\in \C: \Re z \leq 0 \}=\{-1\}\cup \{-\frac12\}. 
    $$
Furthermore, the multiplicity of the eigenvalue $-1$ is equal to $1$ and the multiplicity of the eigenvalue $-\frac{1}{2}$ is $5$. 
\end{theorem}

We make the following two remarks.
\begin{remark}(\emph{On the gauge invariance of the eigenvalues})
Part of the spectrum of $\mathbf L$ arises from infinitesimal gauge transformations, applied to the stationary solution $W$. An infinitesimal gauge transform is represented by an $\sov$-valued $0$-form $s$, that is, a field of antisymmetric $5 \times 5$ matrices on $\mathbb R^5$. The action of $s$ on a $1$-form $\Psi$ is given by $D_\Psi s = d s + [\Psi,s]$. It is not hard to verify that
\begin{align*}
    \mathbf L D s = D (D^\ast D - \frac12 \Lambda) s.
\end{align*}
The operator $D^\ast D - \frac12 \Lambda$ is nonpositive, showing that gauge modes do not give rise to instabilities in the linearized Yang--Mills--DeTurck flow, and thus could in principle be ignored in our analysis. Nevertheless, isolating the eigenvalues arising from such infinitesimal gauge transformations is impractical, and this has not been attempted here.
\end{remark}
\begin{remark}(\emph{On the generalization to higher dimensions})
The Weinkove solution exists not just in dimension $5$, but for all $5\leq d\leq 9$. The symmetry reduction and subsequent analysis pursued here carry over to any dimension, similarly yielding a countable sequence of Schrödinger operators acting on fairly low-dimensional vector bundles on the half-line. However, since the structure theory of the Lie Algebra $\mathfrak{so}(n)$ depends strongly on the dimension, so do the specifics of the symmetry reduction.
\end{remark}

\Cref{thm:main} is an assertion on the spectrum of a self-adjoint linear operator. In the context of the nonlinear Yang--Mills--DeTurck flow, only real-valued initial data is admissible. For the spectral analysis, it makes no difference whether we work with  $\sov$-valued $1$-forms or their complexification, that is, $1$-forms valued in $\mathfrak{so}(5)_{\mathbb C} = \mathfrak{so}(5) \oplus i \, \mathfrak{so}(5)$. Since the latter is preferable from a Lie-theoretic point of view, for the remainder of this paper, all $\sov$-representations, operators and domains are implicitly complexified.

\subsection{Related results}
The study of self-similar solutions to parabolic PDEs and, more concretely their stability properties, is an extremely active area of research. We therefore mention only a handful of results on other closely related geometric models, such as the harmonic map heat flow and the Ricci flow. For the Ricci flow, there is an analogue of the Weinkove solution,
called the Feldman--Ilmanen--Knopf shrinker, whose linear stability without symmetry assumptions was recently established in \cite{NafOzu25}, following earlier nonlinear universality results in the $U(2)$-invariant K\"ahler setting \cite{Max14,BinJia16}.\\ 
In contrast to the Ricci flow, it was shown by Fan \cite{Fan19} that for $3\leq d\leq 6 $ infinitely many shrinkers exist for the corotational harmonic maps heat flow with target manifold $\mathbb S^d$, followed by numerical investigations in \cite{BieBiz11}. A first rigorous construction of a spectrally stable shrinker in 3 dimensions by means of a computer assisted proof followed \cite{BieDon18} with subsequent nonlinear stability \cite{BieDonSch17}. These results were later extended to dimensions $4,5$ and $6$ in \cite{AngKisSch26} with 6 space dimensions being the largest dimension for which corotational shrinkers can exist \cite{BizWas15}. For target manifolds other than $\mathbb
S^d$ stable shrinkers have been constructed in \cite{GloKisSch24}.
\\
Returning to the Yang-Mills flow, there are 3 important articles that are very closely tied to this work. The first of these is by Donninger and Schörkhuber \cite{DonSch19}. In their paper, the authors establish full nonlinear asymptotic stability of the Weinkove solution under the assumption of equivariant symmetry in 5 dimensions. Naturally, this also entails the mode stability of the Weinkove solution within that symmetry class. Later, their result was generalised to $5\leq d \leq 9$ by Glogić and Schörkhuber \cite{GloSch20}. More recently, a finite co-dimensional nonlinear stability result without any symmetry assumptions was shown by Kim and Lee \cite{KimLee24}. However, they do not address the mode stability of the Weinkove solution, hence only obtain a conditional result.
\\
Aside from these,  variational stability of self-similar shrinkers to the parabolic Yang-Mills equations was studied in \cite{CheZha15} and \cite{KelStr16}, while the structure of singular sets was investigated by Kelleher and Streets \cite{KelStr18}.
Lastly, we also want to note that non-self-similar blowup solutions to the Yang-Mills flow have been constructed in \cite{BenDuoGho25} for dimensions $d\geq 11$.

\subsection{Outline of this work}

In \Cref{sec:symmetry_reduction}, we first recall the necessary facts from the representation theory of compact Lie groups, in particular $SO(5)$, and then discuss the symmetry reduction of the spectral problem \eqref{Eq: linearised op}.
The section concludes with a discussion of the accompanying symbolic calculations, which are performed using the Python package SymPy and contained in user-friendly Jupyter notebooks. We obtain a sequence of vector-valued Schr\"odinger operators on the positive real axis indexed by \emph{weights} $\lambda \in \mathbb N^2$.
A related approach to reduce the mode stability analysis of a self-similar solution to a PDE was recently implemented for the wave maps equation, where the authors applied Lie algebra methods to the spaces of homogeneous polynomials to completely decouple the problem into scalar-valued ODEs \cite{WeiKocDon25,DonMos26}. 
In Section~\ref{section:general schr}, we prove some general statements on vector-valued Schrödinger operators. In particular, we prove a vector-valued version of the GGMT bounds and show how one can estimate the spectral gap between the first and second eigenvalues of vector-valued Schrödinger operators with drift. The usefulness of this stems from the fact that, in contrast to the scalar-valued case, one cannot simply perform a supersymmetric removal to get rid of a ground state eigenvalue. One way to mitigate this is to conjugate the operator in such a way that the eigenfunction takes the form $(f,0,\dots,0)^T$, which allows one to perform a supersymmetric removal in the upper-left component. However, such a transformation generally introduces first-order derivatives.

We move on to Section~\ref{section: low angulars}, where we study the finitely many low angular operators that cannot be controlled by straightforward bounds on the lower order terms of the operator $\mathbf L$. Here we apply the previously devised machinery as well as min-max principles and simple positivity arguments to show that all low angular components, aside from two exceptions, have only positive spectrum. One of the remaining two operators corresponds to the equivariant spectral equation and has already been studied in~\cite{GloSch20}. The other one is a $3\times 3$ operator with negative eigenvalue $-\frac12$. By using the tools established in Section~\ref{section:general schr}, we prove that this eigenvalue is simple and that the operator has no other negative spectral point. 

In \Cref{sec:high_modes}, we wrap up the argument by providing a proof that all modes not already handled in \Cref{section: low angulars} do not give rise to negative eigenvalues. This follows from pointwise bounds on the lower order terms in $\Lf$, achieved by a combination of elementary representation theory and finite-dimensional symbolic calculations.

\subsection{Accompanying files}

Computations underpinning the present article are contained in five accompanying files. The two Jupyter notebooks and the two Mathematica files are meant to be read in tandem with the article.

\begin{enumerate}
    \item \texttt{SO5\_basis\_vectors.ipynb}, a Jupyter notebook containing the highest weight vectors spanning the weight space decompositions used in the symmetry reduction of the operator $\mathbf L$ given in \eqref{Eq: linearised op}.
    \item \texttt{SO5\_operators.ipynb} computes the coefficients of the vector--valued Schrödinger operators representing $\mathbf L$ on these weight spaces.
    \item \texttt{SO5\_utility.py} is an auxiliary file containing definitions and methods used by the preceding two notebooks.
    \item \texttt{computations\_low\_momenta.nb} contains symbolic and numerical calculations accompanying \Cref{section: low angulars}.
    \item \texttt{computations\_high\_momenta.nb} accompanies \Cref{sec:high_modes}.
\end{enumerate}

In both Mathematica files, numerical computations are limited to straightforward checks that certain constants are positive. Furthermore, as is elaborated in Section $4$, these checks can be made completely rigorous with elementary but involved computations. 

\subsection{AI statement}

This article is based almost entirely on original ideas of the authors. We acknowledge the contribution of one technical idea in the analysis of a mode of low angular momentum by GPT 5.6. GPT 6.2 was used as an editing assistant and to prepare the tikz figures in this article.

The code accompanying this project is entirely human written. We attempted to capture the feeling of a (machine-aided) pen-and-paper calculation and aimed for instant user-friendliness rather than a high capacity for automatization, which in the present day lies elsewhere entirely. 

\subsection{Acknowledgements}

The first author was partially supported by NSF grant DMS-2247185, awarded to Eugenia Malinnikova, and by the David and Lucile Packard Foundation through the Packard Science and Engineering Fellowship awarded to Alexandr Logunov. He is sincerely grateful to Eugenia Malinnikova and Alexandr Logunov for their support and many valuable discussions and insights. Much of this work was completed while the first author was visiting ETH Zürich and MIT in Fall 2025 and Winter 2026. He is grateful to these institutions and his hosts for their hospitality. The second author was partially supported by by the Austrian Science Fund (FWF) through a Erwin-Schrödinger fellowship (J-4875). We would also like to thank Prof.~Irfan Glogić for his insights regarding this problem.

\section{Symmetry reduction}
\label{sec:symmetry_reduction}

In this section, we discuss the symmetry reduction which transforms the eigenvalue problem \eqref{Eq: linearised op} into a countable list of eigenvalue problems for vector-valued Schr\"odinger operators on the half line. Underlying the symmetry reduction of this problem is the group action $\rho$ of the group $SO(5)$ on the space of $L^2$-sections of $T^\ast \mathbb R^5 \otimes \sov$ on which the eigenvalue problem under consideration is defined.

For a $\sov$-valued $1$-form $\Psi$, $V \in T_x \mathbb R^5$ and $g \in SO(5)$, this action is given in elementary terms by
\begin{align*}
    \langle \rho(g) \Psi (x),  V \rangle = g \cdot \langle \Psi(g^{-1} x), g^{-1} V \rangle \cdot g^{-1}
\end{align*}
This action descends to an action of $\sov$ on $\Gamma(T^\ast \mathbb R^5 \otimes \sov)$ given by 
\begin{align*}
    \rho(X) \Psi (x)_j = [X,\Psi_j] - X^k_j \Psi_k - X^k_\ell x^\ell \partial_k \Psi_j.
\end{align*}

\subsection{Basic structure theory of $\sov$ and its representations}

We begin with a review of the necessary facts from the representation theory of compact Lie groups in general and the Lie group $SO(5)$ in particular. All are contained in any textbook on the representation of compact Lie groups; we direct the reader to the especially concise and well-written book by Sepanski \cite{Sepanski2007}.

The Lie algebra $\sov$ is $10$-dimensional and possesses a Cartan subalgebra $\mathfrak c$, that is, a maximal subalgebra on which the Lie bracket is trivial, spanned by the two matrices
\begin{align*}
    T_1 = \begin{pmatrix}
        0 & 0 & 0 & 0 & 0 \\
        0 & 0 & -1 & 0 & 0 \\
        0 & 1 & 0 & 0 & 0 \\
        0 & 0 & 0 & 0 & 0 \\
        0 & 0 & 0 & 0 & 0 
    \end{pmatrix}, \ 
    T_2 = \begin{pmatrix}
        0 & 0 & 0 & 0 & 0 \\
        0 & 0 & 0 & 0 & 0 \\
        0 & 0 & 0 & 0 & 0 \\
        0 & 0 & 0 & 0 & -1 \\
        0 & 0 & 0 & 1 & 0 
    \end{pmatrix}
\end{align*}
The choice of Cartan subalgebra is not unique, but after it is made, one obtains a special basis of $\sov$ -- strictly speaking, of its complexification $\mathbb C \otimes \sov$ -- by simultaneous diagonalization of the action of $T_1$ and $T_2$ on $\sov$. Besides $T_1$ and $T_2$, on which $\mathfrak c$ acts trivially, one obtains eight more simultaneous eigenvectors, $L_{(\pm 1,0)}$, $L_{(0,\pm 1)}$ and $L_{(\pm 1,\pm 1)}$. Together, these eigenvectors form the so-called \emph{root basis} of $\sov$ with respect to $\mathfrak c$. Except for $(0,0)$, the eigenvalue pairs, that is, the \emph{roots} of $\sov$, are purely imaginary and appear with multiplicity one. Thus, the remaining elements of the root basis besides $T_1$ and $T_2$ are conveniently labeled by the imaginary part of the eigenvalue pair, so that $[T_j,L_\alpha] = i \alpha_j L_\alpha$, $\alpha \in \{(\pm 1,0),(0,\pm 1),(\pm 1,\pm 1)\}$. To avoid confusion with the Laplace operator, we denote the set of roots by $\Lambda(\sov)$.

\subsubsection{Positive roots and the Weyl group}
The set of roots of any Lie group $G$ admits an action of a discrete group, the \emph{Weyl group}, generated by the reflections across the orthogonal complement of each root. The collection of hyperplanes $\alpha^\bot$, $\alpha \in \Lambda$, splits $\mathbb R^{\dim \mathfrak c}$ into convex cones called \emph{Weyl chambers}, each bounded by the orthogonal complements of a finite collection of \emph{simple roots}. Any root can be written as a linear combination of these with exclusively nonnegative or nonpositive coefficients, leading to a decomposition of the set of roots into \emph{positive} and \emph{negative} roots.

In our case, we choose the Weyl chamber
\begin{align*}
    \{\alpha: \alpha_1 \geq \alpha_2 \geq 0\}.
\end{align*}
The corresponding simple roots are $(0,1)$ and $(1,-1)$, and the set of positive roots is thus
\begin{align*}
    \Lambda_+ = \{(0,1),(1,1),(1,0),(1,-1)\}.
\end{align*}

\subsubsection{Representations and highest weights}
Given a representation $\rho: \sov \to \mathfrak{gl}(V)$, we may likewise simultaneously diagonalize the matrices $\rho(T_1),\rho(T_2)$, obtaining an orthogonal decomposition $V = \bigoplus_{\mu} V_\mu$ into \emph{weight spaces}, such that 
\begin{enumerate}
    \item $\rho(T_j) v_\mu = i \mu_j \, v_\mu$ for all $v_\mu \in V_\mu$,
    \item $\rho(L_\alpha) v_\mu \in V_{\mu + \alpha}$ for all $v_\mu \in V_\mu$.
\end{enumerate}
If $\rho$ is the tangent representation of a group representation of $SO(5)$, all weights $\mu$ lie in $\mathbb Z^2$. A \emph{highest weight} vector $v \in V_\mu$ for the representation $(V,\rho)$ is one such that $L_\alpha v = 0$ for all positive roots $\alpha \in \Lambda_+$. The weight $\mu$ is then called a highest weight of $(V,\rho)$.
\begin{enumerate}
  \setcounter{enumi}{2}
    \item Each irreducible representation has exactly one highest weight $\mu$, and $\dim V_\mu = 1$.
    \item Up to isomorphism, the set of irreducible representations of $\sov$ is in one-to-one correspondence with highest weights $\mu \in \mathbb N^2$, $\mu_1 \geq \mu_2 \geq 0$.
    \item If and only if $V$ is irreducible,
    \begin{align*}
        V = \mathrm{span}\left(\left\{\rho(L_{- \alpha_k}) \ldots \rho(L_{- \alpha_1}) v, \, \alpha_1,\ldots,\alpha_k \in \Lambda_+, k \in \mathbb N\right\}\right)
    \end{align*}
    for any nonzero highest weight vector $v$.
\end{enumerate}
In this article, we will denote a general irreducible representation of highest weight $\lambda$ by $V_\lambda$, as well as sometimes specific ones if context makes it clear which representation is meant. Note that $\sov$ is itself an irreducible $\sov$-representation with highest weight $(1,1)$.

\subsubsection{Tensor products of representations}
Given representations $(V^{(1)},\rho^{(1)})$ and $(V^{(2)},\rho^{(2)})$, the tensor product $V^{(1)}\otimes V^{(2)}$ is equipped with the action given by
\begin{align*}
    \rho(X) ( v^{(1)} \otimes v^{(2)} ) = (\rho^{(1)}(X) v^{(1)}) \otimes v^{(2)} + v^{(1)} \otimes (\rho^{(2)}(X) v^{(2)}).
\end{align*}
It follows that the tensor product $V^{(1)}_{\lambda_1} \otimes V^{(2)}_{\lambda_2}$ of the weight spaces $V^{(1)}_{\lambda_1}$ and $V^{(2)}_{\lambda_2}$ lies in the weight space $(V^{(1)} \otimes V^{(2)})_{\lambda_1 + \lambda_2}$. In particular, for irreducible $V^{(1)}$ and $V^{(2)}$ with highest weights $\mu_1$ and $\mu_2$, $(V^{(1)} \otimes V^{(2)})_{\mu_1 + \mu_2} = V^{(1)}_{\mu_1} \otimes V^{(2)}_{\mu_2}$ is one-dimensional and a highest weight space of the -- not usually irreducible -- representation $V^{(1)} \otimes V^{(2)}$.

\subsubsection{Decomposition of a tensor product into irreducibles}
\label{sec:layer_stripping}
All representations occurring in this article are equipped with an $SO(5)$-invariant inner product. We employ this inner product to decompose the tensor product of representations $V^{(1)}$ and $V^{(2)}$ with highest weights $\mu_1$ and $\mu_2$ into irreducible components as follows.

Suppose $V^{(1)} \otimes V^{(2)} = \bigoplus_{k=1}^m E_{\lambda_k}$, where $\lambda_1,\ldots,\lambda_m$ are (not necessarily distinct) weights ordered by nonincreasing height. Then $\lambda_1 = \mu_1 + \mu_2$, and given any weight $\lambda' \le \lambda_1$, one can recover the weight space
\begin{align*}
    \left(E_{\lambda_1}\right)_{\lambda'} = \mathrm{span}\big(\big\{ &\rho(L_{ - \alpha_1}) \ldots \rho(L_{ - \alpha_\ell}) (v^{(1)}_{\mu_1} \otimes v^{(2)}_{\mu_2}) : \\ &\alpha_1,\ldots,\alpha_\ell \in \Lambda_+, \alpha_1 + \ldots + \alpha_\ell = \lambda_1 - \lambda' \big\}\big).
\end{align*}
It is helpful to consider the index set in this computation as the set of \emph{paths} connecting $\lambda_1$ to $\lambda'$ via negative roots. This set is always finite, since $\Lambda_-$ is contained in a strictly convex cone.
Note that $[\rho(L_{- \alpha}), \rho(L_{- \alpha'})] = c_{\alpha,\alpha'} \rho(L_{- \alpha - \alpha'})$ if $\alpha + \alpha'$ is another positive root and zero otherwise. To obtain a set spanning $\left(E_{\lambda_1}\right)_{\lambda'}$, it thus suffices to choose an arbitrary ordering of the positive roots and only consider paths from $\lambda_1$ to $\lambda'$ where all roots occur in the given order, with possible repetitions, preventing combinatorial explosion of the number of paths.

The next highest weight, $\lambda_2$, is simply the first weight such that $(E_{\lambda_1})_{\lambda_2}$ has nontrivial orthogonal complement in $(V^{(1)} \otimes V^{(2)})_{\lambda_2}$. We choose a basis vector $v_{\lambda_2} \in (E_{\lambda_1})_{\lambda_2}^\bot$, the new highest weight vector, and repeat the procedure outlined above until all irreducible components $E_{\lambda_k}$ have been obtained.

The procedure terminates when all weight spaces $(V^{(1)} \otimes V^{(2)})_{\mu_1 + \mu'}$, with $\mu'$ among the weights of $V^{(2)}$, have been filled, as the following lemma shows.

\begin{lemma}
    \label{lem:tensor_weigths}
    All highest weights $\mu$ occurring in the decomposition of $V^{(1)} \otimes V^{(2)}$ into irreducible components lie in 
    \begin{align*}
        \mu_2 + \{\mu: V^{(1)}_{\mu} \neq \{0\} \},
    \end{align*}
    that is, they can be written as the sum of the highest weight of $V^{(2)}$ and a weight of $V^{(1)}$.
\end{lemma}

\begin{proof}
    Write $W^{(1)},W^{(2)} \subseteq \mathbb Z^2$ for the sets of weights with nontrivial weight spaces $V^{(1)}$ and $V^{(2)}$, respectively. For each weight space $V^{(1)}_\mu$, pick a basis $(v_\mu^{j})_{j=1}^{m(\mu)}$, where $m(\mu) = \dim V^{(1)}_\mu$. A vector $v \in (V^{(1)} \otimes V^{(2)})_\lambda$ then decomposes as
    \begin{align*}
        v = \sum_{\mu \in W^{(1)}} \sum_{j=1}^{m(\mu)} v_{\mu}^j \otimes w_{\lambda - \mu}^j,
    \end{align*}
    where $w_{\lambda-\mu}^j$ are uniquely determined vectors in $V^{(2)}_{\lambda-\mu}$, not necessarily forming a basis.

    Assume $v$ is a highest weight vector of $V^{(1)} \otimes V^{(2)}$. Then for all positive roots $\alpha$, that is, for all $\alpha \in \{(0,1),(1,1),(1,0),(1,-1)\}$,
    \begin{align*}
        L_\alpha v = \sum_{\mu \in W^{(1)}} \sum_{j=1}^{m(\mu)} \left( \big(L_\alpha v_{\mu}^j\big) \otimes w_{\lambda - \mu}^j + v_{\mu}^j \otimes \big(L_\alpha w_{\lambda - \mu}^j\big) \right) = 0.
    \end{align*}
    Now pick a lowest weight $\mu$ such that $w_{\lambda-\mu}^j \neq 0$ for at least one $j$. Keeping in mind that $L_\alpha v_\mu^j \in V_{\mu+\alpha}^{(1)}$ and $L_\alpha w_{\lambda - \mu}^j \in V_{\lambda - \mu + \alpha}^{(2)}$, and collecting terms in $V^{(1)}_\mu \otimes V^{(2)}$, we get
    \begin{align*}
        0 = \sum_{j=1}^{m(\mu-\alpha)} \big(L_\alpha v_{\mu-\alpha}^j\big) \otimes w_{\lambda - \mu + \alpha }^j + \sum_{j=1}^{m(\mu)} v_{\mu}^j \otimes \big(L_\alpha w_{\lambda - \mu}^j\big).
    \end{align*}
    By the assumed minimality of $\mu$, the left sum vanishes, so $L_{\alpha} w^j_{\lambda - \mu} = 0$ for all $j$. Since this argument was independent of $\alpha$, $L_{\alpha} w^j_{\lambda - \mu} = 0$ for all $\alpha$. This means that either $w^j_{\lambda - \mu} = 0$ for all $j$, contradicting our assumptions, or that $\lambda - \mu$ is a highest weight of $V^{(2)}$.
\end{proof}

To verify that a given vector $v \in (V^{(1)} \otimes V^{(2)})_{\lambda}$ is a highest weight vector of an irreducible component $E_\lambda$, it suffices to check that $\rho(L_\alpha) v = 0$ for all $\alpha \in \Lambda_+$. This allows for an independent verification, implemented in the accompanying file \texttt{SO5\_basis\_vectors.ipynb}, that the computation of a highest weight vector in a tensor representation did not produce a faulty result.

\subsection{Decomposition of the space of sections}
Let $\mathcal H_\ell$ denote the space of homogeneous harmonic polynomials of degree $\ell$ on $\mathbb R^5$, restricted to $\mathbb S^{4}$. It is well known that the space $L^2(\mathbb S^4)$ admits the orthogonal decomposition 
\begin{align*}
    L^2(\mathbb S^4) = \bigoplus_{\ell=0}^\infty \mathcal H_\ell
\end{align*}
Each space $\mathcal H_\ell$ carries an irreducible $\sov$-action. As an $\sov$-representation, $\mathcal H_\ell \cong V_{(\ell,0)}$.

By the scaling $\int_{\R^5} dx = \int_0^\infty r^4 \int_{\mathbb S^4} d\omega dr$, we may decompose
\begin{align*}
    L^2(\mathbb R^5) = \bigoplus_{\ell=0}^\infty L^2((0,\infty),r^4) \otimes \mathcal H_\ell.
\end{align*}
For the space of square integrable sections of $T^\ast \R^5 \otimes \sov$, we obtain the decomposition
\begin{align*}
    L^2\Gamma(T^\ast \R^5 \otimes \sov) = \bigoplus_{\ell=0}^\infty L^2((0,\infty),r^4) \otimes \mathcal H_\ell \otimes \mathbb R^5 \otimes \sov,
\end{align*}
into representations of $\sov$, where the action on $L^2((0,\infty),r^4)$ is trivial, the action on $\mathcal H_\ell$ is by pull-back as before, the action on $\mathbb R^5$ is the tautological one and the action on $\sov$ is the adjoint action.

Each tensor representation $\mathcal H_\ell \otimes \mathbb R^5 \otimes \sov$ splits into a finite direct sum of irreducible representations. Denote the sum of all irreducible components of $\mathcal H_\ell \otimes \mathbb R^5 \otimes \sov$ corresponding to highest weight $\lambda$ by $E_\lambda^{(\ell)}$, so that
\begin{align*}
    \mathcal H_\ell \otimes \mathbb R^5 \otimes \sov = \bigoplus_{\lambda \in \mathbb N^2} E_\lambda^{(\ell)}.
\end{align*}
Denote $E_\lambda = \bigoplus_{\ell=0}^\infty E_\lambda^{(\ell)}$. As it turns out, each $\dim E_\lambda^{(\ell)} = 0$ for all but finitely many $\ell$, so the direct sum defining $E_\lambda$ is in fact finite. We obtain a decomposition into $\sov$-invariant subspaces,
\begin{align*}
    L^2\Gamma(T^\ast \R^5 \otimes \sov) \cong \bigoplus_{\lambda \in \mathbb N^2} L^2((0,\infty),r^4) \otimes E_{\lambda}.
\end{align*}
For the purpose of analyzing $SO(5)$-invariant differential operators, it is also convenient to note that 
\begin{align*}
    \Gamma(T^\ast \R^5 \otimes \sov) \cong \overline{\bigoplus_{\ell \in \mathbb N} \bigoplus_{\lambda \in \mathbb N^2} C^\infty_\ell([0,\infty)) \otimes E_{\lambda}^{(\ell)}},
\end{align*}
where $C^\infty_\ell([0,\infty)) = \{r^\ell h(r^2), h \in C^\infty([0,\infty))\}$ and the closure is taken with respect to the smooth topology on $\Gamma(T^\ast \R^5 \otimes \sov)$.

\subsection{Symmetry reduction of invariant operators}
Let $\mathbf T$ be an $SO(5)$-invariant differential operator of order $q$ acting on $\Gamma(T^\ast \R^5 \otimes \sov)$. The subspaces $C^\infty([0,\infty)) \otimes E_\lambda$ are invariant under $\mathbf T$ by Schur's lemma. We pick a (non-canonical) decomposition of $E_\lambda$ into irreducible components $E_{\lambda,1},\ldots,E_{\lambda,m} \cong V_\lambda$ and identify $E_{\lambda,j}$ with $E_{\lambda,k}$ by an intertwining map $\phi^{k}_j: E_{\lambda,k} \to E_{\lambda,j}$ (unique up to scalar multiples). To simplify notation, extend $\phi^{k}_j$ to $E_\lambda$ by setting $\phi^{k}_j = 0$ on $E_{\lambda,\ell}$, $\ell \neq k$. This lets us write
\begin{align*}
    \mathbf T|_{C^\infty((0,\infty)) \otimes E_\lambda} \cong \sum_{\ell=0}^q (A_\ell)^j_k \phi^{k}_j \partial_r^\ell
\end{align*}
for some coefficients $(A_\ell)^j_k \in C^\infty((0,\infty))$. To compute these coefficients, we decompose $E_\lambda$ into its weight subspaces $E_{\lambda,\mu}$, $\mu \in \mathbb Z^2$, defined by $\rho(T_j) v = i \mu_j v$ for $v \in E_{\lambda,\mu}$. In particular, $\dim E_{\lambda,\lambda} = m$, since the highest weight subspace in any irreducible component $E_{\lambda,j}$ is one-dimensional. Since $\mathbf T$ commutes with $\rho(T_1)$ and $\rho(T_2)$, the weight subspaces $E_{\lambda,\mu}$ are invariant under $\mathbf T$. Pick normalized elements $v_{\lambda,j} \in E_{\lambda,j}$ such that $\phi^k_j$ maps $v_{\lambda,k}$ to $v_{\lambda,j}$. Then
\begin{align*}
    \mathbf T (f^j v_{\lambda,j})(r) = \sum_{\ell=0}^q (A_\ell)^j_k (\partial_r^\ell f^k(r)) v_{\lambda,j}
\end{align*}
Thus, in the symmetry class corresponding to highest weight $\lambda$, the operator $\mathbf T$ reduces to a differential operator on $C^\infty([0,\infty),\mathbb R^m)$ with coefficients in $C^\infty_{\mathrm{loc}}([0,\infty),\mathbb R^{m \times m})$, where $m$ is the multiplicity of the representation $V_\lambda$ in $E_\lambda$.

\subsection{Symmetry reduction of the linearized Yang--Mills--DeTurck operator}
\label{sec:symmetry_reduction_specific}

The symmetry reduction of the operator \eqref{Eq: linearised op} proceeds along purely algebraic lines. The reason is that, up to a radial factor, the coordinate components of the Weinkove solution $W$ are homogeneous polynomials. We introduce the following notation:
\begin{align*}
    W(r)_j &= f(r^2) \sigma_j, & (\sigma_{j})^\mu_\nu &= \delta_{j\nu} x^\mu - \delta_{j \mu} x^\nu, & \\ f(s) &:= (a s + b)^{-1}, &
    a &:= \frac{\sqrt{3}}{2\sqrt{2}}, & b &:= 9 - 7 \sqrt{\tfrac{3}{2}}.
\end{align*}
With this notation, the linearized operator reads
\begin{align}
    \label{Eq: linearised op reduced}
    \begin{split}
\Lf (\Psi)_j &= \partial_r^2 \Psi_j + \frac{4}{r} \partial_r \Psi_j + \frac1{r^2}\Delta_{\mathbb S^4} \Psi_j + \tfrac1{16}(12 - r^2) \Psi_j \\ &+2 f(r^2) [\sigma_i, \partial^i \Psi_j]+ f(r^2)^2 [\sigma_i,[\sigma^i,\Psi_j]] + 2f(r^2)^2[\Psi_i,[\sigma^i,\sigma_j]] \\ &+ 2 f(r^2) [\Psi_i, \partial^i \sigma_j - \partial_j \sigma_i] + 4 f'(r^2) [\Psi_i, x_i \sigma_j - x_j \sigma_i].
    \end{split}
\end{align}
We denote the various parts of this expression as follows.
\begin{align}
\label{eq:T_operators}
\begin{aligned}
    \mathbf T_1 \Psi_j &= [\sigma_i, \partial^i \Psi_j], &
    \mathbf T_2 \Psi_j &= [\sigma_i,[\sigma^i,\Psi_j]], \\
    \mathbf T_3 \Psi_j &= [\Psi_i, \partial^i \sigma_j - \partial_j \sigma_i], &
    \mathbf T_4 \Psi_j &= [\Psi_i, x_i \sigma_j - x_j \sigma_i], &
    \mathbf T_5 \Psi_j &= [\Psi_i,[\sigma^i,\sigma_j]].
\end{aligned}
\end{align}

Fix a weight $\lambda \in \mathbb N^2_{\lambda_1 \geq \lambda_2}$. Collect the irreducibles with highest weight $\lambda$ in the spaces $\mathcal H_\ell \otimes V_{(1,1)} \otimes V_{(1,0)}$, and pick a highest weight vector $v_k$ in each of them.
Then an element of $C^\infty([0,\infty)) \otimes E_\lambda$ of weight $\lambda$ can be written as
\begin{align*}
    \Psi(r,\omega) = \sum_{k=1}^m g_k(r) v_k(\omega).
\end{align*}
Let us consider the action of the operator $\mathbf T_1$ on an element $g_k(r) v_k(\omega)$ in this decomposition.
\begin{align*}
    [\sigma_i,\partial_i g_k(r) v_k(\omega)] &= g_k'(r)r^{-1}[\sigma_i, x^i v_k(\omega)] + g_k(r) [\sigma_i, \partial_i v_k(\omega)] = g_k(r) [\sigma_i, \partial_i v_k(\omega)],
\end{align*}
the latter since $x^i \sigma_i = 0$. This means that $\mathbf T_1$, although a differential operator of order $1$ on $\sov \otimes T^\ast \mathbb R^5$, reduces to a multiplication operator in each irreducible component, not unlike the spherical part of the Laplacian. As $\mathbf T_1$ commutes with the action of $SO(5)$ on $\sov \otimes T^\ast \mathbb R^5$, Schur's lemma implies
\begin{align*}
    [\sigma_i, \partial_i v_j(\omega)] = (T_1)^k_j v_k(\omega)
\end{align*}
for some matrix $((T_1)^k_j)_{j,k=1}^m$. A priori, this matrix could depend on $r$, but since $v_k(\omega)$ is homogeneous of order $0$ and $\sigma$ is homogeneous of order $1$, $T^k_j$ is constant. Its coefficients are recovered by computing 
\begin{align}
    \label{eq: matrix coefficient}
    (T_1)^k_j = \|v_k\|_{L^2(\mathbb S^4)}^{-1} \|v_j\|_{L^2(\mathbb S^4)}^{-1} \left\langle v_k, [\sigma_i, \partial_i v_j] \right\rangle_{L^2(\mathbb S^4)}.
\end{align}
The following two observations are helpful for efficient symbolic calculation:
\begin{enumerate}
    \item The value of $(T_1)^k_j $ is unaffected by replacing $v_k(\omega)$ by $g_k(r) v_k(\omega)$ in the above calculation. By construction, the coefficients of $v_k(\omega)$ are homogeneous harmonic polynomials of order $\ell$, restricted to $\mathbb S^4$. The evaluation of \eqref{eq: matrix coefficient} simplifies if we choose $g_k(r) = r^\ell$, so that all functions occurring in \eqref{eq: matrix coefficient} are harmonic homogeneous polynomials.
    \item Evaluating the $L^2$-inner product of monomials over $\mathbb S^4$ can be done quickly via a folklore formula popularized by Folland \cite{Folland2001}. For $\alpha_1,\ldots,\alpha_n \in \mathbb N$, write $\beta_j = \frac{1}{2}(\alpha_j + 1)$. Then
    \begin{align*}
        \int_{\mathbb S^{n-1}} x_1^{\alpha_1} \ldots x_n^{\alpha_n} = \begin{cases}
            0 & \text{ if some } \alpha_j \text{ is odd, } \\
            \frac{2\Gamma(\beta_1) \ldots \Gamma(\beta_n)}{\Gamma(\beta_1 + \ldots + \beta_n)}& \text{ if all } \alpha_j \text{ are even. }
        \end{cases}
    \end{align*}
\end{enumerate}

The remaining operators $\mathbf T_2,\ldots,\mathbf T_5$ are immediately seen to be multiplication operators, and the coefficients of the corresponding matrices are computed analogously. Note that while $\mathbf T_1$ and $\mathbf T_3$ are homogeneous of order $0$, the operators $\mathbf T_2,\mathbf T_4$ and $\mathbf T_5$ are homogeneous of order $2$, hence their matrix coefficients come with a pre-factor of $r^2$. Finally, the angular Laplacian acts as the scalar $- \ell(\ell+3)$ on $E_\lambda^{(\ell)}$. Hence, the action of $- \Delta_{\mathbb S^4}$ on $E_\lambda$ reduces to that of a diagonal matrix $S_\lambda$ with $k^{th}$ entry $\ell(\ell+3)$, where $\ell$ is the degree of the basis vector $v_k$.

On $[0,\infty) \times E_\lambda$, \eqref{Eq: linearised op reduced} acts as the operator
\begin{align}
    \label{Eq: linearised op symmetry reduced}
    \begin{split}
    \mathbf L_\lambda = \partial_r^2 &+ \frac4r \partial_r - \frac{1}{r^2} S_\lambda + \tfrac1{16}(12 - r^2) + 2f(r^2) T_1 + f(r^2)^2 r^2 T_2 \\ &+ 2 f(r^2) T_3 + 4 f'(r^2) r^2 T_4 + 2 f(r^2)^2 r^2 T_5 \\
    \end{split} \\
    \label{eq:linearised op symmetry reduced II}
    & \ = \partial_r^2 + \frac4r \partial_r - \frac{1}{r^2} S_\lambda + \tfrac1{16}(12 - r^2) + (a r^2 + b)^{-2}(r^2 A_\lambda + B_\lambda),
\end{align}
where we collected coefficients of $r^2f(r^2)^2$ and $f(r^2)^2$, respectively, obtaining
\begin{align}
\label{eq:def A}
    A_\lambda &= 2 a T_1 + T_2 + 2 a T_3 - 4 a T_4 + 2 T_5 , \\
\label{eq:def B}
    B_\lambda &= 2 b T_1 + 2 b T_3.
\end{align}
The representation $E_\lambda$ furthermore decomposes into parts $E_\lambda^{\mathrm{odd}}$ and $E_\lambda^{\mathrm{even}}$, distinguished by their parity under reflection across the origin. By inspection, each part is invariant under $\mathbf T_1,\ldots,\mathbf T_5$. Our construction yields basis vectors $v_1,\ldots,v_m$ for $E_{\lambda,\lambda}$ with coefficients in $\mathcal H_\ell$, and therefore a basis adapted to this decomposition. The matrices $T_1,\ldots,T_5$, and hence $A_\lambda$ and $B_\lambda$, thus have a block diagonal structure. We denote the corresponding blocks by $A_\lambda^{\mathrm{odd}}, B_\lambda^{\mathrm{odd}}$ and $A_\lambda^{\mathrm{even}}, B_\lambda^{\mathrm{even}}$, respectively.

\subsubsection{Decomposition of tensor products}
\label{sec:tensor_decomposition}

\begin{figure}[t]
    \centering
\begin{tikzpicture}[x=1.4cm,y=1.4cm,font=\small]
  \definecolor{latticegray}{gray}{0.75}
  \path[use as bounding box] (-2.5,-0.5) rectangle (2.5,2.5);
  \clip (-2.5,-0.5) rectangle (2.5,2.5);

  \foreach \x in {-2,...,2}
    \draw[latticegray,line width=0.4pt,dash pattern=on 2pt off 2pt]
      (\x,-0.5) -- (\x,2.5);
  \foreach \y in {0,...,2}
    \draw[latticegray,line width=0.4pt,dash pattern=on 2pt off 2pt]
      (-2.5,\y) -- (2.5,\y);

  \foreach \x in {-2,0,2}
    \fill[latticegray] (\x,0) circle[radius=1.3pt];
  \foreach \x in {-2,0,2}
    \fill[latticegray] (\x,2) circle[radius=1.1pt];

  \draw[black,line width=0.6pt,-{Stealth[length=5pt,width=4pt]}]
    (-2.5,0) -- (2.5,0);

  \draw[black,fill=white,line width=0.6pt] (0,1) circle[radius=3.1pt];
  \foreach \p in {(-1,0),(1,0),(-2,1),(-1,1),(0,1),(1,1),(2,1),(-1,2),(0,2),(1,2)}
    \fill[black] \p circle[radius=2.2pt];

  \node[anchor=north,inner sep=0pt,outer sep=0pt] at (-1,-0.12) {$(\ell\!-\!1,0)$};
  \node[anchor=north,inner sep=0pt,outer sep=0pt] at (1,-0.12) {$(\ell\!+\!1,0)$};
  \node[anchor=south,inner sep=0pt,outer sep=0pt] at (-2,1.12) {$(\ell\!-\!2,1)$};
  \node[anchor=south,inner sep=0pt,outer sep=0pt] at (-1,1.12) {$(\ell\!-\!1,1)$};
  \node[anchor=south,inner sep=0pt,outer sep=0pt] at (0,1.12) {$(\ell,1)$};
  \node[anchor=south,inner sep=0pt,outer sep=0pt] at (1,1.12) {$(\ell\!+\!1,1)$};
  \node[anchor=south,inner sep=0pt,outer sep=0pt] at (2,1.12) {$(\ell\!+\!2,1)$};
  \node[anchor=south,inner sep=0pt,outer sep=0pt] at (-1,2.12) {$(\ell\!-\!1,2)$};
  \node[anchor=south,inner sep=0pt,outer sep=0pt] at (0,2.12) {$(\ell,2)$};
  \node[anchor=south,inner sep=0pt,outer sep=0pt] at (1,2.12) {$(\ell\!+\!1,2)$};
\end{tikzpicture}
    \caption{The highest weights occurring in $\mathbf C^{(\ell)}$ (with multiplicity).}
    \label{fig:C_d}
\end{figure}
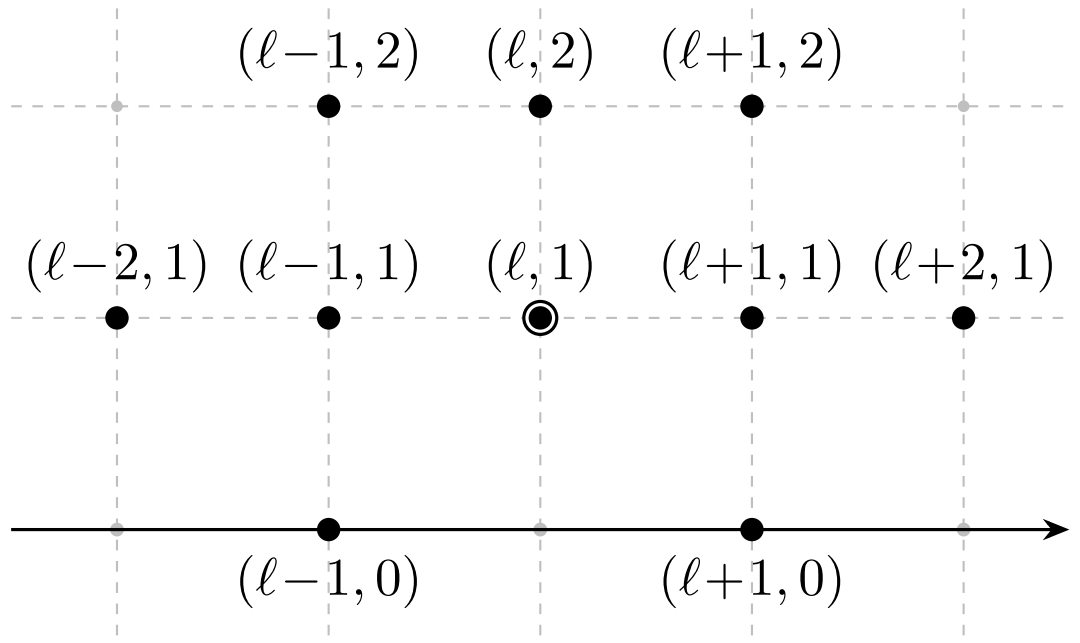

To compute the matrices $A_\lambda$ and $B_\lambda$ for specific weights $\lambda$, we have chosen to begin by decomposing the tensor representation $\sov \otimes T_0^\ast \mathbb R^5 \cong V_{(1,1)} \otimes V_{(1,0)}$, corresponding to forms with constant coefficients. It decomposes into three irreducibles $\mathbf A^{(0)}$, $\mathbf B^{(0)}$ and $\mathbf C^{(0)}$, respectively isomorphic to $V_{(1,0)}$, $V_{(1,1)}$ and $V_{(2,1)}$. We then decompose the tensor product
\begin{align*}
    \mathbf A^{(\ell)} = \mathbf A^{(0)} \otimes \mathcal H_{\ell} \cong V_{(1,0)} \otimes V_{(\ell,0)}
\end{align*}
into irreducible components $\mathbf A^{(\ell)}_\lambda$, and proceed analogously with $\mathbf B^{(\ell)}$ and $\ C^{(\ell)}$. We have employed the Mathematica package \texttt{LieART 2.0} \cite{LieART} to obtain a list of the highest weights occurring among the irreducible components in $\mathbf A^{(\ell)}$, $\mathbf B^{(\ell)}$ and $\mathbf C^{(\ell)}$. This is in principle not necessary, as \Cref{lem:tensor_weigths} shows the computation terminates quickly in the relevant settings. We note that the algorithm implemented in \texttt{LieART 2.0} is based on Klimyk's formula \cite[Section 5.4]{LieART}, which can be easily evaluated by hand as well. The following table lists all occurring irreducibles in $\mathbf A^{(\ell)}$, $\mathbf B^{(\ell)}$ and $\mathbf C^{(\ell)}$, respectively -- in degrees $0,1,2$, the decomposition has a sporadic structure, while in degree $3$ and above, it follows a set pattern, displayed in Figures \ref{fig:A_d}, \ref{fig:B_d} and \ref{fig:C_d}.

\vspace{3mm}
\begin{center}
\begin{tabular}{ |c|c|c|c| } 
\hline
$\ell$ & $\mathbf A^{(\ell)}$ & $\mathbf B^{(\ell)}$ & $\mathbf C^{(\ell)}$ \\
\hline
0 & $(1,0)$ & $(1,1)$ & $(2,1)$ \\
\hline
1 & $(\ell\!-\!1,0), (\ell,1), (\ell\!+\!1,0)$ & $(1,0)$, $(1,1)$, $(2,1)$ & $(1,1)$, $(2,0)$, $(2,1)$, $(2,2)$, $(3,1)$ \\
\hline
2 & \ditto & \begin{tabular}{@{}c@{}}
     $(\ell\!-\!1,1)$, $(\ell,0)$,  \\
     $(\ell,1)$, $(\ell\!+\!1,1)$ 
\end{tabular} & \begin{tabular}{@{}c@{}}
$(1,0)$, $(1,1)$, $2 \times (2,1)$, $(2,2)$, \\
$(3,0)$, $(3,1)$, $(3,2)$, $(4,1)$
\end{tabular} \\
\hline
$3+$ & \ditto & \ditto & \begin{tabular}{@{}c@{}}
$(\ell\!-\!2,1)$, $(\ell\!-\!1,0)$, $(\ell\!-\!1,1)$, \\
$(\ell\!-\!1,2)$, $2 \times (\ell,1)$, $(\ell,2)$,   
$(\ell\!+\!1,0)$, \\
$(\ell\!+\!1,1)$, $(\ell\!+\!1,2)$, $(\ell\!+\!2,1)$
\end{tabular} \\
\hline
\end{tabular}
\end{center}
\vspace{3mm}

The highest weight vectors of $\mathbf A^{(\ell)}_\lambda$ and $\mathbf B^{(\ell)}_\lambda$, which all have multiplicity one, are computed using the layer stripping method outlined in \Cref{sec:layer_stripping}. To obtain the corresponding vectors in $\mathbf C^{(\ell)}_\lambda$, one can use the following computational shortcut. Since $\mathbf T_1$ (as well as $\mathbf T_3$) is an $SO(5)$-invariant and degree-preserving operator, the space $\mathbf A^{(\ell)}_\lambda \oplus \mathbf B^{(\ell)}_\lambda \oplus \mathbf C^{(\ell)}_\lambda$ is invariant under $\mathbf T_1$. Therefore, applying $\mathbf T_1$ to a highest weight vector $a^{(\ell)}_\lambda \in \mathbf A^{(\ell)}_\lambda$ and projecting onto the orthogonal complement of $\{a^{(\ell)}_\lambda,b^{(\ell)}_\lambda \}$ will yield a highest weight vector in $\mathbf C^{(\ell)}_\lambda$. If all but one of the highest weight vectors in $E_\lambda$ have been determined, we may also use one of the operators $\mathbf T_2,\mathbf T_4,\mathbf T_5$ in an analogous way, making sure to project onto the orthogonal complement of all other basis vectors in $E_\lambda$.

One obtains a different basis, however still with homogeneous harmonic polynomials as coefficients, if one instead first takes the tensor product $\sov \otimes \mathcal H_\ell$. One then decomposes this representation, which corresponds to $\sov$-valued $0$-forms with fixed angular momentum, into irreducibles $K^{(\ell)}_\lambda$, and considers the irreducible components of the tensor product $K^{(\ell)}_\lambda \otimes T^\ast_0 \mathbb R^5$. In this decomposition, the operator $\mathbf T_1$ appears as a diagonal matrix, but the components seem less easy to organize. Due to the simpler structure of the matrix representing $\mathbf T_1$, we use this order of operations to handle modes of high angular momentum in \Cref{sec:high_modes}.

\begin{figure}[t]
    \centering
\begin{tikzpicture}[x=1.4cm,y=1.4cm,font=\small]
  \definecolor{latticegray}{gray}{0.75}
  \path[use as bounding box] (-2.5,-0.5) rectangle (2.5,2.5);
  \clip (-2.5,-0.5) rectangle (2.5,2.5);

  \foreach \x in {-2,...,2}
    \draw[latticegray,line width=0.4pt,dash pattern=on 2pt off 2pt]
      (\x,-0.5) -- (\x,2.5);
  \foreach \y in {0,...,2}
    \draw[latticegray,line width=0.4pt,dash pattern=on 2pt off 2pt]
      (-2.5,\y) -- (2.5,\y);

  \foreach \x in {-2,...,2} {
    \foreach \y in {1,2} {
      \ifnum\x=0
        \ifnum\y=1\else
          \fill[latticegray] (\x,\y) circle[radius=1.1pt];
        \fi
      \else
        \fill[latticegray] (\x,\y) circle[radius=1.1pt];
      \fi
    }
  }

  \foreach \x in {-2,0,2}
    \fill[latticegray] (\x,0) circle[radius=1.3pt];
  \draw[black,line width=0.6pt,-{Stealth[length=5pt,width=4pt]}]
    (-2.5,0) -- (2.5,0);
  \foreach \p in {(-1,0),(0,1),(1,0)}
    \fill[black] \p circle[radius=2.2pt];

  \node[anchor=north,inner sep=0pt,outer sep=0pt] at (-1,-0.12) {$(\ell\!-\!1,0)$};
  \node[anchor=south,inner sep=0pt,outer sep=0pt] at (0,1.12) {$(\ell,1)$};
  \node[anchor=north,inner sep=0pt,outer sep=0pt] at (1,-0.12) {$(\ell\!+\!1,0)$};
\end{tikzpicture}
    \caption{The highest weights occurring in $\mathbf A^{(\ell)}$.}
    \label{fig:A_d}
\end{figure}
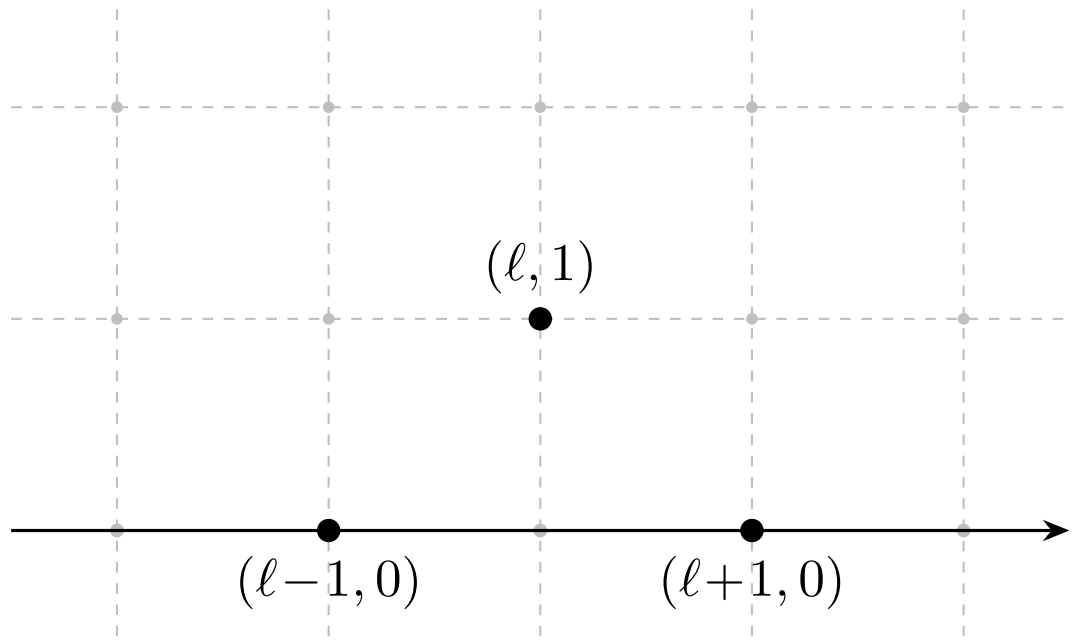

\begin{figure}[t]
    \centering
\begin{tikzpicture}[x=1.4cm,y=1.4cm,font=\small]
  \definecolor{latticegray}{gray}{0.75}
  \path[use as bounding box] (-2.5,-0.5) rectangle (2.5,2.5);
  \clip (-2.5,-0.5) rectangle (2.5,2.5);

  \foreach \x in {-2,...,2}
    \draw[latticegray,line width=0.4pt,dash pattern=on 2pt off 2pt]
      (\x,-0.5) -- (\x,2.5);
  \foreach \y in {0,...,2}
    \draw[latticegray,line width=0.4pt,dash pattern=on 2pt off 2pt]
      (-2.5,\y) -- (2.5,\y);

  \foreach \x in {-2,...,2} {
    \fill[latticegray] (\x,2) circle[radius=1.1pt];
  }
  \foreach \x in {-2,2}
    \fill[latticegray] (\x,1) circle[radius=1.1pt];

  \foreach \x in {-2,-1,1,2}
    \fill[latticegray] (\x,0) circle[radius=1.3pt];
  \draw[black,line width=0.6pt,-{Stealth[length=5pt,width=4pt]}]
    (-2.5,0) -- (2.5,0);
  \foreach \p in {(0,0),(-1,1),(0,1),(1,1)}
    \fill[black] \p circle[radius=2.2pt];

  \node[anchor=north,inner sep=0pt,outer sep=0pt] at (0,-0.12) {$(\ell,0)$};
  \node[anchor=south,inner sep=0pt,outer sep=0pt] at (-1,1.12) {$(\ell\!-\!1,1)$};
  \node[anchor=south,inner sep=0pt,outer sep=0pt] at (0,1.12) {$(\ell,1)$};
  \node[anchor=south,inner sep=0pt,outer sep=0pt] at (1,1.12) {$(\ell\!+\!1,1)$};
\end{tikzpicture}
    \caption{The highest weights occurring in $\mathbf B^{(\ell)}$.}
    \label{fig:B_d}
\end{figure}
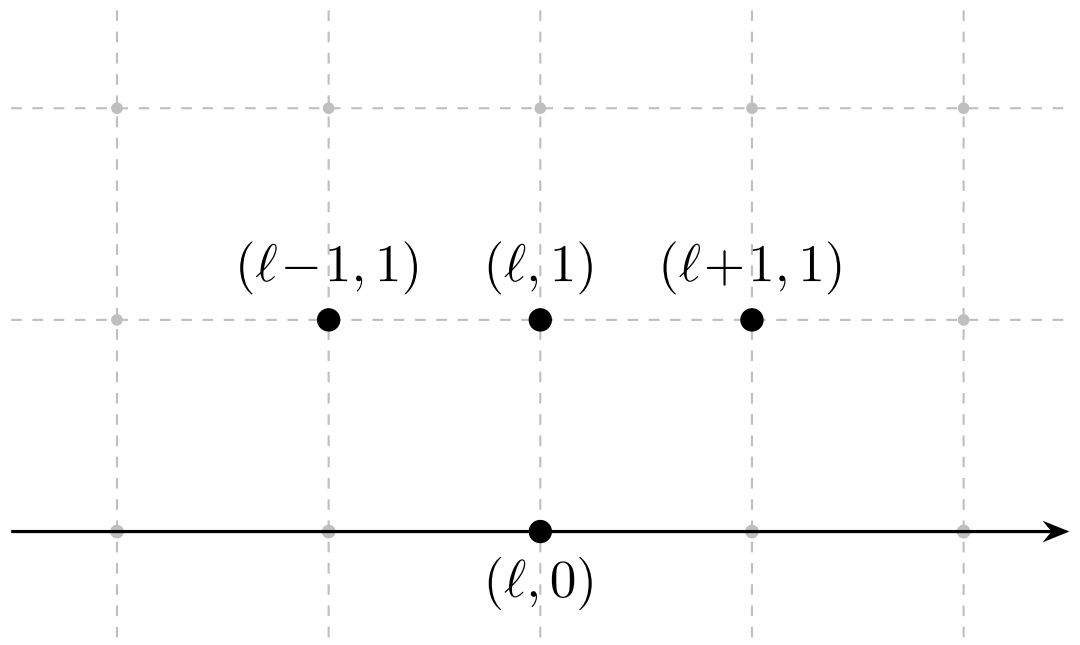

\subsection{Outline of the accompanying code}

The computations necessary for the symmetry reduction of $\mathbf L$ are contained in two annotated Jupyter notebooks, both of which import an auxiliary Python file named \texttt{SO5\_utility.py}. All calculations carried out in these files are purely symbolic and assisted by the Python package SymPy.

The notebook \texttt{SO5\_basis\_vectors.ipynb} contains a list of the highest weight vectors of the components $\mathbf A^{(\ell)}_\lambda$, $\mathbf B^{(\ell)}_\lambda$ and $\mathbf C^{(\ell)}_\lambda$ for all occurring highest weights $\lambda$ up to $\lambda = (5,1)$. It furthermore contains the code used to compute these vectors as outlined in \Cref{sec:tensor_decomposition}.

The notebook \texttt{SO5\_operators.ipynb} computes, for each highest weight $\lambda$ as above, the matrix coefficients $T_1,\ldots,T_5$ representing the invariant operators $\mathbf T_1,\ldots,\mathbf T_5$ defined in \eqref{eq:T_operators}. This yields the (matrix valued) coefficients $A_\lambda, B_\lambda$ of the Schr\"odinger operator $\mathbf L_{\lambda}$ defined in \eqref{eq:linearised op symmetry reduced II}, via the formulas \eqref{eq:def A} and \eqref{eq:def B}. The computation of these coefficients is conceptually straightforward: If $v_1,\ldots,v_m \in C^\infty(\mathbb S^4) \otimes \sov \otimes T_0^\ast \mathbb R^5$ constitute an orthonormal highest weight basis of $E_\lambda$, then for $a \in \{1,\ldots,5\}$
\begin{align*}
    (T_a)^k_\ell = \left\langle v_k, \mathbf T_a v_\ell \right\rangle = \int_{\mathbb S^4} \tr\big( v_k(\omega)_j^\ast (\mathbf T_a v_\ell)(\omega)_j \big) d\omega
\end{align*}
where $\mathbf T_a v_j$ is computed by symbolic calculation, which is possible since the coefficients of $\mathbf T_a$ and $v_j$ are polynomials in $x_1,\ldots,x_5$.

The coefficients of a vector $\Psi \in A^{\ell}_\lambda$ are elements of $\mathcal H_\ell$, and thus only defined on $\mathbb S^4$. In the accompanying Python code, these coefficients are stored as polynomials in the variables $x_1,\ldots,x_5$, so a choice of extension to $\mathbb R^5$ must be made. The most natural choice is to extend the coefficients as homogeneous polynomials of order $\ell$. One benefit of this choice is that one may easily check whether the coefficients of a vector in $\mathbf A^{(\ell)}_\lambda$ indeed lie in $\mathcal H_\ell$, since the homogeneous extension is then harmonic. After this choice has been made, it must be kept in mind when performing an orthogonal projection of a vector $v \in \mathbf A^{(\ell)}_\lambda \oplus \mathbf B^{(\ell)}_\lambda \oplus \mathbf C^{(\ell)}_\lambda$ onto the orthogonal complement of vectors $v' \in A^{(\ell')}_\lambda$, $\ell' < \ell$, a step occurring in the computation of the highest weight vectors of $\mathbf C^{(\ell)}_\lambda$ outlined in \Cref{sec:layer_stripping}. To preserve homogeneity in this step, it suffices to multiply $v'$ by $(x_1^2 + \ldots + x_5^2)^{\frac{\ell-\ell'}{2}}$.

\subsection{Sample computation}

For the reader's convenience, we demonstrate the computation of the highest weight basis of $E_{(1,1)}$ and of the matrices representing the action of the operators $\mathbf T_1,\ldots,\mathbf T_5$ on $E_{(1,1)}$. Since all constituents of $\mathbf L$ commute with point reflection across the origin, $E_{(1,1)}$ furthermore splits into an even and odd subspace, $E_{(1,1)}^{\mathrm{even}}$ and $E_{(1,1)}^{\mathrm{odd}}$, respectively. We shall discuss only the more instructive computation of $E_{(1,1)}^{\mathrm{odd}}$ here.

\subsubsection{Computation of the highest weight vectors}

According to the table in \Cref{sec:tensor_decomposition}, the space $E_{(1,1)}^{\mathrm{odd}}$ nontrivially intersects $\mathbf A^{(1)}$, $\mathbf B^{(1)}$, $\mathbf C^{(1)}$ and $\mathbf C^{(3)}$.

The computation of $\mathbf A^{(1)}$ via the layer stripping method of \Cref{sec:tensor_decomposition} begins with the highest weight vector $a^{(1)}_{(2,0)}$ of $\mathbf A^{(1)}_{(2,0)}$, explicitly given by $a^{(1)}_{(2,0)} = (x_2 - i x_3) \cdot v_{(1,0)}$, where $v_{(1,0)}$ is the highest weight vector in $V_{(1,0)}$ and given by

\vspace{1mm}
\scalebox{0.6}{
\begin{minipage}{\textwidth}
\begin{align*}
    \begin{pmatrix}0 & i & 1 & 0 & 0\\- i & 0 & 0 & 0 & 0\\-1 & 0 & 0 & 0 & 0\\0 & 0 & 0 & 0 & 0\\0 & 0 & 0 & 0 & 0\end{pmatrix}  dx^1 + \begin{pmatrix}0 & 0 & 0 & 0 & 0\\0 & 0 & 1 & 0 & 0\\0 & -1 & 0 & 0 & 0\\0 & 0 & 0 & 0 & 0\\0 & 0 & 0 & 0 & 0\end{pmatrix} dx^2 + \begin{pmatrix}0 & 0 & 0 & 0 & 0\\0 & 0 & - i & 0 & 0\\0 & i & 0 & 0 & 0\\0 & 0 & 0 & 0 & 0\\0 & 0 & 0 & 0 & 0\end{pmatrix} dx^3 + \begin{pmatrix}0 & 0 & 0 & 0 & 0\\0 & 0 & 0 & - i & 0\\0 & 0 & 0 & -1 & 0\\0 & i & 1 & 0 & 0\\0 & 0 & 0 & 0 & 0\end{pmatrix} dx^4 + \begin{pmatrix}0 & 0 & 0 & 0 & 0\\0 & 0 & 0 & 0 & - i\\0 & 0 & 0 & 0 & -1\\0 & 0 & 0 & 0 & 0\\0 & i & 1 & 0 & 0\end{pmatrix} dx^5.
\end{align*}
\end{minipage}}
\vspace{1mm} \newline

One then computes the weight space of weight $(1,1)$ in the irreducible component $\mathbf A^{(1)}_{(2,0)}$ by acting on $a^{(1)}_{(2,0)}$ with $L_{(-1,1)}$. Projecting a generic element in $\mathbf A^{(1)}$ of weight $(1,1)$, say $(x_4 - i x_5) \cdot v_{(1,0)}$ or $(x_2 - i x_3) \cdot v_{(0,1)}$, onto its orthogonal complement yields a highest weight vector $a^{(1)}_{(1,1)}$ of $\mathbf A^{(1)}_{(1,1)}$. Now one computes the weight spaces of the remaining weight $(0,0)$ in $\mathbf A^{(1)}_{(2,0)}$ and $\mathbf A^{(1)}_{(1,1)}$ and projects any vector of weight $(0,0)$ in $\mathbf A^{(1)}$ onto the orthogonal complement to obtain $a^{(1)}_{(0,0)}$. As a welcome sanity check, this vector is the (up to scale) unique $SO(5)$-invariant element of $C^\infty(\mathbb S^4) \otimes \sov \otimes T_0^\ast \mathbb R^5$, namely 

\vspace{1mm}
\scalebox{0.55}{
\begin{minipage}{\textwidth}
\begin{align*}
    \begin{pmatrix}0 & x_{2} & x_{3} & x_{4} & x_{5}\\- x_{2} & 0 & 0 & 0 & 0\\- x_{3} & 0 & 0 & 0 & 0\\- x_{4} & 0 & 0 & 0 & 0\\- x_{5} & 0 & 0 & 0 & 0\end{pmatrix}dx^1 +  \begin{pmatrix}0 & - x_{1} & 0 & 0 & 0\\x_{1} & 0 & x_{3} & x_{4} & x_{5}\\0 & - x_{3} & 0 & 0 & 0\\0 & - x_{4} & 0 & 0 & 0\\0 & - x_{5} & 0 & 0 & 0\end{pmatrix}dx^2 +  \begin{pmatrix}0 & 0 & - x_{1} & 0 & 0\\0 & 0 & - x_{2} & 0 & 0\\x_{1} & x_{2} & 0 & x_{4} & x_{5}\\0 & 0 & - x_{4} & 0 & 0\\0 & 0 & - x_{5} & 0 & 0\end{pmatrix}dx^3 +  \begin{pmatrix}0 & 0 & 0 & - x_{1} & 0\\0 & 0 & 0 & - x_{2} & 0\\0 & 0 & 0 & - x_{3} & 0\\x_{1} & x_{2} & x_{3} & 0 & x_{5}\\0 & 0 & 0 & - x_{5} & 0\end{pmatrix}dx^4 +  \begin{pmatrix}0 & 0 & 0 & 0 & - x_{1}\\0 & 0 & 0 & 0 & - x_{2}\\0 & 0 & 0 & 0 & - x_{3}\\0 & 0 & 0 & 0 & - x_{4}\\x_{1} & x_{2} & x_{3} & x_{4} & 0\end{pmatrix} dx^5.
\end{align*}
\end{minipage}}
\vspace{1mm} \newline

The weight basis of $\mathbf B^{(1)}$ is computed analogously. Next, one obtains $c^{(1)}_{(1,1)}$ by projecting $\mathbf T_1 \, a^{(1)}_{(1,1)}$ onto the orthocomplement of $a^{(1)}_{(1,1)}$ and $b^{(1)}_{(1,1)}$. Finally, $c^{(3)}_{(1,1)}$ is obtained by projecting $\mathbf T_2 \, c^{(1)}_{(1,1)}$, a form with coefficients homogeneous of order $3$, onto the orthocomplement of $|x|^2 \cdot \{a^{(1)}_{(1,1)}, b^{(1)}_{(1,1)}, c^{(1)}_{(1,1)}\}$, yielding the final element of a basis for the highest weight space of $E_{(1,1)}$. The form $c^{(3)}_{(1,1)}$ already carries some complexity, e.g., 
\begin{align*}
    \left\langle c^{(3)}_{(1,1)}, \partial_{2} \right\rangle^2_3 = (x_4 - i x_5) (5 i x_{1}^{2} - 30 i x_{2}^{2} - 35 x_{2} x_{3} + 5 i x_{3}^{2} + 5 i x_{4}^{2} + 5 i x_{5}^{2}).
\end{align*}

\subsubsection{Matrix representations of invariant operators}

With a basis of highest weight vectors $(v_k)_{k=1}^4 = \left(a^{(1)}_{(1,1)}, b^{(1)}_{(1,1)}, c^{(1)}_{(1,1)}, c^{(3)}_{(1,1)}\right)$ of $E_{(1,1)}^{\mathrm{odd}}$ in hand, determining the matrices representing the action of $\mathbf T_j$, $1 \leq j \leq 5$ on $E_{(1,1)}^{\mathrm{odd}}$ is straightforward, via
\begin{align*}
    (T_j)^{k}_{k'} = \frac{ \left\langle v_{k'}, \mathbf T_j v_k\right\rangle}{ \|v_k\| \|v_{k'}\|}.
\end{align*}
The normalization accounts for the fact that the basis obtained above is only orthogonal, but not orthonormal, reducing the computational complexity of the preceding steps considerably. The choice of extension of the coefficients of $v_k$ from $\mathbb S^4$ to $\mathbb R^5$ plays no role in this computation, since all operators under consideration are linear under multiplication by radial factors. However, the extension to homogeneous harmonic polynomials has the benefit of making clear by inspection that $\mathbf T_1$ and $\mathbf T_3$ are degree-preserving, while $\mathbf T_2,\mathbf T_4$ and $\mathbf T_5$ may increase the degree of the coefficients of $v_k$ by at most $2$. Together with the fact that the matrices representing the operators $T_j$ are symmetric, since the operators themselves are, this imposes significant restrictions on the off-diagonal terms in a weight basis with homogeneous components, such as the one chosen here.

\section{General considerations on vector-valued Schrödinger operators} \label{section:general schr}
In order to properly analyze the modes of low angular momentum, we start with some elementary theory on vector-valued Schrödinger operators of the form 
\begin{align}\label{def: basic operator}
    \Lf_0 u(r)= -\partial_r^2u(r)+\frac{1 }{r^2} \text{ diag}(a_j) u(r)+\frac{r^2}{16}u(r)+A(r)u(r).
\end{align}
on $ L^2(\R_+,\C^n)$  where $A:\R_+\to \R^{n\times n}$ is symmetric for all $r\in \R_+$ and such that all its components $ [A]_i^j=a_{ij}\in L^\infty(\R_+)$. As a domain of $\Lf_0$ we set $C_c^\infty((0,\infty),\C^n)$. 
\begin{lem}\label{lem: essentially selfadjoint}
   Assume that $a_j\geq \frac34$ for $j=1,\dots, n$. Then the operator $(\Lf_0,C_c^\infty((0,\infty),\C^n))$ is essentially self-adjoint on $L^2(\R_+,\C^n)$.
\end{lem}
\begin{proof}
  This lemma is a direct consequence of the fact that each of the operators $$L_j u(r)= -\partial_r^2u(r)+\frac{a_j }{r^2}  u(r)+\frac{r^2}{16} u(r)$$
    with domain $C^\infty_c(\R_+)$ is essentially self-adjoint, which follows from the same arguments as Lemma 3.1 in \cite{GloSch20} and since the operator $\Af u(r)= A(r) u(r)$ is a bounded self-adjoint operator on $L^2(\R_+,\C^n)$.
\end{proof}
Later on we will also need to look at operators with drift. 
\begin{lem}
    The operator $(\Lf_g,C_c^\infty((0,\infty),\C^2))$, where
    $$
     \Lf_g =\Lf_0+\Lf'_g
    $$
    where 
    \begin{align*}
      \Lf'_g u(r)=  \begin{pmatrix}
        0 & g(r)
        \\
        -g(r) &0
    \end{pmatrix}\partial_r u(r)
+
    \begin{pmatrix}
        0 & \frac{g'(r)}{2}
        \\
        -\frac{g'(r)}{2} &0
    \end{pmatrix}u(r)
    \end{align*}
    with $g\in C^1_b(\R_+,\R)$ is essentially self-adjoint in $L^2(\R_+,\C^2).$
\end{lem}
\begin{proof}
  Assume that $g\neq 0$ as otherwise the claim is trivial.  A direct computation shows that there exist constants $c_1,c_2>1$ such that 
    \begin{align*}
    \|u'\|^2_{L^2(\R_+)}&\leq (\Lf_0 u,u)+c_1\|u\|_{L^2(\R_+)}^2\leq \frac1{4\|g\|_{L^\infty(\R_+)}^2}  \|\Lf_0 u\|^2_{L^2(\R_+)}+ c_2 \|u\|^2_{L^2(\R_+)}.
    \end{align*}
    So, 
    \begin{align*}
  \|\Lf_g' u(r)\|_{L^2(\R_+)} &\leq \|g\|_{L^\infty(\R_+)}\|u'\|_{L^2(\R_+)}+ c \|u\|_{L^2(\R_+)}\leq \frac1{2}  \|\Lf_0 u\|_{L^2(\R_+)}+ c_3 \|u\|_{L^2(\R_+)}.
    \end{align*}
    Therefore, $\Lf_g'$ is relatively bounded with respect to $\Lf_0$ with relative bound $\frac12$. Hence, the Kato-Rellich theorem implies that $\Lf_g$ is essentially self-adjoint.
\end{proof}
Moving on, we will need to decompose potentials $V$ into their negative and positive parts, $V_+$, $V_-\geq 0$ (i.e. $V=V_+-V_-$). Then, to state a first result, which is concerned with the scalar case, we let $a\geq-\frac{1}{4}$ and let $V\in L^2_{loc}(\R_+,\R)$.
\begin{theorem}[\cite{GloSch20} Theorem A.1] \label{lem: scalar GGMT}
Consider a scalar-valued Schrödinger operator of the form 
\begin{align}
    \Lf u(r)= -\partial_r^2u(r)+\frac{a }{r^2} u(r)+V(r)u(r)
\end{align}
with $a>-\frac
14$ and assume that $V$ is such that the closure of $(\Lf,C^\infty_c(\R_+,\R))$, denoted by ($\mathcal{L},D(\mathcal{L}))$, is self-adjoint. Furthermore, write $V=V_+ -V_-$ where $V_{\pm }\geq 0$, and let $p> 1$. If
\begin{align*}
    \int_{\R^+}r^{2p-1}|V_-(r)|^pdr < \frac{(4 a+1)^{\frac{2p-1}{2}} p^p\Gamma(p)^2}{(p-1)^{p-1}\Gamma(2p)} 
\end{align*}
then $\sigma(\mathcal{L})\subset[0,\infty)$ and $0$ is not an eigenvalue of $\mathcal{L}$.
\end{theorem} 
Note that Theorem \ref{lem: scalar GGMT}  is only concerned with real-valued functions. However, since $V$ is real-valued, the quadratic form  $q_\Lf$ satisfies
$q_\Lf[u]=q_\Lf[\Re u]+q_\Lf[\Im u]$. 
So the inclusion $\sigma(\mathcal{L})\subset[0,\infty)$ remains valid also for complex-valued functions. Additionally, since $\Lf$ is real valued, its kernel is invariant under complex conjugation. Consequently any complex element of its kernel would give rise to a real one. Thus, also in this case $0$ is not an eigenvalue provided the inequality stated in the Theorem holds.
Hence both conclusions extend to
$L^2(\R_+;\mathbb C)$ with the same constant.
For notational convenience we define 
\[
K(a,p):= \frac{(4 a+1)^{\frac{2p-1}{2}} p^p\Gamma(p)^2}{(p-1)^{p-1}\Gamma(2p)} 
\]
We can extend this Theorem to the higher dimensional case in the following way.
\begin{lem}\label{lem: GGMT vector}
Consider the operator $\Lf_0$ as defined in Eq.~\eqref{def: basic operator} with $a_j>\frac{3}{4}$ and let $\mathcal{L}$ be its self-adjoint closure. 
    Assume that for $i,j=1,\dots n$ and $i\neq j$ there exist parameters $w_{ij}>0 $ with $w_{ij}=w_{ji}^{-1} $ and for every $i=1,\dots,n$, a $p_i\in (1,\infty)$ such that 
    \begin{align}   
  \int_{\R_+} r^{2p_i-1}\left[\left( \frac{r^2}{16}+a_{ii}(r)- \sum_{j=1, j\neq i}^n w_{ij}|a_{ij}(r)| \right)_-\right]^{p_i} dr <K(a_i,p_i) .
\end{align}
 Then $\sigma(\mathcal{L})\subset(0,\infty)$.
\end{lem}
\begin{proof}

For $i<j$, Young's inequality and $w_{ij}=w_{ji}^{-1}$ give
\[
    2a_{ij}(r) \Re (u_i u_j)
    \geq -|a_{ij}(r)|
    \bigl(w_{ij}|u_i|^2+w_{ji}|u_j|^2\bigr).
\]
Hence, summing over $i<j$ and using the symmetry of $A(r)$, we obtain
\begin{equation}\label{Eq:comp for lemma}
        (\Lf_0u,u)_{L^2}
    \geq \sum_{i=1}^n (L_i u_i,u_i)_{L^2},
\end{equation}
with
 \[ 
    L_i:=-\partial_r^2+\frac{a_i}{r^2}+W_i(r),
\]
for all $u\in C_c^\infty((0,\infty),\C^n)$, where
\[ W_i(r):=\frac{r^2}{16}+a_{ii}(r)
    -\sum_{\substack{j=1\\j\neq i}}^n w_{ij}|a_{ij}(r)|.\]
Since $W_i(r)-r^2/16$ is bounded and $a_i\geq\frac34$, each $L_i$
is essentially self-adjoint on $C_c^\infty((0,\infty))$ by
Lemma~\ref{lem: essentially selfadjoint}. By Theorem \ref{lem: scalar GGMT}, $\overline{L_i}$ is nonnegative
and has trivial kernel,
which implies that
$
\sigma(\mathcal{L})\subset [0,\infty)$ and $0$ is not an eigenvalue of $ \mathcal{L},$ as Eq.~\eqref{Eq:comp for lemma} implies 
\begin{align*}
     (\mathcal{L} u,u)_{L^2}
    \geq \sum_{i=1}^n (\overline{L_i} u_i,u_i)_{L^2}
\end{align*}
for all $u\in D(\mathcal{L})$. Next, we show that $\mathcal{L}$ has a compact resolvent which implies that its spectrum is discrete and consists of eigenvalues. For this, we note that there exists a $C>0$ such that
\[
    (\mathcal{L}u,u)_{L^2}+C\|u\|_{L^2}^2
    \geq \|u'\|_{L^2}^2+\frac{1}{16}\|ru\|_{L^2}^2
    \qquad \text{for all } u\in D(\mathcal{L}).
\]
By the Rellich compactness theorem on $(0,R)$ and the estimate
\[
    \int_R^\infty |u(r)|^2\,dr
    \leq R^{-2}\|ru\|_{L^2}^2,
\]
we conclude that $D(\mathcal{L})$ embeds compactly into $L^2(\R_+,\C^n)$.
Hence $\mathcal L$ has compact resolvent by
\cite[Theorem~XIII.64]{ReeSim78} which implies that the spectrum of $\mathcal{L}$ consists purely of eigenvalues. Hence $0 \notin \sigma(\mathcal{L}).$
\end{proof}
Moving on, we make some definitions.
\begin{defi}
We define the spaces $\mathcal{H}$ and $\mathcal{H}_{rad}$ as
\begin{align*}
\mathcal{H}&:= L^2(\R^5, e^{-\frac{|x|^2}{4}}dx)\\
\mathcal{H}_{\mathrm{rad}}&:=\{f\in \mathcal{H}: f \text{ is radial}\}.
\end{align*}
and $\mathcal{H}^n :=\bigoplus_{j=1}^n \mathcal{H}$. 
\end{defi}

\begin{defi}
We define the bounded linear operator $C:\mathcal{H}_{\mathrm{rad}}\to L^2(\R_+)$
as
$$Cu(r)=\frac{2\sqrt{2}\pi}{\sqrt3}r^{2}e^{-\frac{r^2}{8}}\widehat u(r)$$
where $\widehat u$ is the radial representative of $u$.
\end{defi}
One readily verifies the following property of the operator $C$.
\begin{lem}
    The operator $C$ is a unitary operator.
\end{lem}
Consider now the operator $(\Lf_0,C^\infty_{c}(\R^5))$
where 
$$
\Lf_0 u(x)=-\Delta u(x)+\frac{1}{2}x\cdot \nabla u(x).
$$
\begin{lem}
The operator \((\Lf_0,C_c^\infty(\mathbb R^5))\) is essentially self-adjoint on \(\mathcal H\). Let \(Y\in\mathcal H_\ell\) be a normalised spherical harmonic and set
\[
\mathcal H_Y:=\{u:u(r\omega)=f(r)Y(\omega),\ f\in\mathcal H_{\mathrm{rad}}\},
\]
where we identify radial functions with their radial profiles. Then the map \(aY\mapsto Ca\) is a unitary equivalence between the restriction of \(\overline{\Lf_0}\) to \(\mathcal H_Y\) and the closure in \(L^2(0,\infty)\) of
\[
\Lf^{(\ell)} u(r)
=-\partial_r^2 u(r)+\left[\frac{2+\ell(\ell+3)}{r^2}
-\frac54+\frac{r^2}{16}\right] u(r),
\]
with $D(\Lf^{(\ell)})= C^\infty_c(0,\infty)$.
\end{lem}
\begin{proof}
    This follows as in the proof of Lemma 3.1 in 
    \cite{GloSch20}.
\end{proof}

\subsection{On Matrix-valued Schrödinger operators with ground state eigenvalues}
Motivated by the presence of the eigenvalue $-\frac12$ in the spectrum of the operator $\Lf_{(1,0)}$ we take a look at general matrix-valued Schrödinger operators with a ground state eigenvalue. Thus, consider a canonical operator 
\begin{align*}
\Lf u(r)= -\partial_r^2 u(r)+\frac{1}{r^2}\text{diag}(a_j)u(r) +A(r) u(r)
\end{align*}
and assume that $0$ is an eigenvalue with eigenfunction $f$.
Let $Q(r) = (g_1(r),g_2(r),g_3(r))$ be a unitary matrix, where $g_1$ is chosen as the Euclidean normalisation of $f$, i.e.~ $f(r)= h(r) g_1(r)$.
For now, we are only interested in the top-left entry of $Q^T L Q$, which is given by
\begin{align*}
    g_1^T \circ L \circ g_1 = g_1^T (L g_1) - 2 g_1^T g_1' \partial_r - \partial_r^2 = - \partial_r^2 + g_1^T (L g_1)
\end{align*}
To determine $Lg_1$, we note that 
\begin{align*}
	0\cdot f =0\cdot g_1 = L f = L (h g_1) = h (L g_1) - 2 h' g_1 ' - h'' g_1.
\end{align*}
Thus 
\begin{align*}
	0= h g_1^T (L g_1) - h'',
\end{align*}
so $g_1^T (L g_1) = \frac{h''}{h}$. Thus, the top left entry of the conjugated operator is
\begin{align*}
	- \partial_r^2 + \frac{h''}{h}.
\end{align*}
Heuristically, the aim of this is to put all the bad contributions into the top left component and then estimate them by means of standard supersymmetric removal, which works as follows:

Let $v$ be such that $-v'' + V v = 0$. Then $\frac{v''}{v} = V$, so 

\[(\partial_r + \frac{v'}{v})\circ(-\partial_r + \frac{v'}{v}) =
- \partial_r^2 + \partial_r(\frac{v'}{v}) + \frac{(v')^2}{v^2} =
- \partial_r^2 + \frac{v'' v - (v')^2}{v^2} + \frac{(v')^2}{v^2} =
- \partial_r^2 + V \]

Interchanging the two factors yields
\[(- \partial_r + \frac{v'}{v}) \circ 	(\partial_r + \frac{v'}{v}) =
-\partial_r^2 - \partial_r(\frac{v'}{v}) + \frac{(v')^2}{v^2} =
-\partial_r^2 + \frac{2(v')^2 - v'' v}{v^2}\]

Applying this to the above shows we need to consider the operator 
\begin{align*}
	L_S = - \partial_r^2 + \frac{2(h')^2 - h'' h}{h^2}
\end{align*}
in order to determine the spectral gap between the ground state and the next eigenvalue of the first component.  
 However, as the conjugation by $Q$ introduces first order derivatives we now consider a (formally) self-adjoint Schrödinger operator on $L^2(\R_+)^2$ with drift. Such an operator is of the form 
\begin{align*} \Lf u(r) :=
    -\partial_r^2 u(r) &+ \frac1{r^2}\text{diag}(a_j)u(r)+A(r) u(r) \\
    &+ \begin{pmatrix}
        0 & g(r)
        \\
        -g(r) &0
    \end{pmatrix}\partial_r u(r)
+
    \begin{pmatrix}
        0 & \frac{g'(r)}{2}
        \\
        -\frac{g'(r)}{2} &0
    \end{pmatrix} u(r),
\end{align*}
where $A(r)$ is a symmetric matrix whose entries $a_{ij}$ are smooth functions on $[0,\infty).$
Furthermore, assume that $0$ is an eigenvalue with eigenfunction $f(r)=(v(r),0)$, where $v$ is non-vanishing on $(0,\infty)$. Our goal is to estimate the spectral gap between $0$ and the next spectral point. Under the assumption that $0$ is the lowest spectral point, this quantity is given by \begin{align*}
    \lambda_2=\inf_{u\in C^\infty_c(\R_+)\times C^\infty_c(\R_+): u\perp f, \|u\|=1} (\Lf u,u). 
\end{align*}
We also note that if we show that $\lambda_2>0$, it will automatically follow that 0 is the smallest spectral point.
To compute $(\Lf u,u)$, we start by making a couple of observations. 
First, we note that
\begin{align*}
  \Re\left(  \begin{pmatrix}
        0 & \frac{g'}{2}
        \\
        -\frac{g'}{2} &0
    \end{pmatrix}u,u\right)=0.
\end{align*}
Further, since $v$ is an eigenfunction, we see that 
\begin{align*}
    -g(r) v'(r)-\frac{g'(r)}2 v(r) +a_{1,2}(r)v(r)=0  \iff  a_{1,2}(r)=g(r)\frac{v'(r)}{v(r)}+\frac{g'(r)}2.
\end{align*}
as well as 
\begin{align*}
   - \partial_r^2 v(r)+a_1/r^2v(r)+a_{11}(r) v(r)=0
\end{align*}
which of course implies that  
\begin{align*}
\frac{v''(r)}{v(r)}=+a_1/r^2+a_{11}(r).
\end{align*}
Therefore, using that $(\Lf u,u)$ is real, we obtain
\begin{align*}
    (\Lf u,u)_{L^2}
    &=\int_0^\infty\left[
        |\partial_r u|^2
        +\frac{v''}{v}|u_1|^2
        +\left(\frac{a_2}{r^2}+a_{22}\right)|u_2|^2
    \right]\,dr\\
    &\quad+\Re\int_0^\infty\left[
        \left(2g\frac{v'}{v}+g'\right)u_1\overline{u_2}
        +g\left(u_2'\overline{u_1}-u_1'\overline{u_2}\right)
    \right]\,dr.
\end{align*}
Now, we integrate by parts to infer that
\begin{align*}
    \int_0^\infty\frac{v''}{v}|u_1|^2\,dr
    &=\int_0^\infty\left[
        \partial_r\left(\frac{v'}{v}\right)
        +\left(\frac{v'}{v}\right)^2
    \right]|u_1|^2\,dr\\
    &=\int_0^\infty\left[
        -2\frac{v'}{v}\Re\left(u_1'\overline{u_1}\right)
        +\left(\frac{v'}{v}\right)^2|u_1|^2
    \right]\,dr.
\end{align*}
Consequently,
\begin{align*}
    \int_0^\infty\left(
        |\partial_r u_1|^2+\frac{v''}{v}|u_1|^2
    \right)\,dr
    =\int_0^\infty
        \left|\partial_r u_1-\frac{v'}{v}u_1\right|^2\,dr.
\end{align*}
Moreover, since
\begin{align*}
    a_{12}+a_{21}=2a_{12}=2g\frac{v'}{v}+g',
\end{align*}
integrating by parts in the drift term gives
\begin{align*}
    &\Re\int_0^\infty\left[
        g\left(u_2'\overline{u_1}-u_1'\overline{u_2}\right)
        +(a_{12}+a_{21})u_1\overline{u_2}
    \right]\,dr\\
    &\qquad=\Re\int_0^\infty\left[
        -2g u_1'\overline{u_2}
        +(a_{12}+a_{21}-g')u_1\overline{u_2}
    \right]\,dr\\
    &\qquad=2\Re\int_0^\infty
        g\left(-u_1'+\frac{v'}{v}u_1\right)\overline{u_2}\,dr.
\end{align*}
Thus,
\begin{align*}
    &\int_0^\infty
        \left|\partial_r u_1-\frac{v'}{v}u_1\right|^2\,dr
      +2\Re\int_0^\infty
        g\left(-u_1'+\frac{v'}{v}u_1\right)\overline{u_2}\,dr\\
    &\qquad=\int_0^\infty\left[
        \left|-\partial_r u_1+\frac{v'}{v}u_1+g u_2\right|^2
        -g^2|u_2|^2
    \right]\,dr,
\end{align*}
and we arrive at
\begin{align*}
    (\Lf u,u)
    =\int_0^\infty\left[
        \left|-\partial_r u_1+\frac{v'}{v}u_1+g u_2\right|^2
        +|\partial_r u_2|^2
        +\left(\frac{a_2}{r^2}+a_{22}-g^2\right)|u_2|^2
    \right]\,dr.
\end{align*}
Now, we define the scalar operators \begin{align*}
    A_1 u(r)&= -\partial_r u(r)+\frac{v'(r)}{v(r)}u(r)
    \\
    L_2u(r)&=- \partial_r^2 u(r)+\frac{a_2}{r^2}u(r)+a_{22}(r) u(r)-g(r)^2 u(r). 
\end{align*}
Assume now that we can provide lower bounds on the quantities
\begin{align*}
    \alpha_1:&=\inf_{u\in C{^\infty_c}(\R_+): u\perp v, \|u\|=1} (A_1u,A_1 u)
    \\
      \alpha_2:&=\inf_{u\in C{^\infty_c}(\R_+):\|u\|=1} (L_2u,u)
      \\
      \beta:&= \left\| g\right\|_{L^\infty(\R_+)}^2
\end{align*}
and that we want $\lambda_2 \geq c$. To achieve this, we infer that 
\begin{align*}
     (\Lf u,u)\geq (1-\varepsilon) \alpha_1 \|u_1\|^2+\left(\alpha_2+(1-\varepsilon^{-1})\beta\right) \|u_2\|^2
\end{align*}
for $0<\varepsilon<1$. Hence, the claim follows if we can pick an $\varepsilon \in (0,1)$ such that
\begin{align*}
    (1-\varepsilon)\alpha_1 \geq c \iff \varepsilon\leq 1-\frac{c}{\alpha_1}
    \\
    \left(\alpha_2+(1-\varepsilon^{-1})\beta\right)\geq c \iff \frac{\beta}{\alpha_2+\beta-c}\leq \varepsilon.
\end{align*}
Thus, if the bounds
\begin{align*}
 \frac{\beta}{\alpha_2+\beta-c}    \leq 1-\frac{c}{\alpha_1}\iff \beta \leq \frac{1}{c}(\alpha_1-c)(\alpha_2-c)
\end{align*}
and 
$$
c< \min\{\alpha_1,\alpha_2\}
$$
hold, then the spectral gap between $\lambda_1$ and $\lambda_2$ is at least $c$.
We summarise this in the following proposition.
\begin{prop}\label{prop:abstract spectral prop}
Consider the operator $(\Lf, C^\infty_c(\R_+,\C^2))$ with 
\begin{align*} \Lf u(r):=
    -\partial_r^2 u(r) &+\frac{1}{r^2}\text{diag} (a_j) u(r)+A(r)u(r) \\
    &+ \begin{pmatrix}
        0 & g(r)
        \\
        -g(r) &0
    \end{pmatrix}\partial_r u(r)
+
    \begin{pmatrix}
        0 & \frac{g'(r)}{2}
        \\
        -\frac{g'(r)}{2} &0
    \end{pmatrix} u(r)
\end{align*}
and assume that all components of $\text{diag} (a_j)+A(r)$ are such that $\Lf$ is essentially self-adjoint with self-adjoint closure $\mathcal{L}$.
 Furthermore, suppose that it has the ground state eigenvalue $0$ with eigenfunction $f(r)=(v(r),0)$, with $v\neq 0$ on $(0,\infty)$, and define 
\begin{align*}
    A_1 u(r)&= -\partial_r u(r)+\frac{v'(r)}{v(r)}u(r)
    \\
    L_2u(r)&=- \partial_r^2 u(r)+\frac{a_2}{r^2}u(r)+a_{22}(r) u(r)-g^2(r) u(r)
\end{align*}
as well as
\begin{align*}
    \alpha_1:&=\inf_{u\in C{^\infty_c}(\R_+)\: u\perp v, \|u\|=1} (A_1u,A_1 u)
    \\
      \alpha_2:&=\inf_{u\in C{^\infty_c}(\R_+):\|u\|=1} (L_2u,u)
      \\
      \beta:&= \left\| g\right\|_{L^\infty(\R_+)}^2.
\end{align*}
Then 
\[
\quad
0<c<\min\{\alpha_1,\alpha_2\},\qquad
\beta\leq \frac1c(\alpha_1-c)(\alpha_2-c)
\quad\Longrightarrow\quad
\inf\sigma\!\left(\mathcal L\big|_{f^\perp}\right)\ge c.
\]
Furthermore, if $\beta < \frac1c(\alpha_1-c)(\alpha_2-c)$, then $\lambda_2>c $.
\end{prop}
With this we return to our explicit operators.
\section{Modes of low angular momentum}\label{section: low angulars}
Before we come to the first of our radial Schrödinger operators, we want to briefly discuss our methods of excluding non positive eigenvalues of the operators that are not strictly positive. We use two methods to study the non-positive spectrum of the radial
Schrödinger operators in this section. The first uses the scalar and
matrix GGMT bounds from Theorem \ref{lem: scalar GGMT}
and Lemma \ref{lem: GGMT vector}, specialized to $p=4$ and potentials of the form
\begin{align*}
    V(r)= -\frac{3}{4}+\frac{r^2}{16}-\frac{\alpha r^2  + \beta}{(ar^2 + b)^2}
\end{align*}
with $\alpha,\beta \in \R$ and $a=\frac{\sqrt{3}}{2\sqrt{2}}, b=9-7\sqrt{3/2}.$ While we in general evaluate the integrals numerically, we note that these integrals can be computed rigorously. Namely, there is a unique partial fraction decomposition of the form  
\begin{align*}
    r^7 V(r)^4
=
\sum_{j=0}^{7} c_j r^{2j+1}
+
\sum_{j=1}^{8}\frac{k_j r}{(ar^2+b)^j},
\end{align*}
with $c_j, k_j \in \R$. 
A second strategy implemented here is to construct potentials such that the spectrum associated Schrödinger $-\partial_r^2+V(r)$  is contained in $[0,\infty)$, and such $V$ provides a strict lower bound for the matrix potential that we actually want to study. All computations related to GGMT bounds or any other estimates from this section can be found in the supplemental Mathematica notebook \texttt{computations\_low\_momenta.nb}.
To begin, we first quickly look at the $(0,0)$ component, which was already dealt with in \cite{GloSch20}.

\subsubsection*{The $(0,0)$-component}
The $(0,0)$ component appears only in the $\ell=1$ sector and  reduces to the scalar-valued equation 
\begin{align*}
    \widetilde\Lf_{(0,0)} u(x) &= \left(-\Delta+\frac12 \Lambda + \frac12 - \frac{|x|^2 (-9 + \frac{9 \sqrt{6}}{2} ) + 162 - 63 \sqrt{6}}{(a|x|^2 + b)^2}\right)u(x)
\end{align*}
Conjugation with respect to $C$ produces the operator 
\begin{align*}
  \Lf_{(0,0)} u(r)=:C \widetilde \Lf_{(0,0)} C^{-1} u(r) &= \left(-\partial_r^2+\frac{6}{r^2}-\frac{3}{4}+\frac{r^2}{16}-\frac{r^2 (-9 + \frac{9 \sqrt{6}}{2} ) + 162 - 63 \sqrt{6}}{(a r^2 + b)^2}\right)u(r)
\end{align*}
\begin{lem}
The operator $(\Lf_{(0,0)},C^\infty_{c}(\R_+,\R))$ is essentially self-adjoint with self-adjoint closure $\mathcal{L}_{(0,0)}$.
Furthermore, the spectrum of $\mathcal{L}_{(0,0)}$ satisfies
\begin{align*}
\sigma(\mathcal{L}_{(0,0)})\setminus \{-1\}\subset (0,\infty).
\end{align*}
\end{lem}
\begin{proof}
    This is part of the statement of Proposition 3.2 in \cite{GloSch20}.
\end{proof}
With this we come to components that require work. To illustrate the general behavior, we start with the scalar-valued odd $(1,0)$ component.
\subsubsection*{The odd $(1,0)$ component}
The odd $(1,0)$ component is given by
\begin{align*}
    \widetilde\Lf_{(1,0)}^{\mathrm{odd}}u(x) &= \left(-\Delta + \frac12 \Lambda + \frac12 - \frac{|x|^2 (-6 + 3\sqrt{6}) + 108 - 42 \sqrt{6}}{(a|x|^2 + b)^2}\right)u(x)
\end{align*}
which, given that it only appears for $\ell=1$, conjugates to
\begin{align*} 
\Lf_{(1,0)}^{\mathrm{odd}} u(r)=C \widetilde \Lf_{(1,0)}^{\mathrm{odd}} C^{-1} u(r)= \left(-\partial_r^2+\frac{6}{r^2}-\frac{3}{4}+\frac{r^2}{16}-\frac{r^2 (-6 + 3\sqrt{6}) + 108 - 42 \sqrt{6}}{(ar^2 + b)^2}\right)u(r)
\end{align*}
\begin{lem}
The operator $(\Lf_{(1,0)}^{\mathrm{odd}},C^\infty_{c}(\R_+,\R))$ is essentially self-adjoint with self-adjoint closure $\mathcal{L}_{(1,0)}^{\mathrm{odd}}$.
Furthermore, the spectrum of $\mathcal{L}_{(1,0)}^{\mathrm{odd}}$ satisfies
\begin{align*}
\sigma(\mathcal{L}_{(1,0)}^{\mathrm{odd}})\subset (0,\infty).
\end{align*}
\end{lem}
\begin{proof}
Essential self-adjointness is a direct consequence of Lemma \ref{lem: essentially selfadjoint}. The positivity of the spectrum follows immediately from Theorem \ref{lem: scalar GGMT}, applied with $p=4$, since 
$$
V_{(1,0)}^{\mathrm{odd}}(r):= -\frac{3}{4}+\frac{r^2}{16}-\frac{r^2 (-6 + 3\sqrt{6}) + 108 - 42 \sqrt{6}}{(ar^2 + b)^2}
$$
is negative up to its first and only root $r_0$ which satisfies $r_0\leq 4.1$
and
$$
\int_0^\infty r^7 V_{(1,0)}^{\mathrm{odd}}(r)_-^4 dr \leq  \int_0^{4.1} r^7 V_{(1,0)}^{\mathrm{odd}}(r)^4 dr \leq 4000 < 5000 \leq  \frac{25^{\frac 72} \cdot 256 \Gamma(4)^2}{27\Gamma(8)} .
$$
\end{proof}
Next, given that the operator associated with even $(1,0)$ component has the symmetry eigenvalue $-\frac{1}{2}$, thus requires special treatment, we will treat it separately at the end of this section and continue 
with the $(1,1)$ components. 
\subsubsection*{The even $(1,1)$ component}
The even $(1,1)$ component is given by
\begin{align*}
   \widetilde  \Lf_{(1,1)}^{\mathrm{even}} &= -\partial_r^2 - \frac{4}{r}\partial_r + \frac1{r^2} \begin{pmatrix}
    0 & 0 & 0 \\
    0 & 10 & 0 \\
    0 & 0 & 10 \\
\end{pmatrix} +\frac12 r\partial_r + \frac12 - \frac{1}{(ar^2 + b)^2}( r^ 2 A_{(1,1)}^{\mathrm{even}}+B_{(1,1)}^{\mathrm{even}})
\end{align*}
where
\begin{align*}
A_{(1,1)}^{\mathrm{even}}=\begin{pmatrix}- \frac{24}{5} + \frac{6 \sqrt{6}}{5} & - \frac{8}{5} + \frac{2 \sqrt{6}}{5} & - \frac{\sqrt{30}}{5} + \frac{2 \sqrt{5}}{5}\\- \frac{8}{5} + \frac{2 \sqrt{6}}{5} & - \frac{26}{5} + \frac{47 \sqrt{6}}{15} & - \frac{2 \sqrt{5}}{5} + \frac{7 \sqrt{30}}{15}\\- \frac{\sqrt{30}}{5} + \frac{2 \sqrt{5}}{5} & - \frac{2 \sqrt{5}}{5} + \frac{7 \sqrt{30}}{15} & -3 + \frac{13 \sqrt{6}}{6}\end{pmatrix}
\\
B_{(1,1)}^{\mathrm{even}}=\begin{pmatrix}72 - 28 \sqrt{6} & 0 & 0\\0 & 132 - \frac{154 \sqrt{6}}{3} & - \frac{14 \sqrt{30}}{3} + 12 \sqrt{5}\\0 & - \frac{14 \sqrt{30}}{3} + 12 \sqrt{5} & 30 - \frac{35 \sqrt{6}}{3}\end{pmatrix}
\end{align*}
which conjugates to 
  \begin{align*}
      \Lf_{(1,1)}^{\mathrm{even}}= -\partial_r^2+   \frac1{r^2} \begin{pmatrix}
    2 & 0 & 0 \\
    0 & 12 & 0 \\
    0 & 0 & 12 
    \end{pmatrix}+\frac{r^2}{16}-\frac{3}{4}- \frac{1}{(ar^2 + b)^2}( r^ 2 A_{(1,1)}^{\mathrm{even}}+B_{(1,1)}^{\mathrm{even}}).
  \end{align*}
\begin{lem}
The operator $(\Lf_{(1,1)}^{\mathrm{even}},C^\infty_{c}(\R_+,\R^3))$ is essentially self-adjoint with self-adjoint closure $\mathcal{L}_{(1,1)}^{\mathrm{even}}$.
Furthermore, the spectrum of $\mathcal{L}_{(1,1)}^{\mathrm{even}}$ satisfies
\begin{align*}
\sigma(\mathcal{L}_{(1,1)}^{\mathrm{even}})\subset (0,\infty).
\end{align*}
\end{lem}
\begin{proof}
 As before, this is  a direct consequence of GGMT bounds. To apply Lemma \ref{lem: GGMT vector} in a succinct fashion, we define $c_{ij}(r)=-\frac{1}{( ar^2+b)^2}\left(r^2 [A_{(1,1)}^{\mathrm{even}}]_i^j+[B_{(1,1)}^{\mathrm{even}}]^j_i\right) $ and
 \begin{align*}
 V_1&:=\frac{r^2}{16}-\frac{3}{4}+ c_{11}(r)- |c_{12}(r)|- |c_{13}(r)|
 \\
 V_2&:= \frac{r^2}{16}-\frac{3}{4}+c_{22}(r)- |c_{12}(r)|- \frac{3}{4}|c_{23}(r)|
 \\
 V_3&:=\frac{r^2}{16}-\frac{3}{4}+ c_{33}(r)- |c_{13}(r)|- \frac{4}{3}|c_{23}(r)|,
 \end{align*}
explicitly
\begin{align*}
V_1(r)&=\frac{r^2}{16}-\frac34
  +\frac{\bigl(16-4\sqrt6+2\sqrt5-\sqrt{30}\bigr)r^2-40b}
        {5(ar^2+b)^2},
\\
V_2(r)&=\frac{r^2}{16}-\frac34 +\frac{\bigl(216-164\sqrt6+18\sqrt5-21\sqrt{30}\bigr)r^2 -(880+60\sqrt5)b}{60(ar^2+b)^2},
\\
V_3(r)&=\frac{r^2}{16}-\frac34+\frac{\bigl(270-195\sqrt6+84\sqrt5-74\sqrt{30}\bigr)r^2 -(300+160\sqrt5)b}{90(ar^2+b)^2}.
\end{align*}
Now, one computes that each potential $V_j$ is negative precisely on $[0,r_j)$  where $r_j$ is its unique positive root.
These roots satisfy $r_1< 3.1, r_2<4.6, r_3< 4.7$. With this at hand, we compute that
\begin{align*}
 \int_0^{3.1} r^7 V_1(r)^4dr \leq  84<  K(2,4)
 \\
  \int_0^{4.6} r^7 V_2(r)^4dr \leq  46000<  55000 \leq K(12,4)
  \\
    \int_0^{4.7} r^7 V_3(r)^4dr \leq  29000< 55000 \leq K(12,4).
\end{align*}
Consequently, the desired spectral inclusion $\sigma( \mathcal{L}_{(1,1)}^{\mathrm{even}})\subset [0,\infty)$ and the fact that $0$ is not an eigenvalue follows from Lemma \ref{lem: GGMT vector}.
\end{proof}
\subsubsection*{The odd (1,1) component}
We continue with the odd $(1,1)$ component which conjugates to
\begin{align*}
    \Lf_{(1,1)}^{\mathrm{odd}}&= -\partial_r^2 +  \frac1{r^2} \begin{pmatrix}
    6 & 0 & 0 & 0 \\
    0 & 6 & 0 & 0 \\
    0 & 0 & 6 & 0 \\
    0 & 0 & 0 & 20 \\
\end{pmatrix}
+\frac{r^2}{16}-\frac 34
- \frac{1}{(ar^2 + b)^2} \left( r^2 A_{(1,1)}^{\mathrm{odd}}+ B_{(1,1)}^{\mathrm{odd}} \right)
\end{align*}
where
\begin{align*}
A_{(1,1)}^{\mathrm{odd}} &= \left(\begin{smallmatrix}- \frac{21}{4} + \frac{15 \sqrt{6}}{8} & \sqrt{3} & - \frac{9 \sqrt{7}}{28} + \frac{19 \sqrt{42}}{56} & - \frac{3 \sqrt{14}}{7} + \frac{3 \sqrt{21}}{14}\\\sqrt{3} & -4 + \frac{4 \sqrt{6}}{3} & - \frac{11 \sqrt{21}}{21} + \frac{2 \sqrt{14}}{7} & - \frac{2 \sqrt{7}}{7} + \frac{\sqrt{42}}{7}\\- \frac{9 \sqrt{7}}{28} + \frac{19 \sqrt{42}}{56} & - \frac{11 \sqrt{21}}{21} + \frac{2 \sqrt{14}}{7} & - \frac{93}{28} + \frac{229 \sqrt{6}}{168} & - \frac{\sqrt{2}}{7} + \frac{3 \sqrt{3}}{14}\\- \frac{3 \sqrt{14}}{7} + \frac{3 \sqrt{21}}{14} & - \frac{2 \sqrt{7}}{7} + \frac{\sqrt{42}}{7} & - \frac{\sqrt{2}}{7} + \frac{3 \sqrt{3}}{14} & - \frac{24}{7} + \frac{24 \sqrt{6}}{7}\end{smallmatrix}\right), & 
\\ 
B_{(1,1)}^{\mathrm{odd}} &= \left(\begin{smallmatrix}\frac{243}{2} - \frac{189 \sqrt{6}}{4} & - 14 \sqrt{3} + 18 \sqrt{2} & - \frac{7 \sqrt{42}}{4} + \frac{9 \sqrt{7}}{2} & 0\\- 14 \sqrt{3} + 18 \sqrt{2} & 96 - \frac{112 \sqrt{6}}{3} & - 6 \sqrt{14} + \frac{14 \sqrt{21}}{3} & 0\\- \frac{7 \sqrt{42}}{4} + \frac{9 \sqrt{7}}{2} & - 6 \sqrt{14} + \frac{14 \sqrt{21}}{3} & - \frac{3}{2} + \frac{7 \sqrt{6}}{12} & 0\\0 & 0 & 0 & 72 - 28 \sqrt{6}\end{smallmatrix}\right), & 
\\
S_{(1,1)}^{\mathrm{odd}} &= \mathrm{diag}(6,6,6,20). &
\end{align*}
Now, these matrix potentials appear slightly too negative for GGMT bounds to be strong enough to prove the absence of negative eigenvalues. Consequently, we will construct a comparison potential with a ground state eigenfunction to obtain that result. Now, one could certainly let AI brute force the construction of such a potential. However, we rather do it somewhat by hand. To motivate our construction, we first compute that the leading order terms of the first 3 rows near the origin are given by 
$$
-\partial_r^2 u(r) +\frac{6}{r^2} u(r).
$$
Since the fundamental solutions of the equation $f''(r)-\frac{6}{r^2}f(r)=0$ are given by $r^3$ and $r^{-2}$, we aim for a ground state with $r^3$ behavior at the origin. Similarly, near $\infty$, the leading order terms are given by 
\[
-\partial_r^2 +\frac{r^2}{16}-\frac34.
\]
Plugging the Ansatz $ r^p e^{-\frac{r^2}{8}}$ into this equation yields  
\begin{align*}
\left(
-\partial_r^2+\frac{r^2}{16}-\frac34
\right)\left(e^{-r^2/8}r^p\right)
=
e^{-r^2/8}r^p\left(
\frac{p-1}{2}-
p(p-1)r^{-2}
\right)
\end{align*}
So, we see that the leading order behavior near $\infty$ should be $r e^{-\frac{r^2}{8}}$. Now the rational structure of the potential therefore suggests the ansatz
\begin{align*}
f_c(r)= r^3e^{-\frac{r^2}{8}}\left((ar^2+b)^{-1}+c( ar^2+b)^{-2}\right)
\end{align*}
where $c$ is a degree of freedom. It turns out that $c=2$ is a choice that works. Thus, 
we define 
\begin{align*}
f(r)= r^3e^{-\frac{r^2}{8}}((ar^2+b)^{-1}+2( ar^2+b)^{-2})
\end{align*}
and $V_{(1,1)}(r)=\frac{f''}{f}$, which explicitly  reads as 
\begin{align*}
V_{(1,1)}(r):=\left(
\frac{3}{r}-\frac{r}{4}
+\frac{2ar}{ar^2+b+2}
-\frac{4ar}{ar^2+b}
\right)^2
-\frac{3}{r^2}-\frac14
+\frac{2a(b+2-ar^2)}{(ar^2+b+2)^2}
-\frac{4a(b-ar^2)}{(ar^2+b)^2}.
\end{align*}
\begin{lem}
The operator $(\Lf_{(1,1)}^{\mathrm{odd}},C^\infty_{c}(\R_+,\R^4))$ is essentially self-adjoint with self-adjoint closure $\mathcal{L}_{(1,1)}^{\mathrm{odd}}$.
Furthermore, the spectrum of $\mathcal{L}_{(1,1)}^{\mathrm{odd}}$ satisfies
\begin{align*}
\sigma(\mathcal{L}_{(1,1)}^{\mathrm{odd}})\subset (0,\infty).
\end{align*}
\end{lem}
\begin{proof}
As elaborated, the idea is to show that $V_{(1,1)}^{\mathrm{odd}}(r) -V_{(1,1)}(r) \mathbf I_4$ where
\begin{align*}
 V_{(1,1)}^{\mathrm{odd}}= \frac1{r^2} \text{diag}(6,6,6,20)
+\frac{r^2}{16}-\frac 34
- \frac{1}{(ar^2 + b)^2} \left( r^2 A_{(1,1)}^{\mathrm{odd}}+ B_{(1,1)}^{\mathrm{odd}} \right).
\end{align*}
is a positive definite matrix for all $r>0$. 
For this, we first turn the coefficients into polynomials by defining
\begin{align*}
M_{(1,1)}^{\mathrm{odd}}(r):= r^2 (a r^2+b)^2(a r^2+b+2)^2[ V_{(1,1)}^{\mathrm{odd}}-V_{(1,1)} \mathbf I_4].
\end{align*}
Then, one can apply Sylvester's criterion to establish that $M_{(1,1)}^{\mathrm{odd}}(r)$ is indeed positive definite for $r>0$. This works out nicely as every principal minor is a polynomial with nonnegative coefficients and at least one positive coefficient.
These computations can be found in the attached notebook \texttt{computations\_low\_momenta.nb}. 
\end{proof}
\subsubsection*{The positive definite low angular momentum operators}
Here, we deal with the remaining low angular momentum operators aside from the even $(1,0)$ one as this one has an unstable eigenvalue. We group the remaining ones together as it turns out that all of these operators are of the form $\Lf u(r)=-\partial_r^2 u(r)+V(r) u(r)$ with positive definite matrix potentials with leading order asymptotic $\frac{r^2}{16}I$ as $r\to \infty$. 

\begin{lem}
For $(i,j)\in \{(2,0),(2,1),(3,0),(3,1),(2,2),(4,0),(4,1),(3,2),(4,2),(5,0),(5,1)\}$, and $p=\mathrm{odd},\mathrm{even},$ the operators $(\Lf_{(i,j)}^p,C^\infty_{c}(\R_+,\R^{d_{i,j,p}}))$ are essentially self-adjoint with self-adjoint closure $\mathcal{L}_{(i,j)}^p$.
Furthermore, the spectrum of $\mathcal{L}_{(i,j)}^p$ satisfies
\begin{align*}
\sigma(\mathcal{L}_{(i,j)}^p)\subset (0,\infty)
\end{align*}
\end{lem}
\begin{proof}
For each of the above operators, define
\[
    M^\mathrm{odd}_{(2,0)}(r) := r^2(ar^2+b)^2(ar^2+b+2)^6 V^\mathrm{odd}_{(2,0)}(r)
\]
and 
\[
    M^p_{(i,j)}(r) := r^2(ar^2+b)^2(ar^2+b+2)^4 V^p_{(i,j)}(r)
\]
for the remaining components. 
Since the prefactor $r^2(ar^2+b)^2(ar^2+b+2)^{j}$ is positive for every $j \in \mathbb N$ and $r>0$, it suffices to prove that
$M^p_{(i,j)}(r)$ is positive definite. As verified in the accompanying
Mathematica notebook, every leading principal minor of $M^p_{(i,j)}(r)$
is a polynomial in $r^2$ with nonnegative coefficients and at least one
strictly positive coefficient. Thus all leading principal minors are
strictly positive for every $r>0$, and the claim follows from
Sylvester's criterion.
\end{proof}

We move on to the component with a negative spectral point, i.e.~the even $(1,0)$ component.
\subsubsection*{The even $(1,0)$ component}
We start by conjugating the even $(1,0)$ component to canonical form to arrive at
\begin{align*}
   \Lf_{(1,0)}^{\mathrm{even}}
&= -\partial_r^2 + \frac{1}{r^2}S_{(1,0)}^{\mathrm{even}} + \frac{r^2}{16}-\frac34- \frac{r^2 A_{(1,0)}^{\mathrm{even}}+B_{(1,0)}^{\mathrm{even}}}
         {(ar^2+b)^2}, \\
S_{(1,0)}^{\mathrm{even}} &= \begin{pmatrix}
    2 & 0 & 0 \\
    0 & 12 & 0 \\
    0 & 0 & 12 \\
\end{pmatrix} \\
A_{(1,0)}^{\mathrm{even}} &= \begin{pmatrix}-6 + \frac{9 \sqrt{6}}{5} & - \frac{3}{2} + \frac{3 \sqrt{6}}{5} & - \frac{3 \sqrt{15}}{10} + \frac{3 \sqrt{10}}{5}\\- \frac{3}{2} + \frac{3 \sqrt{6}}{5} & - \frac{33}{4} + \frac{183 \sqrt{6}}{40} & - \frac{27 \sqrt{10}}{40} + \frac{3 \sqrt{15}}{20}\\- \frac{3 \sqrt{15}}{10} + \frac{3 \sqrt{10}}{5} & - \frac{27 \sqrt{10}}{40} + \frac{3 \sqrt{15}}{20} & - \frac{15}{4} + \frac{21 \sqrt{6}}{8}\end{pmatrix} \\
B_{(1,0)}^{\mathrm{even}} &= \begin{pmatrix}108 - 42 \sqrt{6} & 0 & 0\\0 & \frac{351}{2} - \frac{273 \sqrt{6}}{4} & - \frac{9 \sqrt{15}}{2} + \frac{21 \sqrt{10}}{4}\\0 & - \frac{9 \sqrt{15}}{2} + \frac{21 \sqrt{10}}{4} & \frac{81}{2} - \frac{63 \sqrt{6}}{4}\end{pmatrix}
\end{align*}
with eigenvalue $-\frac12$ and associated eigenfunction 
\begin{align*}
    f_{(1,0)}^{\mathrm{even}}(r)&:=\frac{r^2e^{-\frac{r^2}{8}}}{(ar^2+b)^2} \left(
                                r^2 \begin{pmatrix}\frac{12}{5} - \frac{6 \sqrt{6}}{5}\\
                                \frac{3 \sqrt{6}}{5} - \frac{6}{5}\\
                                \frac{2 \sqrt{15}}{5} - \frac{3 \sqrt{10}}{5}\end{pmatrix}
                                +   \begin{pmatrix}
                                        -72 + 28 \sqrt{6}  \\  0 \\  0
                                    \end{pmatrix}.
                            \right)
\end{align*}
As outlined above we set 
\begin{align*}
  g_1(r)=\frac{ f_{(1,0)}^{\mathrm{even}}(r)}{| f_{(1,0)}^{\mathrm{even}}(r)|}=\mu(r)^{-\frac{1}{2}}\left(
                                r^2 \begin{pmatrix}\frac{12}{5} - \frac{6 \sqrt{6}}{5}\\
                                \frac{3 \sqrt{6}}{5} - \frac{6}{5}\\
                                \frac{2 \sqrt{15}}{5} - \frac{3 \sqrt{10}}{5}\end{pmatrix}
                                +   \begin{pmatrix}
                                        -72 + 28 \sqrt{6}  \\  0 \\  0
                                    \end{pmatrix}
                            \right)
\end{align*}
where \begin{align*}
    \mu(r)=\frac{24}{5} (5 r^4-2 \sqrt{6} r^4  + 64 \sqrt{6} r^2 -156 r^2 + 2060- 840 \sqrt{6} ).
\end{align*}
As a second basis vector, we take
\[
    g_2:=\frac{1}{\sqrt8}
    \begin{pmatrix}0\\ \sqrt5\\ \sqrt3\end{pmatrix}.
\]
One readily computes that $g_1(r)\cdot g_2=0 $ for all $r\geq 0$ and
 by taking the cross product of the two we obtain 
\begin{align*}
    g_3(r)
    :=\frac{1}{5\sqrt{\mu(r)}}
    \begin{pmatrix}
        (12-4\sqrt6)r^2 
        \\
        (9-3\sqrt6)r^2+90\sqrt6-210 
        \\
        (3\sqrt{10}-3\sqrt{15})r^2
        +70\sqrt{15}-90\sqrt{10}
    \end{pmatrix}.
\end{align*}
Thus we want to conjugate $L_{(1,0)}^{\mathrm{even}}$ by the unitary matrix $Q$ whose columns consist of these basis vectors i.e.
$   Q(r):= ( g_1(r)\quad g_2 \quad g_3(r))$
with inverse
$Q^{-1}(r)=Q(r)^T$.
We begin by computing
\begin{align*}
    Q^T\partial_r^2 [Q(r)u(r)]=\partial_r^2 u(r)+2Q^T Q'(r)u'(r)+Q^T Q'' u(r).
\end{align*}
To facilitate the computation of these terms, we note that $(g_j,g_j')=0$ and since \(g_2\) is constant, differentiating \((g_1,g_2)=(g_3,g_2)=0\) gives \((g_1',g_2)=(g_3',g_2)=0\)
Thus,
\begin{align*}
    Q^T Q'=\left(\begin{smallmatrix}
        0 &0 &(g_1, g_3')
        \\
        0 &0 &0
        \\
        -(g_1, g_3') &0 &0
    \end{smallmatrix}\right)
,\qquad 
    Q^TQ''=\left(\begin{smallmatrix}
        (g_1,g_1'') & 0 & (g_1,g_3'')
        \\
        0 &0 &0
        \\
        (g_3,g_1'') &0 &(g_3,g_3'')
    \end{smallmatrix}\right)= \left( \begin{smallmatrix}
 -\|g_1'(r)\|_E^2 & 0 & (g_1,g_3'')
        \\
        0 &0 &0
        \\
        (g_3,g_1'') &0 &-\|g_3'(r)\|_E^2
     \end{smallmatrix}\right).
\end{align*}
With this in mind, we come to the conjugation of the involved matrices. 
We start with 
\begin{align*}
    &\quad Q^T S_{(1,0)}^{\mathrm{even}} Q=E(r)=S_{(1,0)}^{\mathrm{even}}
+
\frac{48(5-2\sqrt6)r^2}{5\mu(r)}
\begin{pmatrix}
2r^2
&
0
&
40-10\sqrt6+\sqrt6\,r^2
\\
0&0&0
\\
40-10\sqrt6+\sqrt6\,r^2
&
0
&
-2r^2
\end{pmatrix}
\end{align*}
while
$$
- \frac{1}{(ar^2 + b)^2}Q^T  ( r^2 A_{(1,0)}^{\mathrm{even}} + B_{(1,0)}^{\mathrm{even}})Q=M(r)=(m_1(r),m_2(r),m_3(r))
$$
with 
\begin{align*}
m_1(r)&=
\begin{pmatrix}
-\frac{1152\left(144720-59080\sqrt6+4(-4534+1851\sqrt6)r^2+(916-374\sqrt6)r^4+(-22+9\sqrt6)r^6\right)}
{5\mu(r)\left(36-14\sqrt6+\sqrt6\,r^2\right)^2}
\\
0
\\
\frac{384\sqrt{15+6\sqrt6}\,r^2\left(-2244+916\sqrt6+(852-348\sqrt6)r^2+(-27+11\sqrt6)r^4\right)}
{5\mu(r)\left(36-14\sqrt6+\sqrt6\,r^2\right)^2}
\end{pmatrix}
\\
m_2(r)&=
\begin{pmatrix}
0
\\
-\frac{48\left(36-14\sqrt6+(-2+\sqrt6)r^2\right)}
{\left(36-14\sqrt6+\sqrt6\,r^2\right)^2}
\\
-\frac{24\left(20(103\sqrt2-84\sqrt3)-4(39\sqrt2-32\sqrt3)r^2+(5\sqrt2-4\sqrt3)r^4\right)}
{\left(36-14\sqrt6+\sqrt6\,r^2\right)^2\sqrt{\frac{5}{24}\mu(r)}}
\end{pmatrix}
\\
m_3(r)&=
\begin{pmatrix}
\frac{384\sqrt{15+6\sqrt6}\,r^2\left(-2244+916\sqrt6+(852-348\sqrt6)r^2+(-27+11\sqrt6)r^4\right)}
{5\mu(r)\left(36-14\sqrt6+\sqrt6\,r^2\right)^2}
\\-\frac{24\left(20(103\sqrt2-84\sqrt3)-4(39\sqrt2-32\sqrt3)r^2+(5\sqrt2-4\sqrt3)r^4\right)}
{\left(36-14\sqrt6+\sqrt6\,r^2\right)^2\sqrt{\frac{5}{24}\mu(r)}}
\\
\frac{48\left(1432560-584840\sqrt6+4(-54866+22399\sqrt6)r^2+(11596-4734\sqrt6)r^4+(-218+89\sqrt6)r^6\right)}
{\left(36-14\sqrt6+\sqrt6\,r^2\right)^2\left(-20380+8320\sqrt6+(1548-632\sqrt6)r^2+(-49+20\sqrt6)r^4\right)}
\end{pmatrix}.
\end{align*}
Next, we note that the first and second components are not coupled together. Motivated by this, we define another differential operator
\begin{align*}
\widetilde{\Lf}_{(1,0)}^{\mathrm{even}} u(r)
&:= \left(-\partial_r^2-2Q(r)^TQ(r)'\partial_r-Q(r)^TQ''(r)
    +\frac{1}{r^2}E(r)
    +\frac{r^2}{16}-\frac34
    +\widetilde M(r)\right)u(r),
\end{align*}

where 
\begin{align*}
    \widetilde M(r)=\begin{pmatrix}
      [M[r]]_1^1  &0 &  [M[r]]_1^3 
        \\
        0  &[M[r]]_2^2+2 [M[r]]_2^3 & 0
        \\
        [M[r]]_1^3 & 0 & [M[r]]_3^3+\frac12 [M[r]]_2^3
    \end{pmatrix}.
\end{align*} Then, as $[M[r]]_2^3$ is negative, a direct computation shows
\begin{align*}
    (Mu,u)-(\tilde Mu,u)\geq 0. 
\end{align*}
Hence, we see that 
\begin{align*}
    (\Qf^T\Lf_{(1,0)}^{\mathrm{even}} \Qf u,u)\geq (\tilde{\Lf}_{(1,0)}^{\mathrm{even}}u,u)
\end{align*}
 and it suffices to study  $\tilde{\Lf}_{(1,0)}^{\mathrm{even}}u$.
The usefulness of this stems from the fact we now only have to estimate a scalar-valued operator and an operator that is exactly of the $2 \times 2$ form with drift as discussed above. 
 The Schrödinger operator associated to the second component is given by 
 \[
\Lf_{\mathrm{dec}}u(r)
:=-\partial_r^2 u(r)+\frac{12}{r^2}u(r)+V_{\mathrm{dec}}(r) u(r),
\qquad
V_{\mathrm{dec}}(r)
:=\frac{r^2}{16}-\frac34+M_{22}(r)+2M_{23}(r).
\]
 Since $V_{\mathrm{dec}}(r)\geq0$ for $r\geq4.5$, we obtain
\[
\int_0^\infty r^7
  \bigl[(V_{\mathrm{dec}}(r))_-\bigr]^4\,dr
\leq
\int_0^{4.5} r^7 V_{\mathrm{dec}}(r)^4\,dr
<50000<K(12,4).
\]
 Thus, the GGMT bound implies that the spectrum of the self-adjoint closure of $\Lf_{\mathrm{dec}}$ denoted by $\mathcal{L}_{\mathrm{dec}}$ satisfies $\sigma(\mathcal{L}_{\mathrm{dec}})\subset [0,\infty)$ and $0$ is not an eigenvalue of that operator. So, only the remaining $2 \times 2$ operator given by
 \begin{align*}
     \widehat \Lf_{(1,0)}^{\mathrm{even}} u(r) :=\left( -\partial_r^2- P_1(r)\partial_r -P_2(r)+ \frac{r^2}{16}-\frac34+ \frac{1}{r^2}\widehat E(r)+ \widehat M(r) \right) u(r),
 \end{align*} 
 where
 \begin{align*}
P_1(r)&=\begin{pmatrix}
        0 &2(g_1, g_3')
        \\
        -2(g_1, g_3') &0 
    \end{pmatrix}
    \\
P_2(r)&= \begin{pmatrix} -\|g_1'(r)\|_E^2 & (g_1,g_3'')
        \\
        (g_3,g_1'') &-\|g_3'(r)\|_E^2
    \end{pmatrix}
    \\
\widehat E(r)&=\begin{pmatrix}
    2 &0
    \\
    0 &12
\end{pmatrix}+
\frac{48(5-2\sqrt6)r^2}{5\mu(r)}
\begin{pmatrix}
2r^2
&
40-10\sqrt6+\sqrt6\,r^2
\\
40-10\sqrt6+\sqrt6\,r^2
&
-2r^2
\end{pmatrix}
\\
    \widehat M(r)&=\begin{pmatrix}
        [M[r]]_1^1 &  [M[r]]_1^3
        \\
        \\
        [M[r]]_1^3  & [M[r]]_3^3+\frac12 [M[r]]_2^3.
    \end{pmatrix}
 \end{align*}
needs to be estimated. To get this operator into the canonical ground state eigenvalue form, we shift it by $\frac{1}{2}$, i.e.~ we study $\Lf_{\frac
12}:=  \widehat \Lf_{(1,0)}^{\mathrm{even}} +\frac12.$
\begin{lem}
    The second eigenvalue of the operator $\Lf_\frac12$, denoted by $\lambda_2$, satisfies
    \begin{align*}
        \lambda_2>\frac12.
    \end{align*}
\end{lem}
\begin{proof}
    Naturally, our goal is to apply Prop.~\ref{prop:abstract spectral prop}. To that end, we first note that the SUSY partner of the first diagonal entry is given by
    \begin{align*}
    \tilde{V}(r)= \frac{6}{r^2}+\frac{r^2}{16}+\frac{1}{4}+\frac{8p(r)}{q(r)}
    \end{align*}
    where
\begin{align*}
p(r)&=
-706586480+288462720\sqrt6
+\bigl(135823376-55449664\sqrt6\bigr)r^2
\\
&\quad
+\bigl(-12723496+5194344\sqrt6\bigr)r^4
+\bigl(658804-268956\sqrt6\bigr)r^6
\\
& \quad
+\bigl(-18239+7446\sqrt6\bigr)r^8
+\bigl(218-89\sqrt6\bigr)r^{10}.
\end{align*}
    and
    \begin{align*}
    q(r)&=\bigl(412-168\sqrt6+4(-7+3\sqrt6)r^2+r^4\bigr)
    \\
    &\quad \times (-2060+840\sqrt6 +(156-64\sqrt6)r^2 +(-5+2\sqrt6)r^4\bigr)^2.
    \end{align*}
    Now, one computes that the "effective'' potential 
    $$
    V_{\text{eff}}(r)=\frac{r^2}{16}+\frac{1}{4}+\frac{8p(r)}{q(r)}
    $$
    is such that $ V_{\text{eff}}\geq \frac65$  for $r\geq 4.36$ and 
    \begin{align*}
    \int_0^\infty r^7\left[(V_{\text{eff}}(r)-\frac65)_{-}\right]^4 dr \leq   \int_0^{4.36} r^7(V_{\text{eff}}(r)-\frac65)^4 dr \leq 4800<  5000\leq  \frac{5^{7} 4^4 \Gamma(4)^2}{27 \Gamma(8)}
    \end{align*}
    and Theorem \ref{lem: scalar GGMT} implies that $\alpha_1\geq \frac65$. To estimate $\alpha_2$, set $\omega(r):=(g_1(r),g_3'(r))$.
Since $g_3'=\omega g_1$ and the drift coefficient in
Proposition \ref{prop:abstract spectral prop} is $g=-2\omega$,
the scalar operator defining $\alpha_2$, shifted by $-\frac7{10}$,
is
\[
\Lf_{\alpha_2}-\frac7{10}
=-\partial_r^2+\frac{12}{r^2}+V_{\alpha_2}(r),
\]
where
\[
V_{\alpha_2}(r):=\frac{r^2}{16}-\frac{19}{20}
      +\frac{E_{33}(r)-12}{r^2}
      +M_{33}(r)+\frac12 M_{23}(r)-3\omega(r)^2.
\]
For $p=6$, we obtain
\[
\int_0^\infty r^{11}\left(V_{\alpha_2}(r)_-\right)^6\,dr
<9842000
<K(12,6).
\]
Thus Theorem \ref{lem: scalar GGMT} gives
$\alpha_2\geq\frac7{10}$. Moreover,
\[
\beta=4\|\omega\|_{L^\infty(\mathbb R_+)}^2\leq0.1652.
\]
Consequently,
\[
\beta\leq0.1652
<\frac7{25}
=2\left(\frac65-\frac12\right)
  \left(\frac7{10}-\frac12\right)
\leq2\left(\alpha_1-\frac12\right)
     \left(\alpha_2-\frac12\right).
\]
Applying Proposition~\ref{prop:abstract spectral prop}
with $c=\frac12$ yields $\lambda_2>\frac12$.
\end{proof}
Translating this result back to the operator $\Lf_{(1,0)}^{\mathrm{even}}$, we arrive at the following.
\begin{prop}
The operator  $(\Lf_{(1,0)}^{\mathrm{even}},C^\infty_c(\R_+,\R^3))$ is essentially self-adjoint on $L^2(\R_+,\R^3)$ and its closure $\mathcal{L}_{(1,0)}^{\mathrm{even}}$ satisfies
\begin{align*}
\sigma (\mathcal{L}_{(1,0)}^{\mathrm{even}})\setminus\{-\frac12\}\subset (0,\infty).
\end{align*}
 Furthermore, the eigenvalue $-\frac12$ is simple.
\end{prop}

This concludes our discussion of the low angular momentum operators.
\section{Modes of high angular momentum}
\label{sec:high_modes}

To attack the modes in representations of large highest weight $\lambda$, it is convenient to choose a different orthonormal basis of $V_\lambda$ than that provided by the decomposition of $C^\infty(\mathbb S^4) \otimes \sov \otimes T^\ast_0 \mathbb R^5$ into $\mathbf A^{(\ell)}_\lambda, \mathbf B^{(\ell)}_\lambda$ and $\mathbf C^{(\ell)}_\lambda$, $\ell \in \mathbb N$. Recall that this decomposition arises from splitting $\sov \otimes T^\ast_0 \mathbb R^5$ into its three irreducible components and taking the tensor product of each component with $\mathcal H_\ell$. Instead, we will take the tensor product $\mathcal K^{(\ell)} := \mathcal H_\ell \otimes \sov$ first, decompose it into irreducibles $\mathcal K^{(\ell)}_{\mu}$ and finally take the tensor product of each with $T_0^\ast \mathbb R^5$.

\subsection{The first order term}

We first consider the operator $\mathbf T_1$ given by $\mathbf T_1 \Psi = [\sigma_i,\partial_i \Psi]$. Since this operator acts separately on each entry $\Psi_j$ of the form $\Psi$, the spaces $\mathcal K^{(\ell)}_\lambda \otimes T^\ast_0 \mathbb R^5$ are invariant under $\mathbf T_1$. To obtain the matrix coefficients of $\mathbf T_1$ in the respective irreducible components of $\mathcal K^{(\ell)}_\lambda \otimes T^\ast_0 \mathbb R^5$, which by Schur's lemma are invariant under $\mathbf T_1$, it suffices to analyze the action of the operator $V \to [\sigma_i,\partial_i V]$ on $V \in \mathcal K^\ell_{\lambda}$.

Recall that a matrix $X \in \sov$ acts on $v \in T^\ast_0 \mathbb R^5$ by $X v$, on $Y \in \sov$ by $[X,Y]$ and on $h \in \mathcal H_\ell$ by $X^{jk} x_k \partial_j$.
Consider the Casimir operator on the representation $\mathcal K^\ell = \sov \otimes \mathcal H_\ell$.
\begin{align*}
    C_{\mathcal K^\ell}V &= \sum_{\mu<\nu} [X(\mu,\nu),[X(\mu,\nu),V]] + \sum_{\mu < \nu} X(\mu,\nu)^{jk} x_k \partial_j( X(\mu,\nu)^{\ell m} x_{\ell} \partial_m V) \\ &- 2 \sum_{\mu<\nu} [X(\mu,\nu)^{j}_k x^k X(\mu,\nu), \partial_j V] 
\end{align*}
Up to sign, the first two summands are Casimir operators of $\sov \cong V_{(1,1)}$ and $\mathcal H_\ell \cong V_{(\ell,0)}$, and therefore coincide with $- 6 \cdot\mathrm{id}$ and $\Delta_{\mathcal S^{4}}|_{\mathcal H_\ell} \cong - \ell(\ell+3) \cdot \mathrm{id}$, respectively. The cross-terms simplify to 
\begin{align*}
    \sum_{\mu<\nu} [X(\mu,\nu)^{jk} x_k X(\mu,\nu), \partial_j V] &= \sum_{\mu<\nu} [x_\nu X(\mu,\nu), \partial_\mu V] - \sum_{\mu<\nu} [x_\mu X(\mu,\nu), \partial_\nu V] \\
    &= \sum_{\mu\neq \nu} [x_\nu X(\mu,\nu), \partial_\mu V] = -[\sigma_\mu,\partial_\mu V].
\end{align*}
The representation $\mathcal K_\ell$ splits into irreducibles as follows:
\begin{align*}
    \mathcal K_\ell = \begin{cases}
        V_{(1,1)} & \ell = 0, \\
        V_{(1,0)} \oplus V_{(1,1)} \oplus V_{(2,1)} & \ell = 1, \\
        V_{(\ell,0)} \oplus V_{(\ell-1,1)} \oplus V_{(\ell,1)} \oplus V_{(\ell+1,1)} & \ell \geq 2.
    \end{cases}
\end{align*}
By Freudenthal's formula, the operator $- C_{\mathcal K^{(\ell)}}$ acts as the scalar $\left\langle (\lambda_1,\lambda_2) , (\lambda_1,\lambda_2) + (3,1) \right\rangle$ on the representation $V_{\lambda}$. In the special cases above, we obtain
\begin{align*}
    - C_{\mathcal K^{(\ell)}}|_{V_{\lambda}} = \begin{cases}
        \ell(\ell+3)=\ell^2+3\ell & \lambda = (\ell,0), \\
        (\ell-1)(\ell+2)+2=\ell^2+\ell & \lambda = (\ell-1,1), \\
        \ell(\ell+3) + 2 = \ell^2 + 3\ell + 2& \lambda = (\ell,1), \\
        (\ell+1)(\ell+4) + 2 = \ell^2 + 5\ell + 6& \lambda = (\ell+1,1). \\
    \end{cases}
\end{align*}

Hence, the operator $V \mapsto [\sigma_j,\partial_j V]$ acts as the following scalars on the irreducible components $\mathcal K_\mu^{(\ell)}$ of $\mathcal K^{(\ell)}$:
\begin{align*}
        [\sigma_j,\partial_j (\cdot)]|_{V^{(\ell)}_{\lambda}} = \begin{cases}
        3 & \lambda = (\ell,0), \\
        \ell + 3 & \lambda = (\ell-1,1), \\
        2 & \lambda = (\ell,1), \\
        -\ell & \lambda = (\ell+1,1). \\
    \end{cases}
\end{align*}
This operator therefore does not have a definite sign, and it must be dealt with carefully.

\subsection{The double commutator}

We turn our attention to the operator $\Psi \mapsto [\sigma_i,[\sigma^i,\Psi]]$. As in the case of $\mathbf T_1$, it acts separately on each coefficient $\Psi_j$ of $\Psi$, but it also does not contain any derivatives. To obtain bounds on $\mathbf T_2$, it thus suffices to analyze the action of $V \mapsto [\sigma_i|_x,[\sigma^i|_x,V]]$ on $\sov$ at any point $x \in \mathbb S^4$, say, at $x = (1,0,0,0,0)$. At this point, $\sigma_j$ specializes to $M_{1,j} := \delta_{j,1} - \delta_{1,j}$. Bounding the operator $X \mapsto [M_{1,j},[M_{1,j},X]]$ on $\sov$ reduces to computing the eigenvalues of a $10\times 10$ symmetric matrix. It turns out that the spectrum of this matrix is $\{-3,-2\}$, and hence, 
\begin{align*}
    - 3\|X\|^2 \leq \left\langle [M_{1,j},[M_{1,j},X]], X \right\rangle \leq -2 \|X\|^2
\end{align*}

\subsection{The curvature commutator}

The last operator to consider is the commutator $\Psi_j \to [\Psi_i, \partial_iW_j-\partial_j W_i+[W^i,W_j]]$. Since it is $SO(5)$-invariant and acts on forms $\Psi \in \Gamma( \sov \otimes T^\ast \mathbb R^5)$ by pointwise multiplication, it suffices to bound it at a single point on each orbit, say $p = (r,0,0,0,0)$, to get a bound on its norm. Writing $W = f(r^2) \sigma$, we obtain

\begin{align*}
    &[\Psi_i, \partial_iW_j-\partial_j W_i+[W^i,W_j]] \\
    = 
     \ &f(r^2) [\Psi_i, \partial_i \sigma_j - \partial_j \sigma_i] + 
    2 f'(r^2) r^2 [\Psi_i, r^{-2} \left( x^i \sigma_j - x^j \sigma_i \right)] + f(r^2)^2 r^2 [\Psi_i, r^{-2}[\sigma_i,\sigma_j]] \\
    =: \ & f(r^2) (T_3 \Psi)_j + 
    2 r^2 f'(r^2) (T_4 \Psi)_j + r^2 f(r^2)^2(T_5 \Psi)_j.
\end{align*}

The operators $T_3$, $T_4$ and $T_5$ are homogeneous of order zero, hence it suffices to consider their action at $p = (1,0,0,0,0)$. Since $\sov \otimes T^\ast_{p} \mathbb R^5$ is a vector space of dimension $50$, one may carry out this computation in any computer algebra system. We used Wolfram Mathematica, see the accompanying notebook \texttt{computations\_high\_momenta.nb}. One obtains 
\begin{align*}
    \mathrm{spec}(T_3|_p) &= \{6,4,-2\}, \\
    \mathrm{spec}(T_4|_p) &= \{-\sqrt{3},-\sqrt{2},0,\sqrt{2},\sqrt{3}\}, \\
    \mathrm{spec}(T_5|_p) &= \{-3,-2,-1,0,1\},
\end{align*}
not counting multiplicity. Note that $f(r^2) > 0$ and $f'(r^2) < 0$. It follows that
\begin{align*}
    \left\langle [\Psi_i, F_{ij}] , \Psi_j \right\rangle|_p \leq \left(6 f(r^2) - 2 \sqrt{3} r^2 f'(r^2) + r^2 f(r^2)^2\right) \| \Psi|_p \|^2
\end{align*}
pointwise at each $p \in \mathbb R^5$.

\subsection{Combining all estimates}

We now collate the results of the previous subsections to obtain a bound from below for the operators $\mathbf L_\lambda$ acting on the spaces $C^\infty([0,\infty)) \otimes E_\lambda$.

To do so, we gather the irreducible components corresponding to highest weight $\lambda$ in the tensor product representations $\mathcal K_\mu^{(\ell)} \otimes V_{(1,0)}$, where $\ell$ ranges over the nonnegative integers and $\mu$ ranges over the highest weights occurring in $\mathcal K^{(\ell)} = \sov \otimes \mathcal H_\ell$. By \Cref{lem:tensor_weigths}, the highest weights occurring in $\mathcal K^{(\ell)}$ are among $(\ell-1,1),(\ell,0),(\ell,1)$ and $(\ell+1,1)$. We employed the Mathematica package LieART \cite{LieART} to verify that each of these weights indeed occurs with multiplicity one. Alternatively, this can be verified in a straightforward computation via the Weyl dimension formula.
Picking a highest weight vector in each irreducible component of $\mathcal K_\mu^{(\ell)} \otimes V_{(1,0)}$ of highest weight $\lambda$ yields a basis of $E_{\lambda,\lambda}$ in which $\mathbf T_1$ is represented by a diagonal matrix, since the spaces $\mathcal K_\mu^{(\ell)} \otimes V_{(1,0)}$ are invariant under $\mathbf T_1$, and on each space, $\mathbf T_1$ acts as a multiple of the identity.

The operators $\mathbf T_2,\mathbf T_3,\mathbf T_4$ and $\mathbf T_5$ are multiplication operators and can be estimated pointwise from above by multiples of the identity
\begin{align*}
    \mathbf T_2 \leq - 2 \, \mathrm{id}, \ \mathbf T_3 \leq 6 \, \mathrm{id}, \ -\sqrt{3} \, \mathrm{id} \leq \mathbf T_4 \leq \sqrt{3} \, \mathrm{id}, \ \mathbf T_5 \leq \mathrm{id}.
\end{align*}
It follows that $-\mathbf L$ can be bounded from below by a diagonal operator, all of whose entries are of the form $-\partial_r^2 + V_{\ell,\mu}$, where 
\begin{align}
\begin{split}
    \label{eq:high_mode_estimate}
    V_{\ell,\mu} = \ &\frac{2+\ell(\ell+3)}{r^2} + \tfrac1{16}(r^2 - 12) + \begin{rcases}\begin{dcases} 2 \ell & \mu = (\ell+1,1) \\ -4 & \mu = (\ell,1) \\ -2(\ell + 3) & \mu = (\ell-1,1) \\ -6 & \mu = (\ell,0) \end{dcases}\end{rcases} f(r^2)  \\ & \ + 2 r^2 f(r^2)^2 - 12 f(r^2) + 4 \sqrt{3} r^2 f'(r^2) - 2 r^2 f(r^2)^2,
\end{split}
\end{align}
and $\ell$, $\mu$ are such that $\mathcal K^{(\ell)}_\mu \otimes V_{(1,0)}$ has a nontrivial irreducible component of highest weight $\lambda$. If $-\partial_r^2 + V_{\ell,\mu}$ is nonnegative for all such $\ell$ and $\mu$, then $-\mathbf L_\lambda$ is bounded below by a nonnegative operator and hence possesses no negative eigenvalues. We next determine all $\mu$ and $\lambda$ where this is the case.

\begin{lemma}
    \label{lem:high_mode_potential}
    With $V_{\ell,\mu}$ defined as in \eqref{eq:high_mode_estimate}, we have that 
    \begin{align*}
        - \partial_r^2 + V_{\ell,\mu} > 0
    \end{align*}
    for all $\ell \in \mathbb N$, $\mu \in \mathbb N^2$ except for
    \begin{enumerate}
        \item $\ell \in \{0,1,2,3\}$,
        \item or $\ell=4$, $\mu \in \{(3,1),(4,0),(4,1)\}$, 
        \item or $\ell=5$, $\mu = (4,1)$.
    \end{enumerate}
\end{lemma}

\begin{proof}
One readily checks that the potentials \(V_{4,(5,1)}\), \(V_{5,(5,1)}\), and \(V_{5,(5,0)}\) are positive for \(r>0\). Since the potentials are increasing in \(\ell\) along each branch for \(\ell\geq4\), it remains to consider \(V_{6,(5,1)}\), which is explicitly given by
\begin{align*}
V_{6,(5,1)}(r)
=
\frac{56}{r^2}+\frac{r^2}{16}-\frac34
-\frac{30}{a r^2+b}
-\frac{3\sqrt2 r^2}
{(a r^2+b)^2}.
\end{align*}
Now, the effective potential
\begin{align*}
V(r):=
\frac{r^2}{16}-\frac34
-\frac{30}{a r^2+b}
-\frac{3\sqrt2 r^2}
{(a r^2+b)^2}
\end{align*}
is positive for $r\geq 6.5$ and one checks that
\[
\int_0^{6.5} r^7
\left(
\frac{r^2}{16}-\frac34
-\frac{30}{ar^2+b}
-\frac{3\sqrt2\,r^2}{(ar^2+b)^2}
\right)^4dr < 11\cdot 10^6< K(56,4).
\]
Thus, the claim follows from Lemma \ref{lem: GGMT vector} and a check with Mathematica can be found in the notebook 
\texttt{computations\_high\_momenta.nb}.
\end{proof}

To apply the previous lemma in our task of estimating $\mathbf L_\lambda$, acting on a component $E_\lambda \in \sov \otimes T_0 \mathbb R^5 \otimes C^\infty(\mathbb S^4)$, recall that $E_\lambda$ collects the irreducible components of highest weight $\lambda$ in the tensor representation $T_0\mathbb R^5 \otimes \mathcal K^{(\ell)}_{\mu}$. By \Cref{lem:tensor_weigths}, and since the list of weights occurring in $T_0 \mathbb R^5 \cong V_{(1,0)}$ is $\{(0,0), \pm (1,0), \pm (0,1)\}$, only weights $\mu \in \lambda \pm \{(0,0), \pm (1,0), \pm (0,1)\}$ must be accounted for.

\begin{proposition}
    The operators $\mathbf L_{\lambda}$ arising from the symmetry reduction of the operator $\mathbf L$ defined in \eqref{Eq: linearised op} are nonnegative for all $\lambda \in \mathbb N^2_{\lambda_1 \geq \lambda_2}$ except possibly
    \begin{align*}
        \lambda \in \{\lambda \in \mathbb N^2_{\lambda_1 \geq \lambda_2}: \lambda_1 \leq 4, \lambda_2 \leq 2\} \cup \{(5,0),(5,1)\}.
    \end{align*}
\end{proposition}

\begin{proof}
    First, assume that $\lambda_1 \geq 6$. By \Cref{lem:tensor_weigths}, $\mathcal K^{(\ell)}_\mu$ has no irreducible components with highest weight $\lambda$ unless $\mu_1 \geq 5$. Now \Cref{lem:high_mode_potential} implies $- \partial_r^2 + V_{\ell,\mu}$ is a nonnegative operator for all such $\mu,\ell$, regardless of the degree $\ell$, so $\mathbf L_\lambda$ is nonnegative.

    Now let $\lambda = (5,2)$. Then the only admissible values of $\mu$ are in $\{(4,2),(5,1),(6,2)\}$. Again by \Cref{lem:high_mode_potential}, the operator $- \partial_r^2 + V_{\ell,\mu}$ is positive for these values of $\mu$.
\end{proof}

For the remaining weights $\lambda$, negative eigenvalues of $\mathbf L_\lambda$ cannot be excluded by the crude arguments of this section. Instead, one must compute the operators $\mathbf L_\lambda$ explicitly, following the procedure laid out in \Cref{sec:symmetry_reduction_specific}, and analyze each resulting Schr\"odinger operator separately. This analysis is contained in \Cref{section: low angulars} of this article.

\appendix
\section{On the self-adjointness of $\Lf$}

In this appendix, we establish the self-adjointness of $\Lf$. For this, we let $\rho(x)=e^{-\frac{|x|^2}{4}}$ and start with a useful equivalence of norms. 
\begin{lem}
    We have that 
    \begin{align*}
     H^2_\rho(\R^5,\Omega(\mathfrak{so}(5)))
     \subset   L^2_{|.|^4\rho}(\R^5,\Omega(\mathfrak{so}(5))).
    \end{align*}
    Consequently, 
    \begin{align*}
        H^2_\rho(\R^5,\Omega(\mathfrak{so}(5)))\cap L^2_{|.|^4\rho}(\R^5,\Omega(\mathfrak{so}(5)))=H^2_\rho(\R^5,\Omega(\mathfrak{so}(5))).
    \end{align*}
\end{lem}
 \begin{proof}
     To see this, let $u\in C^\infty_c(\R^5)$ and note that 
\begin{align*}
    \int_{\R^5} x_i^4 |u|^2(x) \rho(x) dx&= -2\int_{\R^5} x_i^3 |u|^2(x) \partial_i\rho(x) dx
    \\
    &=6\int_{\R^5} x_i^2 |u(x)|^2 \rho(x) dx+4\Re\int_{\R^5} x_i^3 u(x)\partial_i\overline{u(x)} \rho(x) dx.
\end{align*}
Thus, the Cauchy-Schwarz inequality implies 
\begin{align*}
    \|u\|_{L^2_{|.|^4\rho}(\R^5)}&\lesssim   \|u\|_{H^1_{|.|^2\rho}(\R^5)}
\end{align*}
and another integration by parts then shows
\begin{align*}
    \|u\|_{H^1_{|.|^2\rho}(\R^5)}\lesssim \|u\|_{H^2_{\rho}(\R^5)}.
\end{align*}
Hence, by density it follows that $ H^2_\rho(\R^5)\subset L^2_{|.|^4\rho}(\R^5) $ which implies the claim. 
 \end{proof}
 Next, we recall the definition of $\Lf$
\begin{equation*}
        (\Lf U)_j := \Delta U_j -\frac12(\Lambda+1) U_j+2[W^i, \partial_iU_j]+[W^i,[W_i,U_j]]+2[U^i,F_{ij}^W]
\end{equation*}
and
$D(\Lf)= H^2_\rho(\R^5,\Omega(\mathfrak{so}(5)))$.
\begin{lem}\label{lem: self-adjoint}
    The operator $(\Lf,D(\Lf))$ is self-adjoint on $L^2_\rho(\R^5,\Omega (\mathfrak{so}(5)))$. Furthermore, its spectrum is discrete. 
\end{lem}
\begin{proof}
We first study the operator 
$$
    (L_0,D(L_0))=\left(-\Delta +\frac{|x|^2}{16},C^\infty_c(\R^5)\right)
    $$
    According to Theorem X.28 in \cite{ReeSim75}  this operator
    is essentially self-adjoint on $L^2(\R^5)$. 
    To compute its self-adjoint closure, we note that 
    \begin{align*}
    \|L_0u\|_{L^2(\R^5)}^2&= \|\Delta u\|_{L^2(\R^5)}^2+\frac{1}{16^2}\| |.|^2u\|_{L^2(\R^5)}^2-\frac{1}{8}\Re\int_{\R^5} (\Delta u)(x) |x|^2 \overline{u(x)} dx
    \\
    &= \|\Delta u\|_{L^2(\R^5)}^2+\frac{1}{16^2}\||.|^2u\|_{L^2(\R^5)}^2+ \frac{1}{8}\Re\int_{\R^5} (\nabla u)(x) \cdot \overline{\nabla (|.|^2 u)(x) }dx
    \\
    &= \|\Delta u\|_{L^2(\R^5)}^2+\frac{1}{16^2}\||.|^2u\|_{L^2(\R^5)}^2+\frac{1}{8}\||.|\nabla u\|_{L^2(\R^5)}^2- \frac{5}{8}\|u\|_{L^2(\R^5)}^2.
   \end{align*}
   So, the graph norm induced by $L_0$ is equivalent to the norm 
   $$
\|u\|:=\|u\|_{H^2(\R^5)}+\||.|^2 u\|_{L^2(\R^5)}
   $$
   and closing the operator yields that
   $$
\left(\overline{L_0}, H^2(\R^5)\cap L^2(\R^5, |x|^4 dx)\right)
   $$
   is self-adjoint. 
    Thus, also the operator 
    $$
(L,D(L))=\left(\overline{L_0}-\frac34, H^2(\R^5)\cap L^2(\R^5, |x|^4 dx)\right)
    $$
    is self-adjoint.
Consider now the unitary transformation 
\[T:L^2(\R^5)\to L^2_\rho(\R^5),\qquad   Tu= e^{\frac{|\cdot|^2}{8}}u\]
and set $(L_c, D(L_c))= (T L T^{-1}, T(D(L))).$ Then, $(L_c, D(L_c))$ is self-adjoint on $L^2_\rho(\R^5)$,
where 
\[
T(D(L))=H^2_\rho(\R^5)\cap L^2_{|.|^4\rho}(\R^5)=H^2_\rho(\R^5).
\]
Furthermore,     
\[
L_c u(x)= T LT^{-1} u(x)=-\Delta u+\frac12 u+\frac12x\cdot \nabla u. 
\]
Now, note that 
\begin{align*}
    (-\Lf U)_j &=L_c U_j-2[W^i, \partial_iU_j]-[W^i,[W_i,U_j]]-2[U^i,F_{ij}^W]
\end{align*}
and define 
\begin{align*}
    (\Lf' U)_j := 2[W^i, \partial_iU_j]+[W^i,[W_i,U_j]]+2[U^i,F_{ij}^W]
\end{align*}
Then, since $\Lf'$ only contains first order derivatives and given that $W_i$ is smooth and decaying linearly, one readily computes that
$\Lf'$ has zero relative $\Lf_0$ bound where 
$$
(\Lf_0 U)_j=L_c U_j.
$$
Lastly, we show that $\Lf$ is symmetric on $D(\Lf)$. 
Consider $D_W^i:=\partial^i+[W^i,\cdot]$. Then for $X,Y\in C^\infty_c(\R^5, \Omega( \mathfrak{so}(5)))$ we have that 
\begin{align*}
    \int_{\R^5}\tr\left((D_W^i X)^j(x) Y_j(x)^*\right) e^{-\frac{|x|^2}{4}}dx&=   \int_{\R^5}   \left[\tr (\partial^i X^j(x) Y_j(x)^*)+\tr\left([W_i,X^j](x) Y_j(x)^*\right) \right]e^{-\frac{|x|^2}{4}}dx
\\
&= \int_{\R^5}  \left[-\tr(X^j(x) \partial^i Y_j(x)^*)+\frac{x_i}{2}\tr( X^j(x)Y_j(x)^*)\right]\ e^{-\frac{|x|^2}{4}}dx
\\
&\quad - \int_{\R^5} \tr (X^j(x) [W_i,Y_j](x)^*) e^{-\frac{|x|^2}{4}}dx.
\end{align*}
So, we see that, formally, the adjoint of $D_W$, denoted by $D^{W*}$, is given by
\[ D^{W*}_i= -D_W^i+\frac{1}{2}x_i. \]
By using this, and the identities  $\partial^i W_i=x^i W_i=0,$
we can recast 
\begin{align*}
    -(\Lf U)_j=D^{W*}_iD^i_W U_j +\frac{1}{2} U_j-2[U^i , F^W_{ij}]
\end{align*}
for $U\in C^\infty_c(\R^5,\Omega( \mathfrak{so}(5)))$.
Consequently, for $ U,V\in C^\infty_c(\R^5,\Omega( \mathfrak{so}(5)))$
we have that
$$
(-\Lf U,V)_{L^2_\rho}=(U,-\Lf V)_{L^2_\rho}.
$$
Thus, as $C^\infty_c(\R^5, \Omega(\mathfrak{so}(5)))$ is dense in $D(\Lf)$ with respect to the graph norm, $\Lf$ is symmetric on $D(\Lf).$ Consequently, as  $\Lf'$ has zero relative $\Lf_0$ bound, we conclude that 
\[\left(\Lf, H^2_\rho(\R^5,\Omega(\mathfrak{so}(5)))\right) \]is self-adjoint.
Lastly, the discreteness of the spectrum of $-\Lf$ follows from the compactness of the resolvent of $\Lf$ and was already established in Lemma 2.1 of \cite{KimLee24}.
\end{proof}

\section{The low angular momentum matrices}
For completeness, we collect all of the matrix potentials of the operators studied in \Cref{section: low angulars}. To this end, we recall that in canonical form, we are dealing with operators of the form
\begin{flalign*}
\Lf_{(i,j)}^{\mathrm{par}}   u(r)&:= -\partial_r^2 u(r)+ \left[-\frac34 +\frac{r^2}{16}+ \frac{1}{r^2}S_{(i,j)}^{\mathrm{par}}-\frac{1}{(a r^2+b)^2}\left( r^2 A_{(i,j)}^{\mathrm{par}}+B_{(i,j)}^{\mathrm{par}}
\right)\right] u(r), &
\end{flalign*}
where $\mathrm{par}$ denotes the parity of the component's coefficients, either even or odd.

\subsubsection*{The $(0,0)$-component}

\begin{flalign*}
A_{(0,0)}^{\mathrm{odd}} &= -9 + \frac{9 \sqrt{6}}{2}, \
B_{(0,0)}^{\mathrm{odd}} = 162 - 63 \sqrt{6}, \
S_{(0,0)}^{\mathrm{odd}} = 6. &
\end{flalign*}
\subsubsection*{The $(1,0)$-component}

\begin{flalign*}
A_{(1,0)}^{\mathrm{even}} &= \left(\begin{smallmatrix}-6 + \frac{9 \sqrt{6}}{5} & - \frac{3}{2} + \frac{3 \sqrt{6}}{5} & - \frac{3 \sqrt{15}}{10} + \frac{3 \sqrt{10}}{5}\\- \frac{3}{2} + \frac{3 \sqrt{6}}{5} & - \frac{33}{4} + \frac{183 \sqrt{6}}{40} & - \frac{27 \sqrt{10}}{40} + \frac{3 \sqrt{15}}{20}\\- \frac{3 \sqrt{15}}{10} + \frac{3 \sqrt{10}}{5} & - \frac{27 \sqrt{10}}{40} + \frac{3 \sqrt{15}}{20} & - \frac{15}{4} + \frac{21 \sqrt{6}}{8}\end{smallmatrix}\right), & \\
B_{(1,0)}^{\mathrm{even}} &= \left(\begin{smallmatrix}108 - 42 \sqrt{6} & 0 & 0\\0 & \frac{351}{2} - \frac{273 \sqrt{6}}{4} & - \frac{9 \sqrt{15}}{2} + \frac{21 \sqrt{10}}{4}\\0 & - \frac{9 \sqrt{15}}{2} + \frac{21 \sqrt{10}}{4} & \frac{81}{2} - \frac{63 \sqrt{6}}{4}\end{smallmatrix}\right), & \\
S_{(1,0)}^{\mathrm{even}} &= \mathrm{diag}(2,12,12). &
\end{flalign*}
\begin{flalign*}
A_{(1,0)}^{\mathrm{odd}} &= \begin{matrix}-6 + 3 \sqrt{6}\end{matrix}, \
B_{(1,0)}^{\mathrm{odd}} = \begin{matrix}108 - 42 \sqrt{6}\end{matrix}, \
S_{(1,0)}^{\mathrm{odd}} = 6. &
\end{flalign*}

\subsubsection*{The $(1,1)$-component}

\begin{flalign*}
A_{(1,1)}^{\mathrm{even}} &= \left(\begin{smallmatrix}- \frac{24}{5} + \frac{6 \sqrt{6}}{5} & - \frac{8}{5} + \frac{2 \sqrt{6}}{5} & - \frac{\sqrt{30}}{5} + \frac{2 \sqrt{5}}{5}\\- \frac{8}{5} + \frac{2 \sqrt{6}}{5} & - \frac{26}{5} + \frac{47 \sqrt{6}}{15} & - \frac{2 \sqrt{5}}{5} + \frac{7 \sqrt{30}}{15}\\- \frac{\sqrt{30}}{5} + \frac{2 \sqrt{5}}{5} & - \frac{2 \sqrt{5}}{5} + \frac{7 \sqrt{30}}{15} & -3 + \frac{13 \sqrt{6}}{6}\end{smallmatrix}\right), & \\ 
B_{(1,1)}^{\mathrm{even}} &= \left(\begin{smallmatrix}72 - 28 \sqrt{6} & 0 & 0\\0 & 132 - \frac{154 \sqrt{6}}{3} & - \frac{14 \sqrt{30}}{3} + 12 \sqrt{5}\\0 & - \frac{14 \sqrt{30}}{3} + 12 \sqrt{5} & 30 - \frac{35 \sqrt{6}}{3}\end{smallmatrix}\right), & \\ 
S_{(1,1)}^{\mathrm{even}} &= \mathrm{diag}(2,12,12). &
\end{flalign*}

\begin{flalign*}
A_{(1,1)}^{\mathrm{odd}} &= \left(\begin{smallmatrix}- \frac{21}{4} + \frac{15 \sqrt{6}}{8} & \sqrt{3} & - \frac{9 \sqrt{7}}{28} + \frac{19 \sqrt{42}}{56} & - \frac{3 \sqrt{14}}{7} + \frac{3 \sqrt{21}}{14}\\\sqrt{3} & -4 + \frac{4 \sqrt{6}}{3} & - \frac{11 \sqrt{21}}{21} + \frac{2 \sqrt{14}}{7} & - \frac{2 \sqrt{7}}{7} + \frac{\sqrt{42}}{7}\\- \frac{9 \sqrt{7}}{28} + \frac{19 \sqrt{42}}{56} & - \frac{11 \sqrt{21}}{21} + \frac{2 \sqrt{14}}{7} & - \frac{93}{28} + \frac{229 \sqrt{6}}{168} & - \frac{\sqrt{2}}{7} + \frac{3 \sqrt{3}}{14}\\- \frac{3 \sqrt{14}}{7} + \frac{3 \sqrt{21}}{14} & - \frac{2 \sqrt{7}}{7} + \frac{\sqrt{42}}{7} & - \frac{\sqrt{2}}{7} + \frac{3 \sqrt{3}}{14} & - \frac{24}{7} + \frac{24 \sqrt{6}}{7}\end{smallmatrix}\right), & \\ 
B_{(1,1)}^{\mathrm{odd}} &= \left(\begin{smallmatrix}\frac{243}{2} - \frac{189 \sqrt{6}}{4} & - 14 \sqrt{3} + 18 \sqrt{2} & - \frac{7 \sqrt{42}}{4} + \frac{9 \sqrt{7}}{2} & 0\\- 14 \sqrt{3} + 18 \sqrt{2} & 96 - \frac{112 \sqrt{6}}{3} & - 6 \sqrt{14} + \frac{14 \sqrt{21}}{3} & 0\\- \frac{7 \sqrt{42}}{4} + \frac{9 \sqrt{7}}{2} & - 6 \sqrt{14} + \frac{14 \sqrt{21}}{3} & - \frac{3}{2} + \frac{7 \sqrt{6}}{12} & 0\\0 & 0 & 0 & 72 - 28 \sqrt{6}\end{smallmatrix}\right), & \\
S_{(1,1)}^{\mathrm{odd}} &= \mathrm{diag}(6,6,6,20). &
\end{flalign*}

\subsubsection*{The $(2,0)$-component}

\begin{flalign*}
    A_{(2,0)}^{\mathrm{odd}} &= \left(\begin{smallmatrix}- \frac{177}{28} + \frac{87 \sqrt{6}}{56} & \frac{3 \sqrt{15}}{28} + \frac{9 \sqrt{10}}{56} & - \frac{15 \sqrt{10}}{28} + \frac{3 \sqrt{15}}{7} & - \frac{3 \sqrt{15}}{7} + \frac{9 \sqrt{10}}{28}\\\frac{3 \sqrt{15}}{28} + \frac{9 \sqrt{10}}{56} & - \frac{39}{28} + \frac{33 \sqrt{6}}{56} & \frac{3}{7} - \frac{3 \sqrt{6}}{28} & \frac{9}{7} - \frac{11 \sqrt{6}}{28}\\- \frac{15 \sqrt{10}}{28} + \frac{3 \sqrt{15}}{7} & \frac{3}{7} - \frac{3 \sqrt{6}}{28} & - \frac{111}{14} + \frac{135 \sqrt{6}}{28} & - \frac{9}{14} + \frac{33 \sqrt{6}}{28}\\- \frac{3 \sqrt{15}}{7} + \frac{9 \sqrt{10}}{28} & \frac{9}{7} - \frac{11 \sqrt{6}}{28} & - \frac{9}{14} + \frac{33 \sqrt{6}}{28} & - \frac{47}{14} + \frac{71 \sqrt{6}}{28}\end{smallmatrix}\right), & \\
    B_{(2,0)}^{\mathrm{odd}} &= \left(\begin{smallmatrix}\frac{189}{2} - \frac{147 \sqrt{6}}{4} & - \frac{21 \sqrt{10}}{4} + \frac{9 \sqrt{15}}{2} & 0 & 0\\- \frac{21 \sqrt{10}}{4} + \frac{9 \sqrt{15}}{2} & \frac{27}{2} - \frac{21 \sqrt{6}}{4} & 0 & 0\\0 & 0 & 189 - \frac{147 \sqrt{6}}{2} & 27 - \frac{21 \sqrt{6}}{2}\\0 & 0 & 27 - \frac{21 \sqrt{6}}{2} & 45 - \frac{35 \sqrt{6}}{2}\end{smallmatrix}\right), & \\
    S_{(2,0)}^{\mathrm{odd}} &= \mathrm{diag}(6,6,20,20). &
\end{flalign*}
\begin{flalign*}
  A_{(2,0)}^{\mathrm{even}} = \begin{matrix}-6 + 3 \sqrt{6}\end{matrix},\
    B_{(2,0)}^{\mathrm{even}} &= \begin{matrix}108 - 42 \sqrt{6}\end{matrix},\
    S_{(2,0)}^{\mathrm{even}} = 12. & 
\end{flalign*}

\subsubsection*{The $(2,1)$-component}

\begin{flalign*}
 A_{(2,1)}^{\mathrm{even}} &= \left(\begin{smallmatrix}- \frac{3 \sqrt{6}}{5} - \frac{6}{5} & - \frac{3 \sqrt{70}}{35} + \frac{3 \sqrt{105}}{70} & - \frac{2 \sqrt{70}}{35} + \frac{2 \sqrt{105}}{35} & - \frac{17 \sqrt{35}}{70} + \frac{27 \sqrt{210}}{280} & - \frac{9 \sqrt{21}}{70} + \frac{9 \sqrt{14}}{280} & 0\\- \frac{3 \sqrt{70}}{35} + \frac{3 \sqrt{105}}{70} & - \frac{21}{4} + \frac{15 \sqrt{6}}{8} & - \sqrt{6} & - \frac{9 \sqrt{2}}{16} + \frac{19 \sqrt{3}}{16} & - \frac{\sqrt{30}}{80} + \frac{49 \sqrt{5}}{80} & - \frac{\sqrt{42}}{7} + \frac{4 \sqrt{7}}{7}\\- \frac{2 \sqrt{70}}{35} + \frac{2 \sqrt{105}}{35} & - \sqrt{6} & -4 + \frac{4 \sqrt{6}}{3} & - \frac{\sqrt{2}}{2} + \frac{5 \sqrt{3}}{3} & - \frac{\sqrt{30}}{10} + \frac{\sqrt{5}}{5} & - \frac{2 \sqrt{7}}{7} + \frac{\sqrt{42}}{7}\\- \frac{17 \sqrt{35}}{70} + \frac{27 \sqrt{210}}{280} & - \frac{9 \sqrt{2}}{16} + \frac{19 \sqrt{3}}{16} & - \frac{\sqrt{2}}{2} + \frac{5 \sqrt{3}}{3} & - \frac{69}{32} + \frac{173 \sqrt{6}}{192} & - \frac{29 \sqrt{15}}{160} - \frac{59 \sqrt{10}}{320} & - \frac{\sqrt{21}}{7} + \frac{\sqrt{14}}{7}\\- \frac{9 \sqrt{21}}{70} + \frac{9 \sqrt{14}}{280} & - \frac{\sqrt{30}}{80} + \frac{49 \sqrt{5}}{80} & - \frac{\sqrt{30}}{10} + \frac{\sqrt{5}}{5} & - \frac{29 \sqrt{15}}{160} - \frac{59 \sqrt{10}}{320} & - \frac{1687}{1440} + \frac{1111 \sqrt{6}}{960} & - \frac{\sqrt{210}}{21} + \frac{5 \sqrt{35}}{63}\\0 & - \frac{\sqrt{42}}{7} + \frac{4 \sqrt{7}}{7} & - \frac{2 \sqrt{7}}{7} + \frac{\sqrt{42}}{7} & - \frac{\sqrt{21}}{7} + \frac{\sqrt{14}}{7} & - \frac{\sqrt{210}}{21} + \frac{5 \sqrt{35}}{63} & - \frac{29}{9} + \frac{23 \sqrt{6}}{6}\end{smallmatrix} \right), &\\
    B_{(2,1)}^{\mathrm{even}} &= \left(\begin{smallmatrix}-36 + 14 \sqrt{6} & 0 & 0 & 0 & 0 & 0\\0 & \frac{243}{2} - \frac{189 \sqrt{6}}{4} & -36 + 14 \sqrt{6} & - \frac{49 \sqrt{3}}{8} + \frac{63 \sqrt{2}}{8} & - \frac{63 \sqrt{5}}{8} + \frac{27 \sqrt{30}}{8} & 0\\0 & -36 + 14 \sqrt{6} & 96 - \frac{112 \sqrt{6}}{3} & - \frac{56 \sqrt{3}}{3} + 24 \sqrt{2} & 0 & 0\\0 & - \frac{49 \sqrt{3}}{8} + \frac{63 \sqrt{2}}{8} & - \frac{56 \sqrt{3}}{3} + 24 \sqrt{2} & \frac{3}{16} - \frac{7 \sqrt{6}}{96} & - \frac{81 \sqrt{15}}{16} + \frac{189 \sqrt{10}}{32} & 0\\0 & - \frac{63 \sqrt{5}}{8} + \frac{27 \sqrt{30}}{8} & 0 & - \frac{81 \sqrt{15}}{16} + \frac{189 \sqrt{10}}{32} & \frac{549}{16} - \frac{427 \sqrt{6}}{32} & 0\\0 & 0 & 0 & 0 & 0 & 90 - 35 \sqrt{6}\end{smallmatrix} \right), & \\
    S_{(2,1)}^{\mathrm{even}} &= \mathrm{diag}(2,12,12,12,12,30). &
\end{flalign*}
\begin{flalign*}
A_{(2,1)}^{\mathrm{odd}}
&= \frac{1}{42}
\left(\begin{smallmatrix}
-204+38\sqrt{6}       & 48\sqrt{2}-24\sqrt{3} & 12\sqrt{2}+20\sqrt{3} & 24\sqrt{2}-24\sqrt{3} \\
48\sqrt{2}-24\sqrt{3} & -216+144\sqrt{6}      & 12-6\sqrt{6}          & 36-54\sqrt{6}         \\
12\sqrt{2}+20\sqrt{3} & 12-6\sqrt{6}          & -18+\sqrt{6}          & -48+18\sqrt{6}        \\
24\sqrt{2}-24\sqrt{3} & 36-54\sqrt{6}         & -48+18\sqrt{6}        & -108+90\sqrt{6}
\end{smallmatrix}\right), & \\
B_{(2,1)}^{\mathrm{odd}}
&= \frac{1}{3}
\left(\begin{smallmatrix}
180-70\sqrt{6}        & 0                & 36\sqrt{2}-28\sqrt{3} & 0              \\
0                     & 432-168\sqrt{6}  & 0                     & -108+42\sqrt{6} \\
36\sqrt{2}-28\sqrt{3} & 0                & -18+7\sqrt{6}         & 0              \\
0                     & -108+42\sqrt{6}  & 0                     & 108-42\sqrt{6}
\end{smallmatrix}\right), & \\[1ex]
S_{(2,1)}^{\mathrm{odd}}
&= \operatorname{diag}(6,20,6,20). &
\end{flalign*}

\subsubsection*{The $(2,2)$-component.}

\begin{flalign*}
A_{(2,2)}^{\mathrm{even}} = 0, \
B_{(2,2)}^{\mathrm{even}} &= 0, \
S_{(2,2)}^{\mathrm{even}} = 12. &
\end{flalign*}
\begin{flalign*}
A_{(2,2)}^{\mathrm{odd}} = \left(\begin{smallmatrix}- \frac{3 \sqrt{6}}{2} - \frac{3}{7} & \frac{2 \sqrt{3}}{7}\\\frac{2 \sqrt{3}}{7} & - \frac{4}{7} + 2 \sqrt{6}\end{smallmatrix}\right), 
B_{(2,2)}^{\mathrm{odd}} &= \left(\begin{smallmatrix}-54 + 21 \sqrt{6} & 0\\0 & 72 - 28 \sqrt{6}\end{smallmatrix}\right), \
S_{(2,2)}^{\mathrm{odd}} = \mathrm{diag}(6,20). &
\end{flalign*}

\subsubsection*{The $(3,0)$-component}

\begin{flalign*}
 A_{(3,0)}^{\mathrm{even}} &= \left(\begin{smallmatrix}- \frac{13}{2} + \frac{5 \sqrt{6}}{4} & \frac{1}{2} + \frac{5 \sqrt{6}}{12} & - \frac{5 \sqrt{2}}{4} + \sqrt{3} & - \frac{\sqrt{21}}{3} + \frac{\sqrt{14}}{4}\\\frac{1}{2} + \frac{5 \sqrt{6}}{12} & - \frac{29}{18} + \frac{7 \sqrt{6}}{12} & - \frac{\sqrt{2}}{4} + \frac{\sqrt{3}}{3} & - \frac{11 \sqrt{14}}{36} + \frac{\sqrt{21}}{3}\\- \frac{5 \sqrt{2}}{4} + \sqrt{3} & - \frac{\sqrt{2}}{4} + \frac{\sqrt{3}}{3} & - \frac{31}{4} + \frac{41 \sqrt{6}}{8} & - \frac{\sqrt{7}}{4} + \frac{13 \sqrt{42}}{24}\\- \frac{\sqrt{21}}{3} + \frac{\sqrt{14}}{4} & - \frac{11 \sqrt{14}}{36} + \frac{\sqrt{21}}{3} & - \frac{\sqrt{7}}{4} + \frac{13 \sqrt{42}}{24} & - \frac{113}{36} + \frac{61 \sqrt{6}}{24}\end{smallmatrix}\right), & \\
    B_{(3,0)}^{\mathrm{even}} &= \left(\begin{smallmatrix}81 - \frac{63 \sqrt{6}}{2} & 27 - \frac{21 \sqrt{6}}{2} & 0 & 0\\27 - \frac{21 \sqrt{6}}{2} & 9 - \frac{7 \sqrt{6}}{2} & 0 & 0\\0 & 0 & \frac{405}{2} - \frac{315 \sqrt{6}}{4} & - \frac{21 \sqrt{42}}{4} + \frac{27 \sqrt{7}}{2}\\0 & 0 & - \frac{21 \sqrt{42}}{4} + \frac{27 \sqrt{7}}{2} & \frac{99}{2} - \frac{77 \sqrt{6}}{4}\end{smallmatrix}\right), & \\
    S_{(3,0)}^{\mathrm{even}} &= \mathrm{diag}(12,12,30,30). & 
\end{flalign*}

\begin{flalign*}
    A_{(3,0)}^{\mathrm{odd}} = \begin{matrix}-6 + 3 \sqrt{6}\end{matrix}, \
    B_{(3,0)}^{\mathrm{odd}} &= \begin{matrix}108 - 42 \sqrt{6}\end{matrix}, \
    S_{(3,0)}^{\mathrm{odd}} = 20. & 
\end{flalign*}

\subsubsection*{The $(3,1)$-component}
\begin{flalign*}
    A_{(3,1)}^{odd} &= &
\end{flalign*}
\noindent\scalebox{0.65}{
\begin{minipage}{\textwidth}
\begin{flalign*}
    &\begin{pmatrix}- \frac{13 \sqrt{6}}{14} - \frac{11}{7} & - \frac{4 \sqrt{21}}{21} + \frac{\sqrt{14}}{7} & - \frac{\sqrt{210}}{21} + \frac{2 \sqrt{35}}{21} & - \frac{34 \sqrt{106}}{1113} + \frac{58 \sqrt{159}}{1113} & - \frac{2 \sqrt{17490}}{159} + \frac{4 \sqrt{2915}}{159} & 0\\- \frac{4 \sqrt{21}}{21} + \frac{\sqrt{14}}{7} & - \frac{21}{4} + \frac{15 \sqrt{6}}{8} & \sqrt{10} & - \frac{27 \sqrt{2226}}{2968} - \frac{\sqrt{371}}{56} & - \frac{15 \sqrt{40810}}{4664} + \frac{\sqrt{61215}}{88} & - \frac{\sqrt{2310}}{33} + \frac{\sqrt{385}}{22}\\- \frac{\sqrt{210}}{21} + \frac{2 \sqrt{35}}{21} & \sqrt{10} & -4 + \frac{4 \sqrt{6}}{3} & - \frac{26 \sqrt{5565}}{1113} - \frac{\sqrt{3710}}{371} & - \frac{28 \sqrt{4081}}{583} + \frac{3 \sqrt{24486}}{583} & - \frac{2 \sqrt{154}}{33} + \frac{2 \sqrt{231}}{33}\\- \frac{34 \sqrt{106}}{1113} + \frac{58 \sqrt{159}}{1113} & - \frac{27 \sqrt{2226}}{2968} - \frac{\sqrt{371}}{56} & - \frac{26 \sqrt{5565}}{1113} - \frac{\sqrt{3710}}{371} & - \frac{2679}{2968} + \frac{27631 \sqrt{6}}{17808} & - \frac{191 \sqrt{110}}{9328} - \frac{25 \sqrt{165}}{4664} & - \frac{35 \sqrt{2915}}{1749} + \frac{49 \sqrt{17490}}{6996}\\- \frac{2 \sqrt{17490}}{159} + \frac{4 \sqrt{2915}}{159} & - \frac{15 \sqrt{40810}}{4664} + \frac{\sqrt{61215}}{88} & - \frac{28 \sqrt{4081}}{583} + \frac{3 \sqrt{24486}}{583} & - \frac{191 \sqrt{110}}{9328} - \frac{25 \sqrt{165}}{4664} & - \frac{10191}{4664} + \frac{3693 \sqrt{6}}{9328} & - \frac{7 \sqrt{159}}{1749} + \frac{199 \sqrt{106}}{6996}\\0 & - \frac{\sqrt{2310}}{33} + \frac{\sqrt{385}}{22} & - \frac{2 \sqrt{154}}{33} + \frac{2 \sqrt{231}}{33} & - \frac{35 \sqrt{2915}}{1749} + \frac{49 \sqrt{17490}}{6996} & - \frac{7 \sqrt{159}}{1749} + \frac{199 \sqrt{106}}{6996} & - \frac{34}{11} + \frac{47 \sqrt{6}}{11}\end{pmatrix},
\end{flalign*}
\end{minipage}}
\begin{flalign*}
    B_{(3,1)}^{odd} &= &
\end{flalign*}
\noindent\scalebox{0.9}{
\begin{minipage}{\textwidth}
\begin{flalign*}
\left(\begin{smallmatrix}-54 + 21 \sqrt{6} & 0 & 0 & 0 & 0 & 0\\0 & \frac{243}{2} - \frac{189 \sqrt{6}}{4} & - 14 \sqrt{10} + 12 \sqrt{15} & - \frac{69 \sqrt{2226}}{212} + \frac{161 \sqrt{371}}{212} & - \frac{21 \sqrt{61215}}{212} + \frac{27 \sqrt{40810}}{212} & 0\\0 & - 14 \sqrt{10} + 12 \sqrt{15} & 96 - \frac{112 \sqrt{6}}{3} & - \frac{24 \sqrt{3710}}{53} + \frac{56 \sqrt{5565}}{159} & - \frac{12 \sqrt{24486}}{53} + \frac{28 \sqrt{4081}}{53} & 0\\0 & - \frac{69 \sqrt{2226}}{212} + \frac{161 \sqrt{371}}{212} & - \frac{24 \sqrt{3710}}{53} + \frac{56 \sqrt{5565}}{159} & \frac{10455}{212} - \frac{24395 \sqrt{6}}{1272} & - \frac{147 \sqrt{165}}{212} + \frac{343 \sqrt{110}}{424} & 0\\0 & - \frac{21 \sqrt{61215}}{212} + \frac{27 \sqrt{40810}}{212} & - \frac{12 \sqrt{24486}}{53} + \frac{28 \sqrt{4081}}{53} & - \frac{147 \sqrt{165}}{212} + \frac{343 \sqrt{110}}{424} & - \frac{3141}{212} + \frac{2443 \sqrt{6}}{424} & 0\\0 & 0 & 0 & 0 & 0 & 108 - 42 \sqrt{6}\end{smallmatrix}\right), &
\end{flalign*}
\end{minipage}}

\begin{flalign*}
S_{(3,1)}^{\mathrm{odd}} &= \mathrm{diag}(6,20,20,20,20,42). &
\end{flalign*}
\begin{flalign*}
A_{(3,1)}^{\mathrm{even}} &= \left(\begin{smallmatrix}- \frac{44}{9} + \frac{16 \sqrt{6}}{27} & \frac{2 \sqrt{5}}{9} + \frac{7 \sqrt{30}}{27} & - \frac{4 \sqrt{5}}{9} + \frac{8 \sqrt{30}}{27} & - \frac{\sqrt{70}}{9} + \frac{2 \sqrt{105}}{27}\\\frac{2 \sqrt{5}}{9} + \frac{7 \sqrt{30}}{27} & - \frac{2}{3} - \frac{\sqrt{6}}{27} & - \frac{4}{9} + \frac{4 \sqrt{6}}{27} & - \frac{8 \sqrt{21}}{27} + \frac{\sqrt{14}}{3}\\- \frac{4 \sqrt{5}}{9} + \frac{8 \sqrt{30}}{27} & - \frac{4}{9} + \frac{4 \sqrt{6}}{27} & - \frac{46}{9} + \frac{101 \sqrt{6}}{27} & - \frac{22 \sqrt{21}}{27} + \frac{2 \sqrt{14}}{9}\\- \frac{\sqrt{70}}{9} + \frac{2 \sqrt{105}}{27} & - \frac{8 \sqrt{21}}{27} + \frac{\sqrt{14}}{3} & - \frac{22 \sqrt{21}}{27} + \frac{2 \sqrt{14}}{9} & - \frac{7}{3} + \frac{119 \sqrt{6}}{54}\end{smallmatrix}\right), & \\
B_{(3,1)}^{\mathrm{even}} &= \left(\begin{smallmatrix}48 - \frac{56 \sqrt{6}}{3} & - \frac{14 \sqrt{30}}{3} + 12 \sqrt{5} & 0 & 0\\- \frac{14 \sqrt{30}}{3} + 12 \sqrt{5} & -12 + \frac{14 \sqrt{6}}{3} & 0 & 0\\0 & 0 & 156 - \frac{182 \sqrt{6}}{3} & - 12 \sqrt{14} + \frac{28 \sqrt{21}}{3}\\0 & 0 & - 12 \sqrt{14} + \frac{28 \sqrt{21}}{3} & 42 - \frac{49 \sqrt{6}}{3}\end{smallmatrix}\right), & \\
S_{(3,1)}^{\mathrm{even}} &= \mathrm{diag}(12,12,30,30). &
\end{flalign*}

\subsubsection*{The $(3,2)$-component}

\begin{flalign*}
 A_{(3,2)}^{\mathrm{even}} = \left(\begin{smallmatrix}- 2 \sqrt{6} - \frac{4}{9} & \frac{2 \sqrt{5}}{9}\\\frac{2 \sqrt{5}}{9} & - \frac{5}{9} + \frac{5 \sqrt{6}}{2}\end{smallmatrix}\right), \
    B_{(3,2)}^{\mathrm{even}} &= \left(\begin{smallmatrix}-72 + 28 \sqrt{6} & 0\\0 & 90 - 35 \sqrt{6}\end{smallmatrix}\right), \
    S_{(3,2)}^{\mathrm{even}} = \mathrm{diag}(12,30). & 
\end{flalign*}
\begin{flalign*}
    B_{(3,2)}^{\mathrm{odd}} &= \begin{matrix} 0 \end{matrix}, \
    A_{(3,2)}^{\mathrm{odd}} = \begin{matrix} 0 \end{matrix}, \
    S_{(3,2)}^{\mathrm{odd}} = 20. & 
\end{flalign*}

\subsubsection*{The $(4,0)$-component.}

\begin{flalign*}
 A_{(4,0)}^{\mathrm{odd}} &= \left(\begin{smallmatrix}- \frac{291}{44} + \frac{81 \sqrt{6}}{88} & \frac{9 \sqrt{7}}{44} + \frac{21 \sqrt{42}}{88} & - \frac{3 \sqrt{42}}{11} + \frac{15 \sqrt{7}}{22} & - \frac{6 \sqrt{7}}{11} + \frac{3 \sqrt{42}}{22}\\\frac{9 \sqrt{7}}{44} + \frac{21 \sqrt{42}}{88} & - \frac{7}{4} + \frac{47 \sqrt{6}}{88} & \frac{9}{22} - \frac{3 \sqrt{6}}{11} & \frac{18}{11} - \frac{\sqrt{6}}{2}\\- \frac{3 \sqrt{42}}{11} + \frac{15 \sqrt{7}}{22} & \frac{9}{22} - \frac{3 \sqrt{6}}{11} & - \frac{84}{11} + \frac{60 \sqrt{6}}{11} & - \frac{45}{11} + \frac{3 \sqrt{6}}{11}\\- \frac{6 \sqrt{7}}{11} + \frac{3 \sqrt{42}}{22} & \frac{18}{11} - \frac{\sqrt{6}}{2} & - \frac{45}{11} + \frac{3 \sqrt{6}}{11} & -3 + \frac{57 \sqrt{6}}{22}\end{smallmatrix}\right), & \\
    B_{(4,0)}^{\mathrm{odd}} &= \left(\begin{smallmatrix}\frac{135}{2} - \frac{105 \sqrt{6}}{4} & - \frac{21 \sqrt{42}}{4} + \frac{27 \sqrt{7}}{2} & 0 & 0\\- \frac{21 \sqrt{42}}{4} + \frac{27 \sqrt{7}}{2} & \frac{9}{2} - \frac{7 \sqrt{6}}{4} & 0 & 0\\0 & 0 & 216 - 84 \sqrt{6} & 42 - 18 \sqrt{6}\\0 & 0 & 42 - 18 \sqrt{6} & 54 - 21 \sqrt{6}\end{smallmatrix}\right), & \\
    S_{(4,0)}^{\mathrm{odd}} &= \mathrm{diag}(20,20,42,42). & 
\end{flalign*}
\begin{flalign*}
    B_{(4,0)}^{\mathrm{even}} &=  \begin{matrix}108 - 42 \sqrt{6}\end{matrix}, \
    A_{(4,0)}^{\mathrm{even}} = \begin{matrix}-6 + 3 \sqrt{6}\end{matrix}, \
    S_{(4,0)}^{\mathrm{even}} = 30. & 
\end{flalign*}

\subsubsection*{The $(4,1)$-component.}

\begin{flalign*}
    A_{(4,1)}^{even} &= &
\end{flalign*}

\noindent\scalebox{0.9}{
\begin{minipage}{\textwidth}
\begin{flalign*}
    \left(\begin{smallmatrix}- \frac{4 \sqrt{6}}{3} - \frac{16}{9} & - \frac{\sqrt{110}}{11} + \frac{\sqrt{165}}{22} & - \frac{\sqrt{66}}{11} + \frac{2 \sqrt{11}}{11} & - \frac{31 \sqrt{66}}{759} + \frac{901 \sqrt{11}}{4554} & - \frac{16 \sqrt{4290}}{759} + \frac{16 \sqrt{715}}{253} & 0\\- \frac{\sqrt{110}}{11} + \frac{\sqrt{165}}{22} & - \frac{21}{4} + \frac{15 \sqrt{6}}{8} & \sqrt{15} & - \frac{37 \sqrt{10}}{184} - \frac{11 \sqrt{15}}{92} & - \frac{152 \sqrt{39}}{299} + \frac{36 \sqrt{26}}{299} & - \frac{12 \sqrt{286}}{143} + \frac{6 \sqrt{429}}{143}\\- \frac{\sqrt{66}}{11} + \frac{2 \sqrt{11}}{11} & \sqrt{15} & -4 + \frac{4 \sqrt{6}}{3} & - \frac{53}{23} - \frac{2 \sqrt{6}}{23} & - \frac{28 \sqrt{65}}{299} + \frac{158 \sqrt{390}}{897} & - \frac{4 \sqrt{715}}{143} + \frac{2 \sqrt{4290}}{143}\\- \frac{31 \sqrt{66}}{759} + \frac{901 \sqrt{11}}{4554} & - \frac{37 \sqrt{10}}{184} - \frac{11 \sqrt{15}}{92} & - \frac{53}{23} - \frac{2 \sqrt{6}}{23} & - \frac{18371}{19044} + \frac{22691 \sqrt{6}}{12696} & \frac{28 \sqrt{390}}{20631} + \frac{320 \sqrt{65}}{6877} & - \frac{60 \sqrt{4290}}{3289} + \frac{126 \sqrt{715}}{3289}\\- \frac{16 \sqrt{4290}}{759} + \frac{16 \sqrt{715}}{253} & - \frac{152 \sqrt{39}}{299} + \frac{36 \sqrt{26}}{299} & - \frac{28 \sqrt{65}}{299} + \frac{158 \sqrt{390}}{897} & \frac{28 \sqrt{390}}{20631} + \frac{320 \sqrt{65}}{6877} & - \frac{1062}{529} + \frac{2207 \sqrt{6}}{20631} & - \frac{9 \sqrt{66}}{253} + \frac{32 \sqrt{11}}{3289}\\0 & - \frac{12 \sqrt{286}}{143} + \frac{6 \sqrt{429}}{143} & - \frac{4 \sqrt{715}}{143} + \frac{2 \sqrt{4290}}{143} & - \frac{60 \sqrt{4290}}{3289} + \frac{126 \sqrt{715}}{3289} & - \frac{9 \sqrt{66}}{253} + \frac{32 \sqrt{11}}{3289} & -3 + \frac{123 \sqrt{6}}{26}\end{smallmatrix}\right), &
    \end{flalign*}
\end{minipage}
}

\newpage
\begin{flalign*}
    B_{(4,1)}^{even} &= &
\end{flalign*}
\noindent\scalebox{1}{
\begin{minipage}{\textwidth}
\begin{flalign*}
    &\left(\begin{smallmatrix}-72 + 28 \sqrt{6} & 0 & 0 & 0 & 0 & 0\\0 & \frac{243}{2} - \frac{189 \sqrt{6}}{4} & - 14 \sqrt{15} + 18 \sqrt{10} & - \frac{243 \sqrt{15}}{46} + \frac{567 \sqrt{10}}{92} & - \frac{144 \sqrt{26}}{23} + \frac{112 \sqrt{39}}{23} & 0\\0 & - 14 \sqrt{15} + 18 \sqrt{10} & 96 - \frac{112 \sqrt{6}}{3} & \frac{798}{23} - \frac{342 \sqrt{6}}{23} & - \frac{140 \sqrt{390}}{69} + \frac{120 \sqrt{65}}{23} & 0\\0 & - \frac{243 \sqrt{15}}{46} + \frac{567 \sqrt{10}}{92} & \frac{798}{23} - \frac{342 \sqrt{6}}{23} & \frac{59661}{1058} - \frac{46403 \sqrt{6}}{2116} & - \frac{672 \sqrt{65}}{529} + \frac{288 \sqrt{390}}{529} & 0\\0 & - \frac{144 \sqrt{26}}{23} + \frac{112 \sqrt{39}}{23} & - \frac{140 \sqrt{390}}{69} + \frac{120 \sqrt{65}}{23} & - \frac{672 \sqrt{65}}{529} + \frac{288 \sqrt{390}}{529} & - \frac{11580}{529} + \frac{13510 \sqrt{6}}{1587} & 0\\0 & 0 & 0 & 0 & 0 & 126 - 49 \sqrt{6}\end{smallmatrix}\right), & 
\end{flalign*}
\end{minipage}
}

\begin{flalign*}
    S_{(4,1)}^{even} &= \mathrm{diag}(12,30,30,30,30,56). &
\end{flalign*}
\begin{flalign*}
    A_{(4,1)}^{\mathrm{odd}} &= \left(\begin{smallmatrix}- \frac{54}{11} + \frac{3 \sqrt{6}}{11} & \frac{6}{11} + \frac{9 \sqrt{6}}{11} & - \frac{2 \sqrt{30}}{11} + \frac{8 \sqrt{5}}{11} & \frac{8}{11} - \frac{4 \sqrt{6}}{11}\\\frac{6}{11} + \frac{9 \sqrt{6}}{11} & - \frac{9}{11} - \frac{3 \sqrt{6}}{22} & - \frac{\sqrt{30}}{11} + \frac{2 \sqrt{5}}{11} & - \frac{16}{11} + \frac{6 \sqrt{6}}{11}\\- \frac{2 \sqrt{30}}{11} + \frac{8 \sqrt{5}}{11} & - \frac{\sqrt{30}}{11} + \frac{2 \sqrt{5}}{11} & - \frac{56}{11} + \frac{134 \sqrt{6}}{33} & - \frac{26 \sqrt{30}}{33} + \frac{4 \sqrt{5}}{11}\\\frac{8}{11} - \frac{4 \sqrt{6}}{11} & - \frac{16}{11} + \frac{6 \sqrt{6}}{11} & - \frac{26 \sqrt{30}}{33} + \frac{4 \sqrt{5}}{11} & - \frac{24}{11} + \frac{76 \sqrt{6}}{33}\end{smallmatrix}\right), & \\
        B_{(4,1)}^{\mathrm{odd}} &= \left(\begin{smallmatrix}36 - 14 \sqrt{6} & 36 - 14 \sqrt{6} & 0 & 0\\36 - 14 \sqrt{6} & -18 + 7 \sqrt{6} & 0 & 0\\0 & 0 & 168 - \frac{196 \sqrt{6}}{3} & - 24 \sqrt{5} + \frac{28 \sqrt{30}}{3}\\0 & 0 & - 24 \sqrt{5} + \frac{28 \sqrt{30}}{3} & 48 - \frac{56 \sqrt{6}}{3}\end{smallmatrix}\right), &
        \\
    S_{(4,1)}^{\mathrm{odd}} &= \mathrm{diag}(20,20,42,42). &
\end{flalign*}

\subsubsection*{The $(4,2)$-component.}

\begin{flalign*}
    A_{(4,2)}^{\mathrm{even}}= 0, \
    B_{(4,2)}^{\mathrm{even}}& = 0, \
    S_{(4,2)}^{\mathrm{even}} = 30. &
\end{flalign*}
\begin{flalign*}
 A_{(4,2)}^{\mathrm{odd}} = \left(\begin{smallmatrix}- \frac{5 \sqrt{6}}{2} - \frac{5}{11} & \frac{\sqrt{30}}{11}\\\frac{\sqrt{30}}{11} & - \frac{6}{11} + 3 \sqrt{6}\end{smallmatrix}\right), \ 
    B_{(4,2)}^{\mathrm{odd}} &= \left(\begin{smallmatrix}-90 + 35 \sqrt{6} & 0\\0 & 108 - 42 \sqrt{6}\end{smallmatrix}\right), \ 
    S_{(4,2)}^{\mathrm{odd}} = \mathrm{diag}(20,42). &
\end{flalign*}

\subsubsection*{The $(5,0)$-component.}

\begin{flalign*}
    A_{(5,0)}^{\mathrm{even}} &= \left(\begin{smallmatrix}- \frac{87}{13} + \frac{15 \sqrt{6}}{26} & \frac{3 \sqrt{6}}{13} + \frac{27}{13} & - \frac{6 \sqrt{15}}{13} + \frac{15 \sqrt{10}}{26} & - \frac{18}{13} + \frac{9 \sqrt{6}}{26}\\\frac{3 \sqrt{6}}{13} + \frac{27}{13} & - \frac{24}{13} + \frac{6 \sqrt{6}}{13} & - \frac{3 \sqrt{10}}{13} + \frac{3 \sqrt{15}}{26} & - \frac{33}{26} + \frac{9 \sqrt{6}}{13}\\- \frac{6 \sqrt{15}}{13} + \frac{15 \sqrt{10}}{26} & - \frac{3 \sqrt{10}}{13} + \frac{3 \sqrt{15}}{26} & - \frac{393}{52} + \frac{603 \sqrt{6}}{104} & - \frac{153 \sqrt{10}}{104} + \frac{9 \sqrt{15}}{52}\\- \frac{18}{13} + \frac{9 \sqrt{6}}{26} & - \frac{33}{26} + \frac{9 \sqrt{6}}{13} & - \frac{153 \sqrt{10}}{104} + \frac{9 \sqrt{15}}{52} & - \frac{151}{52} + \frac{277 \sqrt{6}}{104}\end{smallmatrix}\right), &\\
        B_{(5,0)}^{\mathrm{even}} &= \left(\begin{smallmatrix}54 - 21 \sqrt{6} & -42 + 18 \sqrt{6} & 0 & 0\\-42 + 18 \sqrt{6} & 0 & 0 & 0\\0 & 0 & \frac{459}{2} - \frac{357 \sqrt{6}}{4} & - \frac{27 \sqrt{15}}{2} + \frac{63 \sqrt{10}}{4}\\0 & 0 & - \frac{27 \sqrt{15}}{2} + \frac{63 \sqrt{10}}{4} & \frac{117}{2} - \frac{91 \sqrt{6}}{4}\end{smallmatrix}\right), &\\
    S_{(5,0)}^{\mathrm{even}} &= \mathrm{diag}(30,30,56,56). & 
\end{flalign*}
\begin{flalign*}
 A_{(5,0)}^{\mathrm{odd}} = \begin{matrix}-6 + 3 \sqrt{6}\end{matrix}, \
    B_{(5,0)}^{\mathrm{odd}} &=  \begin{matrix}108 - 42 \sqrt{6}\end{matrix}, \
    S_{(5,0)}^{\mathrm{odd}} = 42. & 
\end{flalign*}

\subsubsection*{The $(5,1)$-component.}

\begin{flalign*}
    A_{(5,1)}^{\mathrm{even}} &= \left(\begin{smallmatrix}- \frac{64}{13} - \frac{2 \sqrt{6}}{39} & \frac{2 \sqrt{14}}{13} + \frac{22 \sqrt{21}}{39} & - \frac{8 \sqrt{7}}{13} + \frac{2 \sqrt{42}}{13} & - \frac{3 \sqrt{14}}{13} + \frac{2 \sqrt{21}}{13}\\\frac{2 \sqrt{14}}{13} + \frac{22 \sqrt{21}}{39} & - \frac{12}{13} - \frac{10 \sqrt{6}}{39} & - \frac{4 \sqrt{2}}{13} + \frac{4 \sqrt{3}}{13} & \frac{18}{13} - \frac{8 \sqrt{6}}{13}\\- \frac{8 \sqrt{7}}{13} + \frac{2 \sqrt{42}}{13} & - \frac{4 \sqrt{2}}{13} + \frac{4 \sqrt{3}}{13} & - \frac{66}{13} + \frac{57 \sqrt{6}}{13} & - \frac{6 \sqrt{3}}{13} + \frac{45 \sqrt{2}}{13}\\- \frac{3 \sqrt{14}}{13} + \frac{2 \sqrt{21}}{13} & \frac{18}{13} - \frac{8 \sqrt{6}}{13} & - \frac{6 \sqrt{3}}{13} + \frac{45 \sqrt{2}}{13} & - \frac{27}{13} + \frac{63 \sqrt{6}}{26}\end{smallmatrix}\right), &\\
    B_{(5,1)}^{\mathrm{even}} &= \left(\begin{smallmatrix}24 - \frac{28 \sqrt{6}}{3} & - \frac{28 \sqrt{21}}{3} + 12 \sqrt{14} & 0 & 0\\- \frac{28 \sqrt{21}}{3} + 12 \sqrt{14} & -24 + \frac{28 \sqrt{6}}{3} & 0 & 0\\0 & 0 & 180 - 70 \sqrt{6} & - 42 \sqrt{2} + 36 \sqrt{3}\\0 & 0 & - 42 \sqrt{2} + 36 \sqrt{3} & 54 - 21 \sqrt{6}\end{smallmatrix}\right), &\\
    S_{(5,1)}^{\mathrm{even}} &= \mathrm{diag}(30,30,56,56). &
\end{flalign*}
\begin{flalign*}
   A_{(5,1)}^{\mathrm{odd}} &= &
\end{flalign*}
\noindent\scalebox{0.62}{
\begin{minipage}{\textwidth}
\begin{flalign*}
    \begin{pmatrix}- \frac{39 \sqrt{6}}{22} - \frac{21}{11} & - \frac{12 \sqrt{143}}{143} + \frac{3 \sqrt{858}}{143} & - \frac{2 \sqrt{3003}}{143} + \frac{2 \sqrt{2002}}{143} & - \frac{148 \sqrt{14729}}{14729} + \frac{3 \sqrt{88374}}{1133} & - \frac{48 \sqrt{140595}}{14729} + \frac{63 \sqrt{93730}}{14729} & 0\\- \frac{12 \sqrt{143}}{143} + \frac{3 \sqrt{858}}{143} & - \frac{21}{4} + \frac{15 \sqrt{6}}{8} & \sqrt{21} & - \frac{9 \sqrt{103}}{412} + \frac{115 \sqrt{618}}{824} & - \frac{12 \sqrt{118965}}{5665} + \frac{24 \sqrt{79310}}{5665} & - \frac{3 \sqrt{910}}{65} + \frac{3 \sqrt{1365}}{130}\\- \frac{2 \sqrt{3003}}{143} + \frac{2 \sqrt{2002}}{143} & \sqrt{21} & -4 + \frac{4 \sqrt{6}}{3} & - \frac{11 \sqrt{2163}}{309} + \frac{2 \sqrt{1442}}{103} & - \frac{342 \sqrt{5665}}{5665} + \frac{6 \sqrt{33990}}{5665} & - \frac{6 \sqrt{65}}{65} + \frac{3 \sqrt{390}}{65}\\- \frac{148 \sqrt{14729}}{14729} + \frac{3 \sqrt{88374}}{1133} & - \frac{9 \sqrt{103}}{412} + \frac{115 \sqrt{618}}{824} & - \frac{11 \sqrt{2163}}{309} + \frac{2 \sqrt{1442}}{103} & - \frac{669}{412} + \frac{1765 \sqrt{6}}{2472} & - \frac{108 \sqrt{770}}{1133} - \frac{68 \sqrt{1155}}{5665} & - \frac{\sqrt{140595}}{1030} + \frac{3 \sqrt{93730}}{1339}\\- \frac{48 \sqrt{140595}}{14729} + \frac{63 \sqrt{93730}}{14729} & - \frac{12 \sqrt{118965}}{5665} + \frac{24 \sqrt{79310}}{5665} & - \frac{342 \sqrt{5665}}{5665} + \frac{6 \sqrt{33990}}{5665} & - \frac{108 \sqrt{770}}{1133} - \frac{68 \sqrt{1155}}{5665} & - \frac{21818}{16995} + \frac{6517 \sqrt{6}}{5665} & - \frac{24 \sqrt{88374}}{6695} + \frac{176 \sqrt{14729}}{20085}\\0 & - \frac{3 \sqrt{910}}{65} + \frac{3 \sqrt{1365}}{130} & - \frac{6 \sqrt{65}}{65} + \frac{3 \sqrt{390}}{65} & - \frac{\sqrt{140595}}{1030} + \frac{3 \sqrt{93730}}{1339} & - \frac{24 \sqrt{88374}}{6695} + \frac{176 \sqrt{14729}}{20085} & - \frac{44}{15} + \frac{26 \sqrt{6}}{5}\end{pmatrix},
\end{flalign*}
\end{minipage}}
\begin{flalign*}
    B_{(5,1)}^{\mathrm{odd}} &= &
\end{flalign*}
\noindent\scalebox{1}{
\begin{minipage}{\textwidth}
\begin{flalign*}
    &\left( \begin{smallmatrix}-90 + 35 \sqrt{6} & 0 & 0 & 0 & 0 & 0\\0 & \frac{243}{2} - \frac{189 \sqrt{6}}{4} & - 14 \sqrt{21} + 18 \sqrt{14} & - \frac{7 \sqrt{618}}{4} + \frac{9 \sqrt{103}}{2} & 0 & 0\\0 & - 14 \sqrt{21} + 18 \sqrt{14} & 96 - \frac{112 \sqrt{6}}{3} & - \frac{42 \sqrt{1442}}{103} + \frac{98 \sqrt{2163}}{309} & - \frac{36 \sqrt{33990}}{103} + \frac{84 \sqrt{5665}}{103} & 0\\0 & - \frac{7 \sqrt{618}}{4} + \frac{9 \sqrt{103}}{2} & - \frac{42 \sqrt{1442}}{103} + \frac{98 \sqrt{2163}}{309} & \frac{2283}{206} - \frac{5327 \sqrt{6}}{1236} & - \frac{144 \sqrt{1155}}{103} + \frac{168 \sqrt{770}}{103} & 0\\0 & 0 & - \frac{36 \sqrt{33990}}{103} + \frac{84 \sqrt{5665}}{103} & - \frac{144 \sqrt{1155}}{103} + \frac{168 \sqrt{770}}{103} & \frac{2412}{103} - \frac{938 \sqrt{6}}{103} & 0\\0 & 0 & 0 & 0 & 0 & 144 - 56 \sqrt{6}\end{smallmatrix}\right), &
\end{flalign*}
\end{minipage}}

\begin{flalign*}
    S_{(5,1)}^{\mathrm{odd}} &= \mathrm{diag}(20,42,42,42,42,72). &
\end{flalign*}

\bibliography{references}

@book {Sepanski2007,
    AUTHOR = {Sepanski, Mark R.},
     TITLE = {Compact {L}ie groups},
    SERIES = {Graduate Texts in Mathematics},
    VOLUME = {235},
 PUBLISHER = {Springer, New York},
      YEAR = {2007},
     PAGES = {xiv+198},
      ISBN = {978-0-387-30263-8; 0-387-30263-8},
   MRCLASS = {22Exx},
  MRNUMBER = {2279709},
MRREVIEWER = {David\ A.\ Renard},
       DOI = {10.1007/978-0-387-49158-5},
       URL = {https://doi.org/10.1007/978-0-387-49158-5},
}

@misc{KimLee24,
      title={Type-{I} Blowup Solutions for {Yang-Mills} Flow}, 
      author={Jaehwan Kim and Sanghoon Lee},
      year={2024},
      eprint={2411.19293},
      archivePrefix={arXiv},
      primaryClass={math.DG},
      url={https://arxiv.org/abs/2411.19293}, 
}

@misc{Fee16,
      title={Global existence and convergence of solutions to gradient systems and applications to {Yang-Mills} gradient flow}, 
      author={Paul M. N. Feehan},
      year={2016},
      eprint={1409.1525},
      archivePrefix={arXiv},
      primaryClass={math.DG},
      url={https://arxiv.org/abs/1409.1525}, 
}

@article{Rad92,
url = {https://doi.org/10.1515/crll.1992.431.123},
title = {On the {Yang-Mills} heat equation in two and three dimensions.},
author = {Johan Rade},
pages = {123--164},
volume = {1992},
number = {431},
journal = {Journal für die reine und angewandte Mathematik},
doi = {doi:10.1515/crll.1992.431.123},
year = {1992},
lastchecked = {2026-10-06}
}

@article {Str94,
    AUTHOR = {Struwe, Michael},
     TITLE = {The {Y}ang-{M}ills flow in four dimensions},
   JOURNAL = {Calc. Var. Partial Differential Equations},
  FJOURNAL = {Calculus of Variations and Partial Differential Equations},
    VOLUME = {2},
      YEAR = {1994},
    NUMBER = {2},
     PAGES = {123--150},
      ISSN = {0944-2669,1432-0835},
   MRCLASS = {58E15 (35K22 53C07 58G11)},
  MRNUMBER = {1385523},
MRREVIEWER = {J\"urgen\ Eichhorn},
       DOI = {10.1007/BF01191339},
       URL = {https://doi.org/10.1007/BF01191339},
}

@article{SchStrTah98,
 ISSN = {00029327, 10806377},
 URL = {http://www.jstor.org/stable/25098587},
 author = {Andreas E. Schlatter and Michael Struwe and A. Shadi Tahvildar-Zadeh},
 journal = {American Journal of Mathematics},
 number = {1},
 pages = {117--128},
 publisher = {Johns Hopkins University Press},
 title = {Global Existence of the Equivariant {Yang-Mills} Heat Flow in Four Space Dimensions},
 urldate = {2026-10-19},
 volume = {120},
 year = {1998}
}

@article {Wal19,
    AUTHOR = {Waldron, Alex},
     TITLE = {Long-time existence for {Y}ang-{M}ills flow},
   JOURNAL = {Invent. Math.},
  FJOURNAL = {Inventiones Mathematicae},
    VOLUME = {217},
      YEAR = {2019},
    NUMBER = {3},
     PAGES = {1069--1147},
      ISSN = {0020-9910,1432-1297},
   MRCLASS = {53E99 (35R01 58E15)},
  MRNUMBER = {3989258},
MRREVIEWER = {Yu\ Zheng},
       DOI = {10.1007/s00222-019-00877-2},
       URL = {https://doi.org/10.1007/s00222-019-00877-2},
}

@article {OhTat22,
    AUTHOR = {Oh, Sung-Jin and Tataru, Daniel},
     TITLE = {The {Y}ang-{M}ills heat flow and the caloric gauge},
   JOURNAL = {Ast\'erisque},
  FJOURNAL = {Ast\'erisque},
    NUMBER = {436},
      YEAR = {2022},
     PAGES = {viii+128},
      ISSN = {0303-1179,2492-5926},
      ISBN = {978-2-85629-961-6},
   MRCLASS = {35K55 (70S15)},
  MRNUMBER = {4518477},
       DOI = {10.24033/ast.1179},
       URL = {https://doi.org/10.24033/ast.1179},
}

@misc{SirWerZheZho26,
      title={Global oscillatory solutions for the {Yang-Mills} heat flow}, 
      author={Yannick Sire and Juncheng Wei and Youquan Zheng and Yifu Zhou},
      year={2026},
      eprint={2601.21017},
      archivePrefix={arXiv},
      primaryClass={math.AP},
      url={https://arxiv.org/abs/2601.21017}, 
}

@article{Nai94,
author = {Hisashi Naito},
title = {{Finite time blowing-up for the {Y}ang-{M}ills gradient flow in higher dimensions}},
volume = {23},
journal = {Hokkaido Mathematical Journal},
number = {3},
publisher = {Hokkaido University, Department of Mathematics},
pages = {451 -- 464},
year = {1994},
doi = {10.14492/hokmj/1381413099},
URL = {https://doi.org/10.14492/hokmj/1381413099}
}

@article{Gro01,
author={Grotowski, Joseph F.},
title={Finite time blow-up for the {Yang-Mills} heat flow in higher dimensions},
journal={Mathematische Zeitschrift},
year={2001},
month={Jun},
day={01},
volume={237},
number={2},
pages={321-333},
issn={1432-1823},
doi={10.1007/PL00004871},
url={https://doi.org/10.1007/PL00004871}
}

@article{Gas02,
author={Gastel, Andreas},
title={Singularities of first kind in the harmonic map and {Yang-Mills} heat flows},
journal={Mathematische Zeitschrift},
year={2002},
month={Feb},
day={01},
volume={242},
number={1},
pages={47-62},
issn={1432-1823},
doi={10.1007/s002090100306},
url={https://doi.org/10.1007/s002090100306}
}

@ARTICLE{BizWas15,
  author={Bizoń, Piotr and Wasserman, Arthur},
  journal={International Mathematics Research Notices}, 
  title={Nonexistence of Shrinkers for the Harmonic Map Flow in Higher Dimensions}, 
  year={2015},
  volume={2015},
  number={17},
  pages={7757-7762},
  doi={10.1093/imrn/rnu176}}

@article{Wan25,
author={Wang, Guanxiang
and Zhang, Chuanjing},
title={The Finite Time Blow-up of the {Yang-Mills} Flow},
journal={Communications in Mathematics and Statistics},
year={2025},
month={Jan},
day={15},
issn={2194-671X},
doi={10.1007/s40304-024-00405-6},
url={https://doi.org/10.1007/s40304-024-00405-6}
}

@article{Wei04,
author={Weinkove, Ben},
title={Singularity formation in the {Yang-Mills Flow}},
   JOURNAL = {Calc. Var. Partial Differential Equations},
  FJOURNAL = {Calculus of Variations and Partial Differential Equations},
year={2004},
month={Jan},
day={01},
volume={19},
number={2},
pages={211-220},
issn={1432-0835},
doi={10.1007/s00526-003-0217-x},
url={https://doi.org/10.1007/s00526-003-0217-x}
}

@article{DonSch19,
author = {Roland Donninger and Birgit Sch{\"o}rkhuber},
title = {{Stable blowup for the supercritical Yang–Mills heat flow}},
volume = {113},
journal = {Journal of Differential Geometry},
number = {1},
publisher = {Lehigh University},
pages = {55 -- 94},
year = {2019},
doi = {10.4310/jdg/1567216954},
URL = {https://doi.org/10.4310/jdg/1567216954}
}

@article{GloSch20,
author = {Irfan Glogić and Birgit Schörkhuber},
title = {Nonlinear stability of homothetically shrinking {Yang-Mills }solitons in the equivariant case},
journal = {Communications in Partial Differential Equations},
volume = {45},
number = {8},
pages = {887--912},
year = {2020},
publisher = {Taylor \& Francis},
doi = {10.1080/03605302.2020.1743308},
URL = { 
    
        https://doi.org/10.1080/03605302.2020.1743308
},
eprint = { 
    
        https://doi.org/10.1080/03605302.2020.1743308    
}
}

@article{CheZha15,
 ISSN = {00029947},
 URL = {http://www.jstor.org/stable/24513379},
 author = {Zhengxiang Chen and Yongbing Zhang},
 journal = {Transactions of the American Mathematical Society},
 number = {7},
 pages = {5015--5041},
 publisher = {American Mathematical Society},
 title = {STABILITIES OF HOMOTHETICALLY SHRINKING {Y}ANG-{M}ILLS SOLITONS},
 urldate = {2026-10-19},
 volume = {367},
 year = {2015}
}

@article{KelStr16,
author = {Kelleher, Casey and Streets, Jeffrey},
title = {Entropy, stability, and {Yang–Mills flow}},
journal = {Communications in Contemporary Mathematics},
volume = {18},
number = {02},
pages = {1550032},
year = {2016},
doi = {10.1142/S0219199715500327},
URL = { 
    
        https://doi.org/10.1142/S0219199715500327 
},
eprint = { 
    
        https://doi.org/10.1142/S0219199715500327  
}
}

@article{NafOzu25,
      title={Linear stability of the blowdown {R}icci shrinker in 4D}, 
      author={Keaton Naff and Tristan Ozuch},
      year={2025},
      eprint={2509.05470},
      archivePrefix={arXiv},
      primaryClass={math.DG},
      url={https://arxiv.org/abs/2509.05470}, 
}

@article{Max14,
url = {https://doi.org/10.1515/crelle-2012-0080},
title = {On the blow-up of four-dimensional {Ricci} flow singularities},
author = {Davi Máximo},
pages = {153--171},
volume = {2014},
number = {692},
journal = {Journal für die reine und angewandte Mathematik (Crelles Journal)},
doi = {doi:10.1515/crelle-2012-0080},
year = {2014},
lastchecked = {2026-10-19}
}

@article {BinJia16,
    AUTHOR = {Guo, Bin and Song, Jian},
     TITLE = {On {F}eldman-{I}lmanen-{K}nopf's conjecture for the blow-up
              behavior of the {K}\"ahler {R}icci flow},
   JOURNAL = {Math. Res. Lett.},
  FJOURNAL = {Mathematical Research Letters},
    VOLUME = {23},
      YEAR = {2016},
    NUMBER = {6},
     PAGES = {1681--1719},
      ISSN = {1073-2780,1945-001X},
   MRCLASS = {53C44 (53C55)},
  MRNUMBER = {3621103},
       DOI = {10.4310/MRL.2016.v23.n6.a6},
       URL = {https://doi.org/10.4310/MRL.2016.v23.n6.a6},
}

@article{BieBiz11,
doi = {10.1088/0951-7715/24/8/005},
url = {https://doi.org/10.1088/0951-7715/24/8/005},
year = {2011},
month = {Jun},
publisher = {},
volume = {24},
number = {8},
pages = {2211},
author = {Biernat, Paweł and Bizoń, Piotr},
title = {Shrinkers, expanders, and the unique continuation beyond generic blowup in the heat flow for harmonic maps between spheres},
journal = {Nonlinearity}
}

@article{BieDon18,
doi = {10.1088/1361-6544/aabe4c},
url = {https://doi.org/10.1088/1361-6544/aabe4c},
year = {2018},
month = {Jun},
publisher = {IOP Publishing},
volume = {31},
number = {8},
pages = {3543},
author = {Biernat, Pawel and Donninger, Roland},
title = {Construction of a spectrally stable self-similar blowup solution to the supercritical corotational harmonic map heat flow},
journal = {Nonlinearity}
}

@article{BieDonSch17,
author={Biernat, Pawe{\l}
and Donninger, Roland
and Sch{\"o}rkhuber, Birgit},
title={Stable self-similar blowup in the supercritical heat flow of harmonic maps},
   JOURNAL = {Calc. Var. Partial Differential Equations},
  FJOURNAL = {Calculus of Variations and Partial Differential Equations},
year={2017},
month={Nov},
day={04},
volume={56},
number={6},
pages={171},
issn={1432-0835},
doi={10.1007/s00526-017-1256-z},
url={https://doi.org/10.1007/s00526-017-1256-z}
}

@article{GloKisSch24,
author={Glogi{\'{c}}, Irfan
and Kistner, Sarah
and Sch{\"o}rkhuber, Birgit},
title={Existence and stability of shrinkers for the harmonic map heat flow in higher dimensions},
   JOURNAL = {Calc. Var. Partial Differential Equations},
  FJOURNAL = {Calculus of Variations and Partial Differential Equations},
year={2024},
month={Apr},
day={10},
volume={63},
number={4},
pages={96},
issn={1432-0835},
doi={10.1007/s00526-024-02707-7},
url={https://doi.org/10.1007/s00526-024-02707-7}
}

@article{AngKisSch26,
      title={Existence of a stable shrinker for the corotational harmonic map heat flow in higher space dimensions}, 
      author={Johannes Angerer and Sarah Kistner and Birgit Schörkhuber},
      year={2026},
      eprint={2607.27072},
      archivePrefix={arXiv},
      primaryClass={math.AP},
      url={https://arxiv.org/abs/2607.27072}, 
}

@article{Fan19,
author={Fan, Huijun},
title={Existence of the self-similar solutions in the heat flow of harmonic maps},
journal={Science in China Series A: Mathematics},
year={1999},
month={Feb},
day={01},
volume={42},
number={2},
pages={113-132},
issn={1862-2763},
doi={10.1007/BF02876563},
url={https://doi.org/10.1007/BF02876563}
}

@article {BenDuoGho25,
    AUTHOR = {Bensouilah, A. and Duong, G. K. and Ghoul, T. E.},
     TITLE = {Non-self similar blowup solutions to the higher dimensional
              {Y}ang {M}ills heat flows},
   JOURNAL = {J. Differential Equations},
  FJOURNAL = {Journal of Differential Equations},
    VOLUME = {427},
      YEAR = {2025},
     PAGES = {26--142},
      ISSN = {0022-0396,1090-2732},
   MRCLASS = {35K55 (35B44 35K57)},
  MRNUMBER = {4858612},
       DOI = {10.1016/j.jde.2025.01.039},
       URL = {https://doi.org/10.1016/j.jde.2025.01.039},
}

@article {KelStr18,
    AUTHOR = {Kelleher, Casey and Streets, Jeffrey},
     TITLE = {Singularity formation of the {Y}ang-{M}ills flow},
   JOURNAL = {Ann. Inst. H. Poincar\'e{} C Anal. Non Lin\'eaire},
  FJOURNAL = {Annales de l'Institut Henri Poincar\'e{} C. Analyse Non
              Lin\'eaire},
    VOLUME = {35},
      YEAR = {2018},
    NUMBER = {6},
     PAGES = {1655--1686},
      ISSN = {0294-1449,1873-1430},
   MRCLASS = {53C44 (35Q40 81T13)},
  MRNUMBER = {3846240},
MRREVIEWER = {Dmitri\u i\ Vladimir\ Alekseevsky},
       DOI = {10.1016/j.anihpc.2018.01.006},
       URL = {https://doi.org/10.1016/j.anihpc.2018.01.006},
}

@book {ReeSim75,
    AUTHOR = {Reed, Michael and Simon, Barry},
     TITLE = {Methods of modern mathematical physics. {II}. {F}ourier
              analysis, self-adjointness},
 PUBLISHER = {Academic Press [Harcourt Brace Jovanovich, Publishers], New
              York-London},
      YEAR = {1975},
     PAGES = {xv+361},
   MRCLASS = {47-02 (81.47)},
  MRNUMBER = {493420},
MRREVIEWER = {P.\ R.\ Chernoff},
}

@book {ReeSim78,
    AUTHOR = {Reed, Michael and Simon, Barry},
     TITLE = {Methods of modern mathematical physics. {IV}. {A}nalysis of
              operators},
 PUBLISHER = {Academic Press [Harcourt Brace Jovanovich, Publishers], New
              York-London},
      YEAR = {1978},
     PAGES = {xv+396},
      ISBN = {0-12-585004-2},
   MRCLASS = {47-02 (81.47)},
  MRNUMBER = {493421},
MRREVIEWER = {P.\ R.\ Chernoff},
}

@article {LieART,
    AUTHOR = {Feger, Robert and Kephart, Thomas W. and Saskowski, Robert J.},
     TITLE = {Lie{ART} 2.0---a {M}athematica application for {L}ie algebras
              and representation theory},
   JOURNAL = {Comput. Phys. Commun.},
  FJOURNAL = {Computer Physics Communications},
    VOLUME = {257},
      YEAR = {2020},
     PAGES = {107490, 24},
      ISSN = {0010-4655,1879-2944},
   MRCLASS = {17-08 (17Bxx 22-08)},
  MRNUMBER = {4131870},
       DOI = {10.1016/j.cpc.2020.107490},
       URL = {https://doi.org/10.1016/j.cpc.2020.107490},
}

@article {Folland2001,
    AUTHOR = {Folland, Gerald B.},
     TITLE = {How to integrate a polynomial over a sphere},
   JOURNAL = {Amer. Math. Monthly},
  FJOURNAL = {American Mathematical Monthly},
    VOLUME = {108},
      YEAR = {2001},
    NUMBER = {5},
     PAGES = {446--448},
      ISSN = {0002-9890,1930-0972},
   MRCLASS = {26A06},
  MRNUMBER = {1837866},
       DOI = {10.2307/2695802},
       URL = {https://doi.org/10.2307/2695802},
}

@misc{WeiKocDon25,
      title={Mode stability of blow-up for wave maps in the absence of symmetry}, 
      author={Max Weissenbacher and Herbert Koch and Roland Donninger},
      year={2025},
      eprint={2503.02632},
      archivePrefix={arXiv},
      primaryClass={math.AP},
      url={https://arxiv.org/abs/2503.02632}, 
}

@article{DonMos26,
 author = {Donninger, Roland and Moscatelli, Frederick},
 title = {Mode stability of self-similar wave maps without symmetry in higher dimensions},
 fjournal = {Transactions of the American Mathematical Society. Series B},
 journal = {Trans. Am. Math. Soc., Ser. B},
 issn = {2330-0000},
 volume = {13},
 pages = {441--476},
 year = {2026},
 language = {English},
 doi = {10.1090/btran/258},
 zbMATH = {8234715}
}
\bibliographystyle{plain}
\end{document}